# An Alternative Finite Difference WENO-like Scheme with Physical Constraint Preservation for Divergence-Preserving Hyperbolic Systems

By

Dinshaw S. Balsara[1,2], Deepak Bhoriya[1] and Chi-Wang Shu[3]

[1]Physics Department, [2]ACMS Department, University of Notre Dame

[3]Division of Applied Mathematics, Brown University

**Abstract**

Alternative finite difference Weighted Essentially Non-Oscillatory (AFD-WENO) schemes allow us to very efficiently update hyperbolic systems even in complex geometries. Recent innovations in AFD-WENO methods allow us to treat hyperbolic system with non-conservative products almost as efficiently as conservation laws. However, some PDE systems, like computational electrodynamics (CED) and magnetohydrodynamics (MHD) and relativistic magnetohydrodynamics (RMHD), have involution constraints that require divergence-free or divergence-preserving evolution of vector fields. In such situations, a Yee-style collocation of variables proves indispensable; and that collocation is retained in this work. In previous works, only higher order finite volume discretization of such involution constrained systems was possible. In this work, we show that substantially more efficient AFD-WENO methods have been extended to encompass divergence-preserving hyperbolic PDEs.

Our method retains the Yee-style collocation of normal components of the divergence-free/preserving vector field. However, the variables that require zone-centered evolution are evolved with AFD-WENO methods. Since those variables make up the bulk of the primal variables for the PDE of interest, this results in a substantial savings in computational complexity. Even the volumetric reconstruction of the divergence-free/preserving vector field is bypassed. Instead, we realize that any divergence-preserving update of a vector field must have a general form. This general form looks closely like the familiar induction equation that is well-known in CED or MHD. We exploit the generality of that form to extract the edge-centered variables that are needed in the update of the facially-averaged vector field components. The two-dimensional Riemann solver is used to provide us multidimensionally stabilized versions of these update terms.

The generality of our approach is demonstrated by the fact that problems in CED, MHD and RMHD can all be solved by the same general AFD-WENO algorithm that is presented here. Spatial accuracies up to ninth order of accuracy are demonstrated. Several stringent test problems from CED, MHD and RMHD are shown. We also show that the algorithm takes well to the physical constraint preserving (PCP) formulation of AFD-WENO schemes that was presented by the authors. The efficient and time-explicit PCP strategy for divergence-preserving PDEs that we have presented here extends the applicability of the our method to very stringent MHD and RMHD problems.

# I) Introduction

Novel applications in science and engineering have called for higher order schemes for the simulation of hyperbolic PDEs. Most of these PDEs are strongly non-linear which implies that it is beneficial to pick solution methods that are generally drawn from the family of higher order Godunov schemes. Within that general family of Godunov schemes, the two most popular strategies for achieving higher order accuracy consist of Weighted Essentially Non-Oscillatory (WENO) schemes and Discontinuous Galerkin (DG) schemes. Essentially Non-Oscillatory methods were developed in finite volume form by Harten *et al*. [57]. DG schemes were also subsequently developed (Cockburn and Shu [47], [50], Cockburn, Lin and Shu [48], Cockburn, Hou and Shu [49]). Early work on this topic focused on conservation laws. Strong Stability Preserving Runge-Kutta (SSP-RK) methods for high order time-integration (Shu and Osher [78], Spiteri and Ruuth [81], [82]) were developed to match the high order spatial accuracy. Yet it quickly became apparent that practitioners wanted methods that could handle stiff source terms (Pareschi and Russo [72], Kupka *et al*. [67]) and non-conservative products (Baer and Nunziato [6], Andrianov and Warnecke [4], Castro *et al*. [45]). Newer classes of hyperbolic PDEs, such as magnetohydrodynamics (MHD) and computational electrodynamics (CED), came with their own involution constraints. For example, both MHD and CED require the ability to evolve vector fields in a divergence-free or divergence-preserving fashion. The divergence constraint is more of a geometrical constraint, with the result that a collocation of primal variables and update variables that is consistent with a Yee [88] type mesh becomes essential. Fig. 1 shows how the induction equation, that is common to MHD and CED, is discretized. The components of the magnetic fields are collocated, in facially-averaged form at the faces of the mesh; as in Fig. 1. These are updated using edge-averaged electric field components at the edges of the mesh; as in Fig. 1. The induction equation for MHD is generically written as

$$\frac{\partial \mathbf{B}}{\partial t} + \nabla \times \mathbf{E} = 0 \qquad (1)$$

with $\mathbf{E} = -\mathbf{v} \times \mathbf{B}$, where "$\mathbf{v}$" is the fluid velocity. In Fig. 1 we see that the components of the magnetic field "$\mathbf{B}$" (which are the primal variables of the scheme) are collocated at the faces of the mesh and the discrete divergence constraint is only preserved if the electric fields "$\mathbf{E}$" (which provide the update) are collocated at the edges of the mesh. As a result, new styles of discretization and multidimensional upwinding had to be invented. This paper focuses on hyperbolic PDE systems that may indeed have a dominant zone-centered collocation of variables; however, the PDE system includes one or more vector fields that look like the update of the induction equation. Past efforts have focused on individual PDEs. The *overarching goal* of this paper is to discern the common structure of all PDE systems that have vector fields that have an update that looks like the induction equation, and then to develop common methods for the numerical evolution of such systems.

MHD was indeed the first involution-constrained system that was extended to include higher order Godunov methodology while preserving the divergence-free evolution of magnetic fields. Early work on divergence-free Godunov schemes for MHD was restricted to second order (Dai and Woodward [52], Ryu et al. [74], Balsara and Spicer [8]). The need to do adaptive mesh refinement (AMR) at second order forced one to pay attention to the TVD-based divergence-free reconstruction of vector fields (Balsara [10]). Finite volume WENO-based, divergence-free reconstruction, as well as divergence-free AMR, has been extended to higher orders (Balsara [12], [13], Balsara, Samantaray and Subramanian [33], Balsara and Sarris [34]). This vector field reconstruction ensures that if the modes of variation within each face of Fig. 1 are known for the magnetic field components then the entire divergence-free vector field can be reconstructed at all locations within the volume shown in Fig. 1. However, note too from Fig. 1 that the edge-centered electric fields are also needed. In keeping with the realization that upwinding provides stabilization, we realize that the edge-collocated electric fields shown in Fig. 1 will have to be upwinded in their two transverse directions. This insight led to the development of multidimensional Riemann solvers (Balsara [17], [18], [22], [23], Balsara *et al*. [21], Balsara and Dumbser [24], Balsara and Nkonga [27]). Higher order divergence-free reconstruction along with the multidimensional upwinding from the multidimensional Riemann solver opened the door to globally divergence-free higher order WENO-based MHD schemes (Balsara *et al*. [16], [19], [20]). Finite volume WENO-based schemes for CED that used the same two ideas were also invented (Balsara *et al*. [29], [30]). Balsara and Käppeli [28] formulated globally divergence-free DG schemes for MHD; and working DG schemes for MHD, which used WENO methods for stabilization, were presented in Balsara *et al*. [32]. Likewise, Balsara and Käppeli [31] formulated globally divergence-preserving DG schemes for CED; and working DG schemes for CED with up to fifth order of accuracy were presented in Hazra *et al*. [58].

All of the advances in the above paragraphs were reported within the context of finite volume formulations. However, ever since the seminal papers by Shu and Osher [78], [79], we have been aware that higher order finite difference WENO schemes are substantially more efficient compared to finite volume WENO schemes. The first finite difference WENO (FD-WENO) schemes for conservation laws that were suitable for production codes emerged in Jiang and Shu [63] and Balsara and Shu [9]. Modern PDEs also include non-conservative products and stiff source terms. In Balsara *et al*. [35] a way was found for simulating such PDEs with non-conservative products and stiff source terms using FD-WENO schemes. In Balsara *et al*. [36], [36] we developed Alternative FD-WENO (AFD-WENO) methods that can accommodate non-conservative products and are suitable for production codes. The AFD-WENO schemes offer the following advantages:- First, any Riemann solver can be used, which can be useful for several applications where Riemann solvers with special capabilities are desired. Second, the free stream condition can be respected on curvilinear meshes (Jiang, Shu and Zhang [64], [65]). Third, variables that are in flux conservative form can be treated in a manner that respects conservation, while also accommodating non-conservative products. Fourth, stiff source terms can be included quite easily owing to the pointwise nature of the scheme. Fifth, well-balancing can be easily added to the scheme (Xu and Shu [87]). However, a major deficiency up to now, has been the fact that AFD-WENO methods are not available for PDEs that have a globally divergence-free or divergence-preserving involution constraint. The *first important goal* of this paper is to develop schemes that use many of the AFD-WENO principles while including the divergence constraint.

All PDE systems that have a divergence constraint have a structure like the induction equation. Therefore, they will necessarily have to retain the Yee-mesh collocation philosophy shown in Fig. 1. Because Stokes law has to be applied to the induction equation, it necessarily means that the facial magnetic fields in Fig. 1 are face-averaged and the edge-centered electric fields are edge-averaged in Fig. 1. As a result, we realize that the constrained vector fields will have to retain a finite volume style of update. However, we realize that the induction equation forms a very small part of much larger PDE systems. The MHD equations (Alfvén [3], Jeffrey and Taniuti [62], Roe and Balsara [73]), or the Chew Goldberger and Low (CGL) equations (Chew Goldberger and Low [46], Bhoriya *et al*. [40], Singh et al. [80]), or the equations of Relativistic MHD (RMHD) (Anile [5], Balsara [11], del Zanna *et al*. [89], Balsara and Kim [26]) are cases in point. In all those PDEs, we have a much larger hydrodynamical set of equations that are coupled via magnetic fields to one single induction equation for the evolution of the magnetic vector field. As a result, it would still be advantageous to treat the rest of the variables with the highly efficient and versatile AFD-WENO formulation while retaining a concept of area-averaged updates for the facial components of the magnetic vector field. This is the compromise position that we will adopt in this paper.

Such compromise positions, which combine some elements of finite difference methods with some elements of finite volume methods in order to obtain greater computational efficiency, have been adopted before. For example, see the WENO-based work of Buchmüller, Dreher and Helzel [43] who drew their inspiration from McCorquodale and Colella [69], Shi, Hu and Shu [77] and Zhang, Zhang and Shu [90]. Recently, Donnert *et al.* [53] used a fifth order FD-WENO formulation for zone-centered variables along with a second order formulation for the magnetic fields to at least get a more accurate evolution of the fluid variables; even though the evolution of the magnetic fields was only second order accurate. Seo and Ryu [76] improved on that paper by making higher order interpolations of the facial fluxes in order to get higher order edge-centered electric fields. This yielded a fifth order FD-WENO formulation that was indeed fifth order for the flow variables as well as the magnetic field. The above advantages have been made on a case by case basis for each different type of PDE. The *second important goal* of this paper is to show that equations like the induction equation – which give rise to divergence-based involution constraints – all have a common structure. By focusing on that common structure, we find a common solution for the multidimensional upwinding of the induction equation. This common solution, which we develop in this paper, enables us to offer an AFD-WENO based solution for all PDEs that have a divergence-based involution. The benefit from this innovation is vast swathes of divergence-constrained PDEs can all be treated using AFD-WENO schemes of progressively higher order using the same algorithmic framework.

It should also be mentioned that WENO reconstruction and interpolation has seen several recent advances which make it easy to implement such ideas. In Henrick *et al.* [59] and Borges *et al.* [42] it was shown that accuracy can be retained at critical points by modifying the non-linear WENO weights. In Balsara, Garain and Shu [25] it was shown that the smoothness indicators for one-dimensional WENO reconstruction can be written analytically as the sum of perfect squares. In the supplement to Balsara, Samantaray and Subramanian [33] it was shown that the smoothness indicators for two- and three-dimensional WENO reconstruction can also be written analytically as the sum of perfect squares, and the WENO reconstruction was simplified in two and three dimensions. In Balsara *et al*. [36] one-dimensional WENO interpolation was also presented at very high orders. Several authors (Levy, Puppo and Russo [68], Zhu and Qiu [91], Balsara, Garain and Shu [25], Cravero and Semplice [51], Semplice Coco and Russo [75], Zhu and Shu [92]) have shown the value of using multiple stencils of different orders for stabilizing WENO schemes. These advances in WENO reconstruction and WENO interpolation make it easy to formulate and implement the methods presented here. (This distinction between reconstruction and interpolation is very important. In Balsara, Garain and Shu [25] several WENO reconstruction formulae at different orders have been provided that are useful for classical FD-WENO. In Balsara *et al*. [36] analogous pointwise WENO interpolation formulae have been provided at different orders for use in AFD-WENO.)

It is also worth mentioning for the interested reader that there is a very attractive ecosystem of ideas developing around the AFD-WENO methods. Well-balanced AFD-WENO methods for conservation laws have been documented in Xu and Shu [87]. It is also important to have higher order methods that can preserve positivity of density and pressure, and keep other variables within physical bounds. This goes under the rubric of physical constraint preserving (PCP) methods. PCP methods for AFD-WENO have also been documented in Bhoriya *et al*. [41]. Other noteworthy work on PCP schemes in this area included papers by Wu and Shu [85], [86]. The methods have also been extended to extremely large hyperbolic systems like the equations of general relativity in Balsara *et al*. [38]. This illustrates that along with the basic algorithm, there are numerous additional capabilities available for AFD-WENO methods that make them generally useful in a large number of contexts.

The plan of this paper is as follows. In Section 2 we motivate the idea that it is possible to mix finite difference and finite volume WENO approaches in order to get the best advantages of computational efficiency. In Section 3 we discuss the structure of induction equation and the multidimensional dissipation it requires in order to achieve multidimensionally upwinded schemes. In Section 4 we provide a step by step implementation of the scheme. Section 5 describes the additional steps needed for obtaining an efficient and time-explicit PCP algorithm for divergence-preserving PDEs. Section 6 presents accuracy analysis for several hyperbolic PDEs that have a divergence-based involution constraint. Section 7 presents numerous stringent test problems for the same PDE systems. Section 8 ends with some conclusions.

## 2) Motivating a Mixed, i.e. Finite Difference and Finite Volume, Form for Divergence-Preserving PDEs

Practically all users of higher order Godunov schemes would be comfortable with the idea of a higher order finite difference approximation of a conservation law $\partial_t \mathbf{u} + \partial_x \mathbf{f} + \partial_y \mathbf{g} + \partial_z \mathbf{h} = 0$. We consider the spatial discretization of this equation on a uniform mesh with zone size $\Delta x$, $\Delta y$ and $\Delta z$ in the x-, y- and z-directions. The multidimensional finite difference update that is discrete in space but continuous in time looks like:-

$$\partial_t \mathbf{u}_{i,j,k} = -\frac{1}{\Delta x}\left(\mathbf{f}_{i+1/2,j,k} - \mathbf{f}_{i-1/2,j,k}\right) - \frac{1}{\Delta y}\left(\mathbf{g}_{i,j+1/2,k} - \mathbf{g}_{i,j-1/2,k}\right) - \frac{1}{\Delta z}\left(\mathbf{h}_{i,j,k+1/2} - \mathbf{h}_{i,j,k-1/2}\right) \quad (2)$$

The above update term can then be combined with an SSP-RK update to achieve higher order accuracy in time. Although the variables in the above equation don't have the same meaning/interpretation as they would have in a finite volume approximation, we understand that the flux does, nevertheless, satisfy a telescoping property with the result that a notion of conservation can still be asserted. AFD-WENO schemes do have the selling point that they offer high accuracy along with high processing speed for multidimensional problems. Consequently, we begin by claiming that we would like to have ultra-efficient mimetic AFD-WENO schemes for involution-constrained PDEs.

Here we will devise AFD-WENO-like schemes that update zone-centered variables as point values. However, it proves very valuable to retain the integral sense in which eqn. (1) holds. In other words, we want the discrete in space, but continuous in time, version of eqn. (1) to be given by the finite volume-style approximation:-

$$\partial_t \overline{B}_{x;\, i+1/2,j,k} = -\frac{1}{\Delta y \Delta z}\left(\Delta z \overline{E}^{num}_{z;\, i+1/2,j+1/2,k} - \Delta z \overline{E}^{num}_{z;\, i+1/2,j-1/2,k} + \Delta y \overline{E}^{num}_{y;\, i+1/2,j,k-1/2} - \Delta y \overline{E}^{num}_{y;\, i+1/2,j,k+1/2}\right) \quad (3)$$

$$\partial_t \overline{B}_{y;\, i,j-1/2,k} = -\frac{1}{\Delta x \Delta z}\left(\Delta x \overline{E}^{num}_{x;\, i,j-1/2,k+1/2} - \Delta x \overline{E}^{num}_{x;\, i,j-1/2,k-1/2} + \Delta z \overline{E}^{num}_{z;\, i-1/2,j-1/2,k} - \Delta z \overline{E}^{num}_{z;\, i+1/2,j-1/2,k}\right) \quad (4)$$

$$\partial_t \overline{B}_{z;\, i,j,k+1/2} = -\frac{1}{\Delta x \Delta y}\left(\Delta x \overline{E}^{num}_{x;\, i,j-1/2,k+1/2} - \Delta x \overline{E}^{num}_{x;\, i,j+1/2,k+1/2} + \Delta y \overline{E}^{num}_{y;\, i+1/2,j,k+1/2} - \Delta y \overline{E}^{num}_{y;\, i-1/2,j,k+1/2}\right) \quad (5)$$

The magnetic field components with overbars, $\bar{B}_x$, $\bar{B}_y$ and $\bar{B}_z$ in the formulae above, are area averages over the faces of the mesh; see Fig. 1. The numerically stabilized electric field components with overbars, $\bar{E}_x^{num}$, $\bar{E}_y^{num}$ and $\bar{E}_z^{num}$ in the formulae above, are high order accurate line averages along the edges of the mesh; see Fig. 1. The overbars in the above three equations are intended to highlight this averaging. Please realize that without this integral interpretation, Stokes law will not work. The divergence-constraint on the magnetic field will, therefore, be exactly preserved by the above three equations in an integral sense.

We realize, therefore, that eqn. (2) will be satisfied in the sense of point values. Usually, most PDE systems will have a large number of zone-centered conserved variables. Furthermore, these PDEs need to satisfy a free stream condition even when curvilinear logically Cartesian meshes are used. For that reason, AFD-WENO is an optimal solution choice for eqn. (2). This is especially true if stiff sources are also involved, see Balsara *et al.* [35], [36], [37]. So our choice of AFD-WENO will keep the computational cost down to a minimum for the part of the calculation that is computationally very costly. On the other hand, it is rare to have multiple equation sets of the form shown in eqn. (1). Typically, most PDEs of interest have only one divergence-preserving vector field. The only counter-example that we know of is Maxwell's equations, which has two such sets. We do not seek a multidimensional constraint-preserving reconstruction for the divergence-preserving vector field because that could increase the cost if it is not coupled to an efficient ADER predictor step. Instead, we choose a nominal finite volume style discretization for eqn. (1), but we do this using the full set of finite difference tricks at our disposal, thereby reducing the computational cost for evolving all types of divergence-preserving PDEs. Thus our solution strategy will be a hybrid:- We will use AFD-WENO methods for conservation law-like structures in the solution vector; but we will also nominally use a finite volume approach for the evolution of the divergence-preserving parts of the PDE.

**3) Understanding the Structure of the Multidimensional Dissipation in Divergence-Preserving PDEs**

From the discussion in Section 1, and from Fig. 1, we have seen that the electric field components (that are parallel to the edges of the mesh) have to be stabilized in the two directions that are transverse to the field (and indeed transverse to the edges of the mesh). This calls for a two-dimensional Riemann solver to give us a numerically stabilized value for the electric field components along each edge. Here we focus on the z-component of the electric field that lies along the z-edges of the mesh; the analogous expressions for the x- and y-edges can be obtained by cyclic rotations.

In this Section we derive the explicit expressions of the numerically stabilized edge-aligned electric field in the LLF and HLL limits. The HLL limit is very useful because it ensures that we have a good supersonic limit for the electric field that we get from the multidimensional Riemann solver. Our practical experience has been that for stringent MHD test problems it is sometimes beneficial to have a well-designed supersonic limit in the Riemann solver. This is true both for the 1D and 2D Riemann solvers. The Section is divided into three Sub-sections. In Sub-section 3.1 we do some stage-setting. In Sub-section 3.2 we derive the edge-aligned electric field for the two-dimensional LLF Riemann solver. In Sub-section 3.3 we derive the edge-aligned electric field for the two-dimensional HLL Riemann solver. Strategies for endowing the multidimensional HLL Riemann solver with internal sub-structure have also been described in Balsara and Nkonga [27], so we do not need to describe that here.

**3.1) Stage-Setting for the Evaluation of the Edge-centered Electric Field**

We focus on blending AFD-WENO with equations (3), (4) and (5) in the text. The AFD-WENO algorithm already requires us to evaluate the fluxes as point values at each zone center. To that end, it is worth noting that the electric fields are also available from the fluxes, which are evaluated at the zone-centers. Therefore, the electric field components can be interpolated with higher order accuracy and those higher order interpolants can be evaluated at the edges of the mesh. However, note that these zone-centered point values for the electric fields do not have any contribution from the numerical dissipation that is needed to stabilize the scheme. Therefore, in order to efficiently evaluate the numerically stabilized electric field at the edges of the mesh, we have to make two innovations:- First, we have to efficiently interpolate the dissipation-free zone-centered electric field components with high accuracy to the appropriate edges of the mesh. Second, we have to find a way to add numerical dissipation in a way that is highly accurate and proportional only to the jumps in the reconstructed magnetic field components at the edges of the mesh. Both these innovations are needed to stabilize the update in eqn. (1).

Let us illustrate the above two innovations by briefly focusing on the edge-collocated z-component of the electric field in Fig. 1. It has four zones surrounding it, so we will have four z-components of the electric field at those four zone centers. We will need to interpolate those z-components of the electric field with high accuracy to the center of the zone edge shown in Fig. 1. (A strategy for efficiently turning point values into line-integrals will be described later.) This is accomplished via a two-dimensional WENO interpolation in the xy-plane and is, therefore, quite inexpensive. We will soon show that for all equations that are formally similar to eqn. (1) the numerical dissipation is governed entirely by the jumps in the facially reconstructed x- and y-components of the magnetic field at the z-edge being considered in Fig. 1. Notice that the x-component of the facial magnetic field in Fig. 1 only needs to be reconstructed two-dimensionally from the neighboring faces in the yz-directions, whereas the y-component of the facial magnetic field in Fig. 1 only needs to be reconstructed two-dimensionally from the neighboring faces in the xz-directions. These are also inexpensive two-dimensional WENO reconstructions and the result of those reconstructions can be re-used at other neighboring edges. Therefore, this step is also very inexpensive. We pay attention to the very important topic of multidimensional dissipation for all equations that are formally similar to eqn. (1) in this Section.

Eqn. (1) can be written in flux form as:-

$$\frac{\partial}{\partial t}\begin{pmatrix} B_x \\ B_y \\ B_z \end{pmatrix} + \frac{\partial}{\partial x}\begin{pmatrix} 0 \\ -E_z \\ E_y \end{pmatrix} + \frac{\partial}{\partial y}\begin{pmatrix} E_z \\ 0 \\ -E_x \end{pmatrix} + \frac{\partial}{\partial z}\begin{pmatrix} -E_y \\ E_x \\ 0 \end{pmatrix} = 0 \qquad (6)$$

In general, eqn. (6) will be a sub-portion of a larger PDE system. The fluxes in eqn. (6) always have a very special anti-symmetrical form. Notice from eqn. (6) that the second component of the x-flux is just the negative of the first component of the y-flux; the third component of the y-flux is just the negative of the second component of the z-flux; likewise the first component of the z-flux is just the negative of the third component of the x-flux. At a deeper level, these anti-symmetries are an inevitable consequence of the tensorial invariance of Maxwell's equations, but we do not delve into that idea any further here. This form of the fluxes in eqn. (6) is very central to update equations that have the structure shown in eqn. (1). Therefore, we do not posit any constitutive relationship between the vector field "**E**" and the vector field "**B**". (For example, in MHD we have the constitutive relation $\mathbf{E} = -\mathbf{v} \times \mathbf{B}$, with "**v**" as the velocity vector of the fluid. But for the purposes of the argument developed in this Section, we do not use this fact.) Even without asserting a constitutive relationship, our intention is to show that the very structure of eqn. (1) naturally gives us a very special structure for the multidimensional Riemann solver. The anti-symmetry that is central to eqn. (1) gives rise to the dualism that was originally exploited in numerical schemes (Balsara and Spicer [8]). Furthermore, we wish to show that the structure is such as to give us a centered average term along with a very general structure of the dissipation term. (Recall that the 1D LLF Riemann solver can also be written as a centered average flux and a dissipation term. Recall too that this convenient split was then put to good use as a building block for a classical FD-WENO scheme.) In a similar fashion, we wish to use the multi-dimensional Riemann solver of Balsara [17], [18], [22] to write the electric field at the edges as a centered average term along with a multidimensional dissipation term.

To keep the discussion general, we will assume that eqn. (1) is a sub-part of a larger PDE system which sets the wave speeds in all directions. In 2D, the input states are shown in Fig. 2 and the resulting wave structure of the 2D Riemann problem is shown in Fig. 3. Those figures also serve to explain the notation that is used here. The subscript "RU" stands for right-upper; "LU" stands for left-upper; "LD" stands for left-down and "RD" stands for right-down. Fig. 3 shows a wave model with speeds that span $[S_L, S_R] \times [S_D, S_U]$. While we use the HLL version of the Riemann solver, we will also be interested in the LLF variant of this Riemann solver given by setting $S_R = -S_L = S_U = -S_D = S$; where "$S$" is some estimate of the maximal wave speed for the waves around the edge of interest. We will also build in the continuity of the normal component of the magnetic field at the faces of the mesh from the very start so that we have $B_{xRD} = B_{xLD} \equiv B_{xD}$, $B_{xRU} = B_{xLU} \equiv B_{xU}$, $B_{yRU} = B_{yRD} \equiv B_{yR}$ and $B_{yLU} = B_{yLD} \equiv B_{yL}$. This is also shown in Fig. 2. As a result, we realize that the only inputs to the 2D Riemann solver are the four facial magnetic fields $B_{xD}$, $B_{xU}$, $B_{yL}$ and $B_{yR}$ along with the four electric field components $E_{zRU}$, $E_{zLU}$, $E_{zLD}$ and $E_{zRD}$. This dramatically simplifies subsequent derivations. Consequently, the one-dimensional HLL Riemann solver in the x-direction between the two upper states in Fig. 2 gives us the resolved state and flux, which we use to get

$$B_{xU}^* = B_{xU} \ ;$$
$$B_{yU}^* = \left(S_R B_{yR} - S_L B_{yL}\right)/\left(S_R - S_L\right) + \left(E_{zRU} - E_{zLU}\right)/\left(S_R - S_L\right) \ ; \quad (7\text{HLL})$$
$$E_{zU}^* = \left(S_R E_{zLU} - S_L E_{zRU}\right)/\left(S_R - S_L\right) - S_R S_L \left(B_{yR} - B_{yL}\right)/\left(S_R - S_L\right) \ .$$

We can also use $S_R = -S_L = S$ to write the LLF variant of the above equation as

$$B_{xU}^* = B_{xU} \ ;$$
$$B_{yU}^* = \left(B_{yR} + B_{yL}\right)/2 + \left(E_{zRU} - E_{zLU}\right)/(2S) \ ; \quad (7\text{LLF})$$
$$E_{zU}^* = \left(E_{zLU} + E_{zRU}\right)/2 + S\left(B_{yR} - B_{yL}\right)/2 \ .$$

To get the analogous formulae for the x-directional Riemann solver between the two lower states in Fig. 2 we just set $U \to D$ in eqn. (7). When the one-dimensional HLL Riemann solver in the y-direction is applied between the two right states in Fig. 2, we use the resolved state and flux to get

$$B_{xR}^* = \left(S_U B_{xU} - S_D B_{xD}\right)/\left(S_U - S_D\right) - \left(E_{zRU} - E_{zRD}\right)/\left(S_U - S_D\right) \ ;$$
$$B_{yR}^* = B_{yR} \ ; \quad (8\text{HLL})$$
$$E_{zR}^* = \left(S_U E_{zRD} - S_D E_{zRU}\right)/\left(S_U - S_D\right) + S_U S_D \left(B_{xU} - B_{xD}\right)/\left(S_U - S_D\right) \ .$$

We can also use $S_U = -S_D = S$ to write the LLF variant of the above equation as

$$B_{xR}^* = \left(B_{xU} + B_{xD}\right)/2 - \left(E_{zRU} - E_{zRD}\right)/(2S) \ ;$$
$$B_{yR}^* = B_{yR} \ ; \quad (8\text{LLF})$$
$$E_{zR}^* = \left(E_{zRD} + E_{zRU}\right)/2 - S\left(B_{xU} - B_{xD}\right)/2 \ .$$

To get the analogous formulae for the y-directional Riemann solver between the two left states in Fig. 2 we just set $R \to L$ in eqn. (8). Eqns. (7) and (8) give us the resolved states in Fig. 3. These are the states that circumscribe the strongly-interacting state in Fig. 3.

The resolved state in Fig. 3 is most easily obtained by applying eqns. (12), (13) and (14) from Balsara [22]. We synopsize the essential results in this paragraph. We consider a PDE of the form $\partial_t \mathbf{U} + \partial_x \mathbf{F} + \partial_y \mathbf{G} = 0$. When four states come together at the edge of the mesh, four one-dimensional Riemann problems are formed. The one-dimensional Riemann problems evolve self-similarly. However, those four one-dimensional Riemann problems interact amongst themselves, resulting in a strongly interacting state that also evolves self-similarly. The strongly-interacting state is shown by the variables with the double starred superscript in Fig. 3. Since the subsonic case is the most interesting case, and the one that occurs most frequently in the code, we are interested in that strongly interacting state. In the three ensuing equations, the strongly interacting state in Fig. 3 is thought to be covered by a coordinate system in the similarity variables $\tilde{\xi} \equiv x/t$ and $\tilde{\psi} \equiv y/t$. This coordinate system spans $(\tilde{\xi},\tilde{\psi}) \in [S_L, S_R] \times [S_D, S_U]$. The math becomes easier if we pick a coordinate system that is centered on the strongly interacting state. This is done with the help of linear transformations $\xi \equiv (\tilde{\xi} - \xi_c)/\Delta\xi$ and $\psi \equiv (\tilde{\psi} - \psi_c)/\Delta\psi$; with the definitions $\Delta\xi \equiv S_R - S_L$, $\xi_c \equiv (S_R + S_L)/2$, $\Delta\psi \equiv S_U - S_D$ and $\psi_c \equiv (S_U + S_D)/2$. This linear rescaling of the coordinates maps the wave speeds from $[S_L, S_R] \times [S_D, S_U]$ to $[-1/2, 1/2]^2$. For ease of use, we specialize the equations for the strongly interacting state vector $\mathbf{U}^{**}$ and its associated fluxes $\mathbf{F}^{**}$ and $\mathbf{G}^{**}$ for our current needs in this paper as:-

$$\mathbf{U}^{**} = -\left[\begin{array}{l} \dfrac{1}{2\Delta\xi} \int_{-1/2}^{1/2} \left(\mathbf{F}(1/2,\psi) - S_R \mathbf{U}(1/2,\psi)\right) d\psi - \dfrac{1}{2\Delta\xi} \int_{-1/2}^{1/2} \left(\mathbf{F}(-1/2,\psi) - S_L \mathbf{U}(-1/2,\psi)\right) d\psi \\ + \dfrac{1}{2\Delta\psi} \int_{-1/2}^{1/2} \left(\mathbf{G}(\xi,1/2) - S_U \mathbf{U}(\xi,1/2)\right) d\xi - \dfrac{1}{2\Delta\psi} \int_{-1/2}^{1/2} \left(\mathbf{G}(\xi,-1/2) - S_D \mathbf{U}(\xi,-1/2)\right) d\xi \end{array}\right]$$

(9)

$$\mathbf{F}^{**} = \xi_c \mathbf{U}^{**} + \left[\begin{array}{l} \dfrac{1}{2} \int_{-1/2}^{1/2} \left(\mathbf{F}(1/2,\psi) - S_R \mathbf{U}(1/2,\psi)\right) d\psi + \dfrac{1}{2} \int_{-1/2}^{1/2} \left(\mathbf{F}(-1/2,\psi) - S_L \mathbf{U}(-1/2,\psi)\right) d\psi \\ + \dfrac{\Delta\xi}{\Delta\psi} \int_{-1/2}^{1/2} \xi \left(\mathbf{G}(\xi,1/2) - S_U \mathbf{U}(\xi,1/2)\right) d\xi - \dfrac{\Delta\xi}{\Delta\psi} \int_{-1/2}^{1/2} \xi \left(\mathbf{G}(\xi,-1/2) - S_D \mathbf{U}(\xi,-1/2)\right) d\xi \end{array}\right]$$

(10)

$$\mathbf{G}^{**} = \psi_c \mathbf{U}^{**} + \left[\begin{array}{l} \dfrac{\Delta\psi}{\Delta\xi} \int_{-1/2}^{1/2} \psi \left(\mathbf{F}(1/2,\psi) - S_R \mathbf{U}(1/2,\psi)\right) d\psi - \dfrac{\Delta\psi}{\Delta\xi} \int_{-1/2}^{1/2} \psi \left(\mathbf{F}(-1/2,\psi) - S_L \mathbf{U}(-1/2,\psi)\right) d\psi \\ + \dfrac{1}{2} \int_{-1/2}^{1/2} \left(\mathbf{G}(\xi,1/2) - S_U \mathbf{U}(\xi,1/2)\right) d\xi + \dfrac{1}{2} \int_{-1/2}^{1/2} \left(\mathbf{G}(\xi,-1/2) - S_D \mathbf{U}(\xi,-1/2)\right) d\xi \end{array}\right]$$

(11)

$E_z^{**}$ can now be obtained from taking the negative of the second component of $\mathbf{F}^{**}$ or by taking the first component of $\mathbf{G}^{**}$. Here we have simply stated the above three equations with the minimal requisite explanation. In Balsara [22] detailed explanations for the derivations are given.

**3.2) Edge-Aligned Electric Field for the Multidimensional LLF Riemann Solver**

We are first interested in the LLF limit. Notice that $E_z^{**}$ can be obtained as the negative of the second component of the resolved x-flux. This would require use of eqn. (10). Alternatively, $E_z^{**}$ can be obtained as the first component of the resolved y-flux. This would require use of eqn. (11). Regardless of the two alternative ways for obtaining $E_z^{**}$, the mathematics gives us the same result. We get

$$E_z^{**} = \left(E_{zRU} + E_{zLU} + E_{zLD} + E_{zRD}\right)/4 + S\left[B_{xD} - B_{xU} + B_{yR} - B_{yL}\right]/2 \tag{12}$$

For a 2D LLF Riemann solver, $E_z^{**}$ is also the same as he numerically stabilized z-component of the electric field given by $E_z^{num}$. We begin by pointing out that when the variation in the solution is restricted entirely to the x-direction or y-direction, eqn. (12) reduces correctly to the one-dimensional limit. It is, therefore, consistent. Eqn. (12) then shows us that the arithmetically-averaged round bracket, i.e. $\left(E_{zRU} + E_{zLU} + E_{zLD} + E_{zRD}\right)/4$, is the central term. It is free of dissipation and carries the z-component of the electric field that we need at the z-edges of the mesh. This is the electric field that can be evaluated at the zone centers and interpolated in a two-dimensional WENO fashion to the centers of the z-edges, as documented in the figure caption of Fig. 2. WENO is useful for doing this because it ensures that we pick up the largest 2D stencil that is smooth while avoiding possible stencils that might be non-smooth. To do this, we will have to design a pointwise 2D interpolation strategy for WENO; an analogous 2D reconstruction strategy for WENO has been presented in Balsara, Samantaray and Subramanian [33]. The round bracket in eqn. (12) restores consistency to the z-component of the electric field; however, to stabilize it we need a contribution from the dissipation. The dissipation in eqn. (12) is entirely carried by the square bracket, i.e. $\left[B_{xD} - B_{xU} + B_{yR} - B_{yL}\right]/2$. If simple face-centered values are used, we will get first order, dissipation. This dissipation can be very high. However, we can do much better if the facially-reconstructed magnetic field components are evaluated at the zone edges. This can be obtained from the 2D WENO reconstruction of the normal components of the magnetic fields within the 2D faces. We, therefore, see that the additional WENO interpolation and reconstruction strategies that need to be invented are very few or already in hand. They are also very light-weight because the interpolation doesn't need to be projected into characteristic space; rather, it is only applied to a single component of the solution vector or to a single component of the flux vector.

We can additionally write out the x-, y- and z-components of the resolved magnetic field at the z-edge using eqn. (9). For $B_x^{**}$ we get:-

$$B_x^{**} = \left(B_{xU} + B_{xD}\right)/2 + \left(E_{zLD} - E_{zLU} + E_{zRD} - E_{zRU}\right)/(4S) \tag{13}$$

We can see the familiar structure that emerges in a Riemann solver where we see that the flux terms also influence the resolved states in the Riemann problem. For $B_y^{**}$ we get:-

$$B_y^{**} = \left(B_{yR} + B_{yL}\right)/2 + \left(E_{zRD} - E_{zLD} + E_{zRU} - E_{zLU}\right)/(4S) \tag{14}$$

For $B_z^{**}$ we get:-

$$\begin{aligned} B_z^{**} &= \left(B_{zRU} + B_{zLU} + B_{zLD} + B_{zRD}\right)/4 + \left(-E_{xD}^* + E_{xU}^* + E_{yL}^* - E_{yR}^*\right)/(4S) \\ &+ \left[\left(E_{xRU} + E_{xLU}\right) - \left(E_{xRD} + E_{xLD}\right) - \left(E_{yRU} + E_{yRD}\right) + \left(E_{yLU} + E_{yLD}\right)\right]/(8S) \end{aligned} \tag{15}$$

Eqn. (15) above shows that for an involution constrained system, $B_z^{**}$ cannot be obtained from the Riemann solver *per se*, unless more components of the parent PDE are used. However, eqns. (12), (13) and (14) show us that $E_z^{**}$, $B_x^{**}$ and $B_y^{**}$ can indeed be obtained from the multidimensional Riemann solver and that is all that we actually need to make progress with the scheme design. In fairness, even $B_x^{**}$ and $B_y^{**}$ are not useful for the scheme that we design in this paper. However, we provide them here because they could become very useful as and when the inclusion of parabolic terms is considered.

**3.3) Edge-Aligned Electric Field for the Multidimensional HLL Riemann Solver**

For some PDEs, such as the Maxwell equations where the maximal signal speed is always the speed of light, we will always have a subsonic Riemann solver. For other PDEs, such as the MHD equations, the CGL equations, and the relativistic MHD equations, the flow speed can sometimes be strongly supersonic. The LLF Riemann solver can become problematic in such situations because the associated wave model opens up the Riemann fan much more than is warranted. The HLL Riemann solver is very useful in such situations because it has a well-defined supersonic limit. The same is true for the multidimensional HLL Riemann solver. For this reason, we describe the two-dimensional HLL Riemann solver as it applies to eqn. (6).

Using eqns. (7HLL) and (8HLL) and their analogues in eqn. (9), we get the resolved x-component of the magnetic field in the strongly interacting state as:-

$$B_x^{**} = (S_U B_{xU} - S_D B_{xD})/(S_U - S_D) + (E_{zLD} - E_{zLU} + E_{zRD} - E_{zRU})/(2(S_U - S_D)) \qquad (16)$$

The first term in the above equation has a nice geometrical interpretation. It represents the area-weighted average of the x-component of the magnetic field over the wave model shown in Fig. 3. The second term in eqn. (16) can be interpreted as a contribution from the dissipation and is quite analogous to the second term in eqn. (13). We can also obtain the resolved y-component of the magnetic field in the strongly interacting state as:-

$$B_y^{**} = (S_R B_{yR} - S_L B_{yL})/(S_R - S_L) + (-E_{zLD} - E_{zLU} + E_{zRD} + E_{zRU})/(2(S_R - S_L)) \qquad (17)$$

The first and second terms of eqn. (17) can be similarly interpreted. As seen from eqn. (15), it is not valuable to try and obtain $B_z^{**}$, nor is it valuable for the developments that are to follow. The resolved z-component of the electric field in the strongly interacting state can be obtained from the negative of the second component of the x-flux as obtained from eqn. (10). We then get our first approximation for the resolved z-component of the electric field:-

$$\begin{aligned} E_{z;1}^{**} = &-(S_R + S_L) B_y^{**}/2 + (S_U(E_{zLD} + E_{zRD}) - S_D(E_{zLU} + E_{zRU}))/(2(S_U - S_D)) \\ &- S_U S_D (B_{xD} - B_{xU})/(S_U - S_D) + (S_R B_{yR} + S_L B_{yL})/2 \end{aligned} \qquad (18)$$

The resolved z-component of the electric field in the strongly interacting state can also be obtained from the first component of the y-flux as obtained from eqn. (11). We then get our second approximation for the resolved z-component of the electric field:-

$$\begin{aligned} E_{z;2}^{**} = &(S_U + S_D) B_x^{**}/2 + (S_R(E_{zLD} + E_{zLU}) - S_L(E_{zRD} + E_{zRU}))/(2(S_R - S_L)) \\ &- (S_U B_{xU} + S_D B_{xD})/2 - S_R S_L (B_{yR} - B_{yL})/(S_R - S_L) \end{aligned} \qquad (19)$$

Notice that when the wave speeds revert to their LLF limits, $S_R = -S_L = S_U = -S_D = S$, eqns. (18) and (19) revert to the LLF formula in eqn. (12). For the HLL case, the expressions in eqns. (18) and (19) are not identical and it is best to take the arithmetic average of the two expressions as our resolved z-component of the electric field in the strongly interacting state.

The LLF Riemann solver in 1D or 2D is sometimes a slightly weaker option because it lacks a supersonic limit. The above expressions show us how to obtain the edge-aligned z-component of the electric field from the 2D HLL Riemann solver in the fully subsonic case where $S_L \leq 0 \leq S_R$ and $S_D \leq 0 \leq S_U$, as shown in Fig. 3. However, there are also eight possible supersonic cases. They are shown in Fig. 4. The resolved z-component of the electric field that is to be used for numerical work, $E_z^{num}$, can then by obtained by the following nine fold stack of conditionals that can be evaluated as:-

$$
\begin{aligned}
&if\ (S_L \geq 0\ \text{and}\ S_D \geq 0)\ then \\
&\quad E_z^{num} = E_{zLD}\ ;\ return\ ; \\
&elseif\ (S_R \leq 0\ \text{and}\ S_D \geq 0)\ then \\
&\quad E_z^{num} = E_{zRD}\ ;\ return\ ; \\
&elseif\ (S_R \leq 0\ \text{and}\ S_U \leq 0)\ then \\
&\quad E_z^{num} = E_{zRU}\ ;\ return\ ; \\
&elseif\ (S_L \geq 0\ \text{and}\ S_U \leq 0)\ then \\
&\quad E_z^{num} = E_{zLU}\ ;\ return\ ; \\
&elseif\ (S_L \geq 0)\ then \\
&\quad E_z^{num} = E_{zL}^*\ ;\ return\ ; \\
&elseif\ (S_R \leq 0)\ then \\
&\quad E_z^{num} = E_{zR}^*\ ;\ return\ ; \\
&elseif\ (S_D \geq 0)\ then \\
&\quad E_z^{num} = E_{zD}^*\ ;\ return\ ; \\
&elseif\ (S_U \leq 0)\ then \\
&\quad E_z^{num} = E_{zU}^*\ ;\ return\ ; \\
&else \\
&\quad E_z^{num} = \left(E_{z;1}^{**} + E_{z;2}^{**}\right)/2\ ;\ return\ ; \\
&endif
\end{aligned} \quad (20)
$$

The subsonic limit, provided by the strongly interacting state in the last conditional above, is the most important case in the design of a multidimensional Riemann solver since most problems operate in that limit most of the time. The stack of conditionals in the above equation is designed to provide optimal efficiency in computer code. The conditionals in the above stack are designed to be mutually exclusive. As soon as one of the conditionals checks out as positive, the subroutine should evaluate that numerical $E_z^{num}$ and return without considering the subsequent conditionals. Therefore, only the variable that occurs within one of the conditional statements needs to be evaluated when that particular conditional is invoked.

**4) Step by Step Implementation of the Scheme**

All the requisite pieces that go into the algorithm are now in hand. Here we describe an AFD-WENO scheme that uses the pointwise collocated zone-centered variables and facially averaged normal components of vector fields that are divergence-constraint preserving. (If the PDE system does not have any non-conservative products, the AFD-WENO formulation may even be replaced with a FD-WENO formulation. This is accomplished by replacing Step 3 below with a standard flux-reconstruction-based FD-WENO algorithm to obtain fluxes.) In general, the conserved fluid variables are zone-centered point values while the divergence-preserving normal components of the vector fields are facially averaged. These are the primal variables of the scheme. This is the most economical choice and it allows us to design a scheme where the zone-centered variables follow the AFD-WENO update while the facially collocated variables will be updated in a finite volume sense. While this split approach may seem odd, it has two major advantages. First, it retains most of the efficiency advantage of the AFD-WENO scheme; since this is a costly step, it yields a substantial savings. Second, since high order divergence-preserving AMR of vector fields is a solved problem (see Balsara Samantaray and Subramanian [33], Balsara and Sarris [34]), all those advantages are retained. We assume a 3D Cartesian mesh with $N_x \times N_y \times N_z$ zones with uniform zone sizes $\Delta x$, $\Delta y$ and $\Delta z$ in the x-, y- and z-directions. To keep the notation tractable, we will establish a convention that applies to this Section. We will denote averaged quantities with capital letters with an overbar while we denote point values with corresponding small letters. Thus we start with zone-centered, pointwise, conserved variables $\mathbf{u}_{i,j,k}$ for $(i,j,k) \in [1, N_x] \times [1, N_y] \times [1, N_z]$ and facially-averaged variables for the magnetic field denoted by $\overline{B}_{x;i+1/2,j,k}$ for $(i,j,k) \in [0, N_x] \times [1, N_y] \times [1, N_z]$, $\overline{B}_{y;i,j+1/2,k}$ for $(i,j,k) \in [1, N_x] \times [0, N_y] \times [1, N_z]$ and $\overline{B}_{z;i,j,k+1/2}$ for $(i,j,k) \in [1, N_x] \times [1, N_y] \times [0, N_z]$. Notice the difference between the zone-centered pointwise small letters and the facially-averaged capital letters with overbars. Here we describe the algorithm within the context of the ultra-simple two-dimensional LLF Riemann solver because that is what most people will implement first. (Since this is an AFD-WENO scheme, it can accommodate any type of pointwise Riemann solver. So one can initially implement a two-dimensional LLF Riemann solver and subsequently replace it with a two-dimensional HLL Riemann solver.)

We illustrate the algorithm at fifth order, but all the steps are fully generalizable to all orders. The algorithm proceeds with the following steps:-

**Step 1) 2D Reconstruction of Facial B Field Components to Obtain Face-centered Point Values:** At each face, we can make the area-averaged, two-dimensional reconstruction for the facial component of the magnetic field. The 2D high order WENO-based reconstruction strategy that is needed for this step has already been described in the Supplement to Balsara Samantaray and Subramanian [33]. For example, at the x-face, and up to $5^{th}$ order, we can make a 2D finite volume WENO reconstruction of all the facial modes as:-

$$
\begin{aligned}
B_{x;i+1/2,j,k}(y,z) &= \overline{B}_{x;i+1/2,j,k} + b_y y + b_z z + b_{yy}(y^2 - 1/12) + b_{zz}(z^2 - 1/12) + b_{yz} yz \\
&+ b_{yyy}(y^3 - 3y/20) + b_{zzz}(z^3 - 3z/20) + b_{yyz}(y^2 - 1/12)z + b_{yzz} y(z^2 - 1/12) \\
&+ b_{yyyy}(y^4 - 3y^2/14 + 3/560) + b_{zzzz}(z^4 - 3z^2/14 + 3/560) + b_{yyyz}(y^3 - 3y/20)z \\
&+ b_{yzzz} y(z^3 - 3z/20) + b_{yyzz}(y^2 - 1/12)(z^2 - 1/12)
\end{aligned}
\quad (21)
$$

Note that because this is a 2D area-averaged WENO reconstruction, the area-averaged quantity, $\overline{B}_{x;i+1/2,j,k}$, is set. The WENO reconstruction is only tasked with finding the best approximations for the other coefficients on the right hand side of eqn. (21). From eqn. (21) we can evaluate and store the point value of the x-component of the magnetic field at the center of the x-face as:-

$$b_{x;i+1/2,j,k} = B_{x;i+1/2,j,k}(y=0, z=0) = \overline{B}_{x;i+1/2,j,k} - b_{yy}/12 - b_{zz}/12 + 3\,b_{yyyy}/560 + 3\,b_{zzzz}/560 + b_{yyzz}/144$$
(22)

Analogously, we can obtain $b_{y;i,j+1/2,k}$ at the center of the y-face and $b_{z;i,j,k+1/2}$ at the center of the z-face.

Eqn. (21) can also be used to evaluate the facial magnetic field component anywhere within a face. Specifically, we also use eqn. (21) to evaluate $B_{x;i+1/2,j,k}(y=0, z=1/2)$, $B_{x;i+1/2,j,k}(y=0, z=-1/2)$, $B_{x;i+1/2,j,k}(y=1/2, z=0)$ and $B_{x;i+1/2,j,k}(y=-1/2, z=0)$ at the edge-centers of the four edges that surround the x-face in question. These four values should be stored for each x-face. Similar evaluations can be made and stored at the y- and z-faces. These edge-centered evaluations of the magnetic fields will help us with the evaluation of the dissipation terms for the numerical electric field, as shown in eqn. (23) of Step 6 below.

**Step 2) 1D Interpolation from Facial Point Values to Zone-Centered Point Values:** With these facial point values in hand from eqn. (22), we can use the special type of 1D WENO interpolation that was developed in Section 4 of Balsara *et al*. [36] to obtain point values of $b_{x;i,j,k}$, $b_{y;i,j,k}$ and $b_{z;i,j,k}$ at the zone centers. (To take MHD as an example, $\mathbf{u}_{i,j,k}$ can be used to obtain the pointwise fluid density, velocity and pressure at each zone center. The interpolated $b_{x;i,j,k}$, $b_{y;i,j,k}$ and $b_{z;i,j,k}$ are also available through this step, so we have all the MHD variables that are needed at the zone centers.) This allows us to evaluate pointwise values for the fluxes, electric fields and any other variables that we might need at the zone centers.

**Step 3) 2D Interpolation of the Three Electric Field Components; Electric Fields Input to 2D Riemann Solver:** Recall that the zone-centered electric fields can always be obtained as various components of the zone-centered fluxes that were obtained in Step 2. Now that zone-centered pointwise values of the electric field, i.e. $e_{x;i,j,k}$, $e_{y;i,j,k}$ and $e_{z;i,j,k}$ have been evaluated, we make three 2D WENO interpolations of those variables in the yz-plane, the xz-plane and the xy-plane respectively. Specifically, within each zone, we interpolate $e_{x;i,j,k}$ in the yz-plane; we interpolate $e_{y;i,j,k}$ in the xz-plane; and we interpolate $e_{z;i,j,k}$ in the xy-plane. This enables us to obtain those same variables at edge-centered locations. (Such an efficient 2D WENO interpolation has been documented in Appendix A of this paper, at $3^{rd}$, $5^{th}$ and $7^{th}$ order; and higher order extensions, such as ninth order, are easily made. Figs. 5 and 6 show the kinds of stencils that are needed for the high order interpolation.) The point of this interpolation is that it gives us the edge-centered point values of the electric field. At each edge we have four abutting zones, and from each of those zones we will get one edge-aligned component of the electric field. Say we consider the z-edge given by $(i+1/2, j+1/2, k)$, as shown in Fig. 2. At that z-edge, we will have four different values of the z-component of the electric field. We have $e_{z;RU;i+1/2,j+1/2,k}$, which is interpolated with high order accuracy from the "RU" zone relative to this edge. We have $e_{z;LU;i+1/2,j+1/2,k}$, which is interpolated with high order accuracy from the "LU" zone relative to this edge. We have $e_{z;LD;i+1/2,j+1/2,k}$, which is interpolated with high order accuracy from the "LD" zone relative to this edge. And we have $e_{z;RD;i+1/2,j+1/2,k}$, which is interpolated with high order accuracy from the "RD" zone relative to this edge. Realize, therefore, from eqn. (12) that we have obtained everything that we need for evaluating the centered part of the electric field components in a pointwise fashion at the edge-centers. This centered part is shown by the round bracket in eqn. (12). Realize that this averaged electric field is not yet stabilized because we don't yet have the multidimensional dissipation term that is needed in the multidimensional Riemann solver. The process of obtaining the stabilizing terms, i.e. the square bracket in eqn. (12), will be described in the next point.

**Step 4) Reconstructed Magnetic Fields to Provide Inputs for the Dissipation Terms in the 2D Riemann solver:** Now focus on Fig. 2. Recall that the divergence-preserving condition ensures that the normal component of the magnetic field is continuous across the faces of the mesh. Therefore, along the x = constant surfaces of the mesh, we can use our 2D WENO reconstruction (from Step 1) of $\bar{B}_{x;i+1/2,j+1,k}$ in the face $(i+1/2, j+1, k)$ to obtain the point value $b_{xU;i+1/2,j+1/2,k}$ from above the black dot in Fig. 2. Similarly, we can use our 2D WENO reconstruction of $\bar{B}_{x;i+1/2,j,k}$ in the face $(i+1/2, j, k)$ to obtain the point value $b_{xD;i+1/2,j+1/2,k}$ from below the black dot in Fig. 2. Likewise, along the y = constant surfaces of the mesh, we can use our 2D WENO reconstruction (from Step 1) of $\bar{B}_{y;i+1,j+1/2,k}$ in the face $(i+1, j+1/2, k)$ to obtain the point value $b_{yR;i+1/2,j+1/2,k}$ from the right of the black dot in Fig. 2 and we can use another 2D WENO reconstruction of $\bar{B}_{y;i,j+1/2,k}$ in the face $(i, j+1/2, k)$ to obtain $b_{yL;i+1/2,j+1/2,k}$ from the left of the black dot in Fig. 2. When the solution is smooth, and when the interpolation is high order, we will have $b_{xU;i+1/2,j+1/2,k} \to b_{xD;i+1/2,j+1/2,k}$ and $b_{yR;i+1/2,j+1/2,k} \to b_{yL;i+1/2,j+1/2,k}$, with the result that the dissipation is significantly reduced. The end result is that at the edge-centers of the mesh we have the centered part of the z-component of the electric field, i.e. the round bracket in eqn. (12). We also have the dissipation term, i.e. the square bracket in eqn. (12).

**Step 5) Standard AFD-WENO Algorithm to Obtain Fluxes (or Fluctuations as needed):** With the above items in hand, it is easy to see that the AFD-WENO algorithm from Balsara *et al*. [36] or Balsara *et al*. [37] can be used to obtain the time rate of change, $(\partial_t \mathbf{u})_{i,j,k}$ for all the zone-centered variables. This is usually the most expensive step in the algorithm because most PDEs have a large number of zone-centered primal variables. (This step usually becomes rather costly due to the fact that the first interpolation in this step has to be done in the characteristic variables.) Our use of the very economical AFD-WENO algorithm reduces the overall cost of the scheme.

**Step 6) Use Steps 1 and 3, along with Speeds from Step 5 to Obtain Numerical Electric Field:** Notice that in Step 1, the facial magnetic field components were reconstructed and in Step 4 we showed how this reconstruction can be used to evaluate the numerical dissipation terms in the electric field. In Step 3 the electric fields were interpolated so that we can obtain the centered part of the electric field at each edge. This enables us to obtain the four input electric fields that make up the centered part of the numerical electric field at the z-edge. The jumps in the magnetic field components at the z-edge also gives us the numerical dissipation. Step 5 also gave us the speeds in all the directions at a z-edge. As a result, the pointwise, stabilized electric field (that is useful for numerical schemes) can be written at the center of the z-edge as:-

$$e^{num}_{z;i+1/2,j+1/2,k} \equiv \left(e_{z;RU;i+1/2,j+1/2,k} + e_{z;LU;i+1/2,j+1/2,k} + e_{z;LD;i+1/2,j+1/2,k} + e_{z;RD;i+1/2,j+1/2,k}\right)/4$$
$$+ S\left[b_{xD;i+1/2,j+1/2,k} - b_{xU;i+1/2,j+1/2,k} + b_{yR;i+1/2,j+1/2,k} - b_{yL;i+1/2,j+1/2,k}\right]/2 \quad (23)$$

This numerical electric field $e^{num}_{z;i+1/2,j+1/2,k}$ is what we were seeking in this step and the previous step. Observe that this $e^{num}_{z;i+1/2,j+1/2,k}$ is still a point value. However, it is a numerically stabilized point value. Similar edge-centered electric fields can be obtained at the x- and y-edges.

**Step 7) From Point Values of the Electric Fields to Edge-Averaged Values:** Notice that the update in eqns. (3), (4) and (5) are based on an area-averaged interpretation of eqn. (1). Up to this point in the discussion, we only have point values of the numerically stabilized electric field components at the centers of each edge; see Fig. 7. Focusing on the z-component of the numerical electric field, we can now invoke pointwise 1D WENO interpolation along each z-edge. This pointwise 1D WENO interpolation algorithm has been described in Section 3 of Balsara et al. [36]. For instance, at 5$^{th}$ order, the resulting interpolating polynomial can be written as

$$E_{z;i+1/2,j+1/2,k}(z) = \bar{E}^{num}_{z;i+1/2,j+1/2,k} + e_z z + e_{zz}(z^2 - 1/12) + e_{zzz}(z^3 - 3z/20) + e_{zzzz}(z^2 - 3z^2/14 + 3/560)$$

(24)

Note that because this is a 1D pointwise WENO interpolation, all the coefficients on the right hand side of eqn. (24), including $\bar{E}^{num}_{z;i+1/2,j+1/2,k}$, have to be evaluated. We can do this same operation along the x- and y-edges too. The term $\bar{E}^{num}_{z;i+1/2,j+1/2,k}$ will naturally be the edge-averaged electric field along the z-edge. This is the variable that we seek in the update in eqns. (3), (4) and (5). Note that we are using interpolation, not reconstruction; as a result, $\bar{E}^{mum}_{z;i+1/2,j+1/2,k}$ from eqn. (24) usually will not be same as $e^{num}_{z;i+1/2,j+1/2,k}$ from eqn. (23).

This completes our description of the scheme.

The scheme, as described here, can do several nice test problems and several scientific applications without further modification. For that reason, we have only described the baseline scheme here. However, there exist stringent MHD problems where one has very strong magnetic fields or very high speed flows, resulting in very low plasma-$\beta$. Likewise, in several relativistic MHD problems, the presence of a very strong magnetic field results in very high magnetization of the plasma. We have also been able to extend the physical constraint preserving (PCP) property from our work in Bhoriya *et al*. [41] to encompass systems with a divergence-preserving constraint. This is described in the next Section.

In this work we have retained the zone-centered point values as the primal variables of our scheme. In some applications it may be very useful to start with zone-centered point values and obtain a full finite volume reconstruction within a zone. This could be useful if additional physics needs to be specified within a zone, or if the zonal values have to be projected from one mesh to another mesh. In light of McCorquodale and Colella [69] and Buchmüller, Dreher and Helzel [43], one might be given to believe that this can only be accomplished up to fourth order. This is not true. To drive this point home, Appendix B illustrates how the WENO formulation can be used to accomplish this at fifth order; but the procedure displayed in Appendix B can be carried out at all orders. One may also desire to obtain area averaged fluxes from point values of fluxes at any order. In Appendix C we again show that the WENO method can be used to do that up to fifth order; and the method is extensible to all orders. Our goal in these Appendices is to show that there might be a deeper connection between finite volume formulations and finite difference formulations at all orders and to illustrate how WENO methods can be used malleably to establish that connection.

**5) An Efficient and Time-Explicit Physical Constraint Preserving (PCP) Strategy for Divergence-Preserving PDEs**

The physical constraint preserving (PCP) property is built around the notion that there exists a first order scheme that preserves the PCP property. Therefore, whenever we are likely to lose physical realizability in a higher order scheme, we are willing to trade a decrease in accuracy in favor of physical realizability. A step by step implementation strategy for PCP in the numerical solution of conservation laws is given in Section 2 of Bhoriya *et al*. [41]. The core idea derives from Hu, Adams and Shu [61]. Section 5 of Bhoriya *et al*. [41] also shows how a higher order in time SSP-RK update can be written as a linear combination of forward Euler steps. All those ideas are useful for the next two Sub-Sections.

**5.1) A First Order PCP Scheme for Divergence-Free MHD**

It is important to recall that all PCP methods for finite difference applications have in some way or form used the Hu, Adams and Shu [61] idea of hybridization between a low order scheme with a high order scheme. But the plan of carrying out such a hybridization for a divergence-constrained PDE is more intricate because the facial magnetic field variables are indeed the primal variables. A first order, zone-centered scheme for MHD that is updated with facial fluxes from an HLL or LLF Riemann solver will indeed be PCP; as was first proved by Gurski [56]. Such a scheme is very lightweight because it does not entail any reconstruction or interpolation steps. The first, and most important, step is to obtain a PCP scheme for MHD that retains the facial magnetic fields as the primal variables.

For any first order scheme we can always average facial magnetic fields to obtain zone-centered magnetic fields. Because the facial magnetic fields are our primal magnetic fields, those fields will have to be updated using first order accurate edge-centered electric fields, as described in eqns. (18) to (20), or more simply in eqn. (12). The evaluation of the 2D Riemann solver at each edge is also very lightweight. The zone-centered magnetic field will then have to be obtained by arithmetic averaging of the facial magnetic field components to the zone center. Even for the most stringent of problems, such a scheme will remain PCP in most of the zones. There will be a rare few zones where the averaging of the facial magnetic fields to the zone centers will indeed destroy the PCP property. To get a bullet-proof scheme, we have to find a way to cure this problem. (Documenting such a low order PCP scheme is the task of this Sub-Section. Describing how it is hybridized with a higher order scheme that is not PCP is the task of the next Sub-Section.)

Let us explain the previous paragraph in more mathematical terms so that the reader can appreciate the issue that we discuss here. The discussion here will use a first order in time forward Euler approach because a higher order in time can be achieved using an SSP-RK method which combines forward Euler methods. We start with a first order version of eqn. (2) to obtain $\mathbf{u}_{i,j,k}^{n;LO}$. (The superscript "LO" stands for low order; the superscript "HO" stands for high order.) However, for a divergence-preserving scheme with facial magnetic fields as the primal variables, the sixth, seventh and eighth components of $\mathbf{u}_{i,j,k}^{n;LO}$ will be replaced by $\left(\overline{B}_{x;i+1/2,j,k}^{n;LO} + \overline{B}_{x;i-1/2,j,k}^{n;LO}\right)/2$, $\left(\overline{B}_{y;i,j+1/2,k}^{n;LO} + \overline{B}_{y;i,j-1/2,k}^{n;LO}\right)/2$ and $\left(\overline{B}_{z;i,j,k+1/2}^{n;LO} + \overline{B}_{z;i,j,k-1/2}^{n;LO}\right)/2$ respectively. We start by assuming that $\mathbf{u}_{i,j,k}^{n;LO}$ is PCP. We then make the forward Euler updates:-

$$\tilde{\mathbf{u}}_{i,j,k}^{n+1;LO} = \mathbf{u}_{i,j,k}^{n;LO} - \frac{\Delta t}{\Delta x}\left(\mathbf{f}_{i+1/2,j,k}^{LO} - \mathbf{f}_{i-1/2,j,k}^{LO}\right) - \frac{\Delta t}{\Delta y}\left(\mathbf{g}_{i,j+1/2,k}^{LO} - \mathbf{g}_{i,j-1/2,k}^{LO}\right) - \frac{\Delta t}{\Delta z}\left(\mathbf{h}_{i,j,k+1/2}^{LO} - \mathbf{h}_{i,j,k-1/2}^{LO}\right) \quad (25)$$

$$\overline{B}_{x;i+1/2,j,k}^{n+1;LO} = \overline{B}_{x;i+1/2,j,k}^{n;LO} - \frac{\Delta t}{\Delta y \Delta z}\left(\Delta z \overline{E}_{z;i+1/2,j+1/2,k}^{LO} - \Delta z \overline{E}_{z;i+1/2,j-1/2,k}^{LO} + \Delta y \overline{E}_{y;i+1/2,j,k-1/2}^{LO} - \Delta y \overline{E}_{y;i+1/2,j,k+1/2}^{LO}\right)$$

(26)

$$\overline{B}_{y;i,j-1/2,k}^{n+1;LO} = \overline{B}_{y;i,j-1/2,k}^{n;LO} - \frac{\Delta t}{\Delta x \Delta z}\left(\Delta x \overline{E}_{x;i,j-1/2,k+1/2}^{LO} - \Delta x \overline{E}_{x;i,j-1/2,k-1/2}^{LO} + \Delta z \overline{E}_{z;i-1/2,j-1/2,k}^{LO} - \Delta z \overline{E}_{z;i+1/2,j-1/2,k}^{LO}\right)$$

(27)

$$\overline{B}_{z;i,j,k+1/2}^{n+1;LO} = \overline{B}_{z;i,j,k+1/2}^{n;LO} - \frac{\Delta t}{\Delta x \Delta y}\left(\Delta x \overline{E}_{x;i,j-1/2,k+1/2}^{LO} - \Delta x \overline{E}_{x;i,j+1/2,k+1/2}^{LO} + \Delta y \overline{E}_{y;i+1/2,j,k+1/2}^{LO} - \Delta y \overline{E}_{y;i-1/2,j,k+1/2}^{LO}\right)$$

(28)

Now realize that because $\mathbf{u}_{i,j,k}^{n;LO}$ started off PCP, and because eqn. (25) is a first order update with a PCP-preserving HLL or LLF flux, we are guaranteed that the updated zone-centered variable $\tilde{\mathbf{u}}_{i,j,k}^{n+1;LO}$ will be PCP. It is only when the sixth, seventh and eighth components of $\tilde{\mathbf{u}}_{i,j,k}^{n+1;LO}$ are replaced with $\left(\overline{B}_{x;i+1/2,j,k}^{n+1;LO} + \overline{B}_{x;i-1/2,j,k}^{n+1;LO}\right)/2$, $\left(\overline{B}_{y;i,j+1/2,k}^{n+1;LO} + \overline{B}_{y;i,j-1/2,k}^{n+1;LO}\right)/2$ and $\left(\overline{B}_{z;i,j,k+1/2}^{n+1;LO} + \overline{B}_{z;i,j,k-1/2}^{n+1;LO}\right)/2$ from eqns. (26), (27) and (28) that we have the possibility of losing the PCP property in a few rare zones! After the replacement of the sixth, seventh and eighth components of $\tilde{\mathbf{u}}_{i,j,k}^{n+1;LO}$ by the facial averages, let us call the zone-centered variables $\mathbf{u}_{i,j,k}^{n+1;LO}$. So it is this replacement by the facial averages that may, in a few rare zones, turn $\tilde{\mathbf{u}}_{i,j,k}^{n+1;LO}$ (which is PCP) into $\mathbf{u}_{i,j,k}^{n+1;LO}$ (which may not be PCP). Our task is to fix this situation.

Two clear fixes are obvious. The first fix draws upon work by Abgrall [1], [2]. Realize that when $\mathbf{u}_{i,j,k}^{n+1;LO}$ loses the PCP property, it does so because of the slightest discretization error; i.e. the facial fluxes that were used to update the zone-centered energy did not bring in enough energy from neighboring zones to keep $\mathbf{u}_{i,j,k}^{n+1;LO}$ PCP. Let $[\mathbf{X}]_5$ denote the fifth component of a vector $\mathbf{X}$; this is just a choice of notation. We can ever so slightly modify the fifth component of the flux for the MHD case as follows:-

$$\begin{aligned}
&\left[\tilde{\mathbf{f}}_{i+1/2,j,k}^{LO}\right]_5 = \left[\mathbf{f}_{i+1/2,j,k}^{LO}\right]_5 - \overline{w}_{i+1/2,j,k}^{-}\alpha \quad ; \quad \left[\tilde{\mathbf{f}}_{i-1/2,j,k}^{LO}\right]_5 = \left[\mathbf{f}_{i-1/2,j,k}^{LO}\right]_5 + \overline{w}_{i-1/2,j,k}^{+}\alpha \quad ; \\
&\left[\tilde{\mathbf{g}}_{i,j+1/2,k}^{LO}\right]_5 = \left[\mathbf{g}_{i,j+1/2,k}^{LO}\right]_5 - \overline{w}_{i,j+1/2,k}^{-}\alpha \quad ; \quad \left[\tilde{\mathbf{g}}_{i,j-1/2,k}^{LO}\right]_5 = \left[\mathbf{g}_{i,j-1/2,k}^{LO}\right]_5 + \overline{w}_{i,j-1/2,k}^{+}\alpha \quad ; \quad (29)\\
&\left[\tilde{\mathbf{h}}_{i,j,k+1/2}^{LO}\right]_5 = \left[\mathbf{h}_{i,j,k+1/2}^{LO}\right]_5 - \overline{w}_{i,j,k+1/2}^{-}\alpha \quad ; \quad \left[\tilde{\mathbf{h}}_{i,j,k-1/2}^{LO}\right]_5 = \left[\mathbf{h}_{i,j,k-1/2}^{LO}\right]_5 + \overline{w}_{i,j,k-1/2}^{+}\alpha
\end{aligned}$$

The $\overline{w}$ variables will be described shortly and they are just a measure of a neighboring zone's ability to give enough energy to the troubled zone so that the troubled zone has a positive pressure. We can then quantify the variable "$\alpha$" in the above equation by setting it by using the equation:-

$$\begin{aligned}
\frac{1}{8\pi}&\left\{\left(\left(\overline{B}_{x;i+1/2,j,k}^{n+1;LO} + \overline{B}_{x;i-1/2,j,k}^{n+1;LO}\right)/2\right)^2 + \left(\left(\overline{B}_{y;i,j+1/2,k}^{n+1;LO} + \overline{B}_{y;i,j-1/2,k}^{n+1;LO}\right)/2\right)^2 \right. \\
&\left. + \left(\left(\overline{B}_{z;i,j,k+1/2}^{n+1;LO} + \overline{B}_{z;i,j,k-1/2}^{n+1;LO}\right)/2\right)^2 - \left(\left[\tilde{\mathbf{u}}_{i,j,k}^{n+1;LO}\right]_6\right)^2 - \left(\left[\tilde{\mathbf{u}}_{i,j,k}^{n+1;LO}\right]_7\right)^2 - \left(\left[\tilde{\mathbf{u}}_{i,j,k}^{n+1;LO}\right]_8\right)^2\right\} = \quad (30)\\
&\alpha\left\{\frac{\Delta t}{\Delta x}\left(\overline{w}_{i+1/2,j,k}^{-} + \overline{w}_{i-1/2,j,k}^{+}\right) + \frac{\Delta t}{\Delta y}\left(\overline{w}_{i,j+1/2,k}^{-} + \overline{w}_{i,j-1/2,k}^{+}\right) + \frac{\Delta t}{\Delta z}\left(\overline{w}_{i,j,k+1/2}^{-} + \overline{w}_{i,j,k-1/2}^{+}\right)\right\}
\end{aligned}$$

Then we can set

$$w^{-}_{i+1/2,j,k}=\left(\max\left(P^{n;LO}_{i+1,j,k}-P_{\min},0\right)\right)^2 \;;\; w^{+}_{i-1/2,j,k}=\left(\max\left(P^{n;LO}_{i-1,j,k}-P_{\min},0\right)\right)^2 \;;\; w^{-}_{i,j+1/2,k}=\left(\max\left(P^{n;LO}_{i,j+1,k}-P_{\min},0\right)\right)^2 \;;$$
$$w^{+}_{i,j-1/2,k}=\left(\max\left(P^{n;LO}_{i,j-1,k}-P_{\min},0\right)\right)^2 \;;\; w^{-}_{i,j,k+1/2}=\left(\max\left(P^{n;LO}_{i,j,k+1}-P_{\min},0\right)\right)^2 \;;\; w^{+}_{i,j,k-1/2}=\left(\max\left(P^{n;LO}_{i,j,k-1}-P_{\min},0\right)\right)^2$$

(31)

Here $P_{\min}$ is a user-settable small parameter. (We have set $P_{\min}=10^{-3}$ for all the tests reported here.) The weights with overbars can then be normalized so that they add to unity in a WENO-like normalization that goes as follows

$$\overline{w}^{-}_{i+1/2,j,k}=w^{-}_{i+1/2,j,k}\Big/\left(w^{-}_{i+1/2,j,k}+w^{+}_{i-1/2,j,k}+w^{-}_{i,j+1/2,k}+w^{+}_{i,j-1/2,k}+w^{-}_{i,j,k+1/2}+w^{+}_{i,j,k-1/2}\right)$$
$$\overline{w}^{+}_{i-1/2,j,k}=w^{+}_{i-1/2,j,k}\Big/\left(w^{-}_{i+1/2,j,k}+w^{+}_{i-1/2,j,k}+w^{-}_{i,j+1/2,k}+w^{+}_{i,j-1/2,k}+w^{-}_{i,j,k+1/2}+w^{+}_{i,j,k-1/2}\right)$$
$$\overline{w}^{-}_{i,j+1/2,k}=w^{-}_{i,j+1/2,k}\Big/\left(w^{-}_{i+1/2,j,k}+w^{+}_{i-1/2,j,k}+w^{-}_{i,j+1/2,k}+w^{+}_{i,j-1/2,k}+w^{-}_{i,j,k+1/2}+w^{+}_{i,j,k-1/2}\right)$$
$$\overline{w}^{+}_{i,j-1/2,k}=w^{+}_{i,j-1/2,k}\Big/\left(w^{-}_{i+1/2,j,k}+w^{+}_{i-1/2,j,k}+w^{-}_{i,j+1/2,k}+w^{+}_{i,j-1/2,k}+w^{-}_{i,j,k+1/2}+w^{+}_{i,j,k-1/2}\right)$$
$$\overline{w}^{-}_{i,j,k+1/2}=w^{-}_{i,j,k+1/2}\Big/\left(w^{-}_{i+1/2,j,k}+w^{+}_{i-1/2,j,k}+w^{-}_{i,j+1/2,k}+w^{+}_{i,j-1/2,k}+w^{-}_{i,j,k+1/2}+w^{+}_{i,j,k-1/2}\right)$$
$$\overline{w}^{+}_{i,j,k-1/2}=w^{+}_{i,j,k-1/2}\Big/\left(w^{-}_{i+1/2,j,k}+w^{+}_{i-1/2,j,k}+w^{-}_{i,j+1/2,k}+w^{+}_{i,j-1/2,k}+w^{-}_{i,j,k+1/2}+w^{+}_{i,j,k-1/2}\right)$$

(32)

As long as a troubled zone has as at least one neighboring zone that can lend it some thermal energy, it will be able to obtain some extra energy to restore its positive pressure. This is a fully conservative fix.

It is also possible that a zone may not have any von Neumann neighbors that have sufficient energy to lend it some thermal energy. It is only in those rare occasions that we take a different approach. Realize that $\tilde{u}^{n+1;LO}_{i,j,k}$ is still PCP, so it still has positive pressure. Therefore, we already have a pressure in the troubled zone that is positive, except that it was obtained via the update in eqn. (25). We can reset the troubled zone to have that same positive pressure as follows:-

$$\left[\mathbf{u}^{n+1;LO}_{i,j,k}\right]_5 = \left[\tilde{\mathbf{u}}^{n+1;LO}_{i,j,k}\right]_5 + \left[\mathbf{s}^{LO}_{i,j,k}\right]_5 \quad \text{with the definition}$$

$$\left[\mathbf{s}^{LO}_{i,j,k}\right]_5 \equiv \frac{1}{8\pi}\left\{\begin{array}{l}\left(\left(\overline{B}^{n+1;LO}_{x;i+1/2,j,k}+\overline{B}^{n+1;LO}_{x;i-1/2,j,k}\right)/2\right)^2+\left(\left(\overline{B}^{n+1;LO}_{y;i,j+1/2,k}+\overline{B}^{n+1;LO}_{y;i,j-1/2,k}\right)/2\right)^2\\ +\left(\left(\overline{B}^{n+1;LO}_{z;i,j,k+1/2}+\overline{B}^{n+1;LO}_{z;i,j,k-1/2}\right)/2\right)^2-\left(\left[\tilde{u}^{n+1;LO}_{i,j,k}\right]_6\right)^2-\left(\left[\tilde{u}^{n+1;LO}_{i,j,k}\right]_7\right)^2-\left(\left[\tilde{u}^{n+1;LO}_{i,j,k}\right]_8\right)^2\end{array}\right\}$$

(33)

The low order source term $\left[\mathbf{s}^{LO}_{i,j,k}\right]_5$ contributes only to the fifth component, which is the energy equation, for MHD. Only the 5$^{th}$ component of the source term vector $\mathbf{s}^{LO}_{i,j,k}$ is non-zero and that too only for zones that have lost the PCP property in the first order update. The really big advantage of this lower order formulation is that we always have a PCP formulation that is divergence-free which can be used to guide the higher order scheme so that it always remains divergence-free and PCP. In the next Sub-Section we will show how we do this in the most unobtrusive of ways so that for the most part we only use the high order scheme, only resorting to the lower order scheme from this Sub-Section in zones where the PCP property may be lost.

**5.2) Hybridizing the First Order PCP Scheme with the Higher Order Scheme**

In each zone $(i,j,k)$ we define a variable $\theta_{i,j,k}$. Following the philosophy of Hu, Adams, Shu [63], we will design a method such that when $\theta_{i,j,k}=1$ for all the zones, the scheme will exclusively be a high order scheme; this is the default when the MHD problem is not too stringent. When $\theta_{i,j,k}=0$ for all the zones, the scheme will exclusively be a first order scheme. Of course, it is not our intent/desire to have $\theta_{i,j,k}=0$ in any of the zones, but for stringent problems, we might have $0 \leq \theta_{i,j,k} < 1$ in some of the zones. At each zone boundary we can define a flux by

$$\mathbf{f}^{\theta}_{i+1/2,j,k} = \left(1-\theta^{face}_{i+1/2,j,k}\right)\tilde{\mathbf{f}}^{LO}_{i+1/2,j,k} + \theta^{face}_{i+1/2,j,k}\mathbf{f}^{HO}_{i+1/2,j,k} \quad \text{with} \quad \theta^{face}_{i+1/2,j,k} \equiv \min\left(\theta_{i,j,k},\theta_{i+1,j,k}\right) ;$$

$$\mathbf{g}^{\theta}_{i,j+1/2,k} = \left(1-\theta^{face}_{i,j+1/2,k}\right)\tilde{\mathbf{g}}^{LO}_{i,j+1/2,k} + \theta^{face}_{i,j+1/2,k}\mathbf{g}^{HO}_{i,j+1/2,k} \quad \text{with} \quad \theta^{face}_{i,j+1/2,k} \equiv \min\left(\theta_{i,j,k},\theta_{i,j+1,k}\right) ; \quad (34)$$

$$\mathbf{h}^{\theta}_{i,j,k+1/2} = \left(1-\theta^{face}_{i,j,k+1/2}\right)\tilde{\mathbf{h}}^{LO}_{i,j,k+1/2} + \theta^{face}_{i,j,k+1/2}\mathbf{h}^{HO}_{i,j,k+1/2} \quad \text{with} \quad \theta^{face}_{i,j,k+1/2} \equiv \min\left(\theta_{i,j,k},\theta_{i,j,k+1}\right)$$

The left panel of Fig. 8 shows how these facial values of $\theta$ are collocated. Likewise, at each edge we can define

$$\overline{E}^{\theta}_{z;i+1/2,j+1/2,k} = \left(1-\theta^{edge}_{i+1/2,j+1/2,k}\right)\overline{E}^{LO}_{z;i+1/2,j+1/2,k} + \theta^{edge}_{i+1/2,j+1/2,k}\overline{E}^{HO}_{z;i+1/2,j+1/2,k}$$
$$\text{with} \quad \theta^{edge}_{i+1/2,j+1/2,k} \equiv \min\left(\theta_{i,j,k},\theta_{i+1,j,k},\theta_{i,j+1,k},\theta_{i+1,j+1,k}\right) ;$$

$$\overline{E}^{\theta}_{y;i+1/2,j,k+1/2} = \left(1-\theta^{edge}_{i+1/2,j,k+1/2}\right)\overline{E}^{LO}_{y;i+1/2,j,k+1/2} + \theta^{edge}_{i+1/2,j,k+1/2}\overline{E}^{HO}_{y;i+1/2,j,k+1/2} \quad (35)$$
$$\text{with} \quad \theta^{edge}_{i+1/2,j,k+1/2} \equiv \min\left(\theta_{i,j,k},\theta_{i+1,j,k},\theta_{i,j,k+1},\theta_{i+1,j,k+1}\right) ;$$

$$\overline{E}^{\theta}_{x;i,j+1/2,k+1/2} = \left(1-\theta^{edge}_{i,j+1/2,k+1/2}\right)\overline{E}^{LO}_{x;i,j+1/2,k+1/2} + \theta^{edge}_{i,j+1/2,k+1/2}\overline{E}^{HO}_{x;i,j+1/2,k+1/2}$$
$$\text{with} \quad \theta^{edge}_{i,j+1/2,k+1/2} \equiv \min\left(\theta_{i,j,k},\theta_{i,j+1,k},\theta_{i,j,k+1},\theta_{i,j+1,k+1}\right)$$

The left panel of Fig. 9 shows how these edge values of $\theta$ are collocated. The right panel of Fig. 9 shows how these can be used to obtain electric fields at the edges. The above two equations show that the order of accuracy of the fluxes and electric fields can be locally lowered, as needed, to the point where the first order scheme always guarantees PCP behavior.

We now describe an iterative update strategy that does just that. We iterate over the whole mesh, starting the iteration process with $\theta_{i,j,k} = 1$ for all the zones. This $\theta_{i,j,k}$ will be sequentially lowered for any zone that is troubled; but realize that it will only be lowered for the few zones that are troubled. For each iteration, eqns. (2) to (5) can be modified to become

$$\partial_t \mathbf{u}^\theta_{i,j,k} = -\frac{1}{\Delta x}\left(\mathbf{f}^\theta_{i+1/2,j,k} - \mathbf{f}^\theta_{i-1/2,j,k}\right) - \frac{1}{\Delta y}\left(\mathbf{g}^\theta_{i,j+1/2,k} - \mathbf{g}^\theta_{i,j-1/2,k}\right) - \frac{1}{\Delta z}\left(\mathbf{h}^\theta_{i,j,k+1/2} - \mathbf{h}^\theta_{i,j,k-1/2}\right)$$
$$+ \frac{1}{\Delta t}\left(1 - \theta_{i,j,k}\right)\mathbf{s}^{LO}_{i,j,k} \qquad (36)$$

$$\partial_t \overline{B}^\theta_{x;\, i+1/2,j,k} = -\frac{1}{\Delta y \Delta z}\left(\Delta z \overline{E}^\theta_{z;\, i+1/2,j+1/2,k} - \Delta z \overline{E}^\theta_{z;\, i+1/2,j-1/2,k} + \Delta y \overline{E}^\theta_{y;\, i+1/2,j,k-1/2} - \Delta y \overline{E}^\theta_{y;\, i+1/2,j,k+1/2}\right) \quad (37)$$

$$\partial_t \overline{B}^\theta_{y;\, i,j-1/2,k} = -\frac{1}{\Delta x \Delta z}\left(\Delta x \overline{E}^\theta_{x;\, i,j-1/2,k+1/2} - \Delta x \overline{E}^\theta_{x;\, i,j-1/2,k-1/2} + \Delta z \overline{E}^\theta_{z;\, i-1/2,j-1/2,k} - \Delta z \overline{E}^\theta_{z;\, i+1/2,j-1/2,k}\right) \quad (38)$$

$$\partial_t \overline{B}^\theta_{z;\, i,j,k+1/2} = -\frac{1}{\Delta x \Delta y}\left(\Delta x \overline{E}^\theta_{x;\, i,j-1/2,k+1/2} - \Delta x \overline{E}^\theta_{x;\, i,j+1/2,k+1/2} + \Delta y \overline{E}^\theta_{y;\, i+1/2,j,k+1/2} - \Delta y \overline{E}^\theta_{y;\, i-1/2,j,k+1/2}\right) \quad (39)$$

Fig. 7 shows us schematically how our update strategy allows us to access any accuracy of update for the electric fields, going from the highest order accuracy to the lowest order accuracy. In this figure and the text that follows, the high order components are superscripted with "HO"; and low order components are superscripted with "LO". Using these, the forward Euler approximants can be constructed for any SSP-RK timestepping strategy. If a zone in those forward Euler approximants is not within the PCP domain, its local $\theta_{i,j,k}$ will be lowered even further in each successive iteration till it becomes PCP. Realize that this procedure is fully explicit; i.e. it does not require any fluxes or electric fields that are obtained through an implicit process. Realize too that this process is guaranteed to ensure that all zones are brought into the PCP domain.

### 6) Accuracy Analysis for Divergence Involution Constrained Hyperbolic PDE Systems

We focus on three PDE systems which have an involution constraint involving one or more vector fields. The first PDE system is drawn from computational electrodynamics (CED) which has two involution constraints associated with two vector fields. The vector fields are the electric flux density and the magnetic flux density. These two flux densities evolve in response to the curl of the magnetic field intensity and the electric field intensity respectively. CED is based on Maxwell's equations which are fundamental equations of nature and cannot be derived from any other equation. The constitutive relations that relate the field intensities to the flux densities are related to the permittivity and permeability of the material being considered. The second PDE system is drawn from classical magnetohydrodynamics (MHD) and the third PDE system is drawn from relativistic magnetohydrodynamics (RMHD). The latter two PDE systems have an involution constraint involving the divergence of the magnetic field. The magnetic field evolves in response to the curl of an electric field for the latter two PDE systems. The MHD and RMHD equations can be derived from a combination of the equations of fluid dynamics and Maxwell's equations. The electric field can, therefore, be written constitutively in terms of a cross product between the velocity field and the magnetic field.

For the time-update, we used a third-order SSP-RK scheme from Shu and Osher [78] and a five-stage fourth-order SSP-RK scheme from Spiteri and Ruuth [81]. For the fifth, seventh, and ninth order schemes, we had to reduce the timestep as the mesh was refined. The base level grid for all of the accuracy tests was run with a CFL of 0.3. For the spatially third-order accurate scheme, we use third order SSP-RK scheme, and for the fifth, seventh, and ninth-order accurate scheme, we use fourth-order SSP-RK time stepping. Consequently, for the fifth, seventh and ninth order scheme we reduce the time-step size as the mesh was refined so that the temporal error remain dominated by the spatial error. When a spatially fifth-order scheme is used with a temporally fourth-order accurate time-stepping strategy, then every doubling of the mesh requires a reduction in the timestep that goes as $\Delta t \rightarrow \Delta t \, (1/2)^{5/4}$. Similarly, when a spatially seventh-order scheme is used with a temporally fourth-order time-stepping strategy, then every doubling of the mesh requires a reduction in the timestep that goes as $\Delta t \rightarrow \Delta t \, (1/2)^{7/4}$. The same strategy generalizes to the ninth order.

**6.1) Accuracy Analysis for Computational Electrodynamics (CED)**

In this Sub-section, we perform an accuracy analysis for the CED system. The equations of CED can be written as two evolutionary curl-type equations for the electric displacement (extended Ampere's law) and the magnetic induction (Faraday's law).

$$\frac{\partial \mathbf{D}}{\partial t} - \nabla \times \mathbf{H} = -\mathbf{J} \quad ; \quad \frac{\partial \mathbf{B}}{\partial t} + \nabla \times \mathbf{E} = 0. \tag{40}$$

The electric displacement and the magnetic induction vector fields also satisfy the following two non-evolutionary constraint equations given by

$$\nabla \cdot \mathbf{D} = \rho_E \quad ; \quad \nabla \cdot \mathbf{B} = 0 \tag{41}$$

Here $\rho_E$ is the electric charge density and "**J**" is the current density. In material media, the electric displacement vector is also related to the electric field vector, and the magnetic induction vector is related to the magnetic field vector. The constitutive relations are given by

$$\mathbf{D} = \epsilon \mathbf{E} \quad ; \quad \mathbf{B} = \mu \mathbf{H} \tag{42}$$

where, in general, $\epsilon$ is a symmetric 3 × 3 permittivity tensor and $\boldsymbol{\mu}$ is a symmetric 3 × 3 magnetic permeability tensor that depends on material properties.

We perform an accuracy study for the above system using a three-dimensional plane-polarized electromagnetic wave that has been formulated in Balsara *et al.* [29]. This test problem consists of a plane-polarized electromagnetic wave propagating in a vacuum along the diagonal of a three-dimensional Cartesian mesh spanning $[-0.5, 0.5]^3$. Periodic boundary conditions are enforced. The details setup is given in Balsara *et al.* (2016) and therefore we do not provide the details here. Table 1 provides the accuracy analysis for the $B_y$ variable. We see that all the presented schemes reach their design accuracies efficiently.

**Table 1 shows the accuracy of the 3D Electromagnetic wave problem for the CED system using the divergence-preserving AFD-WENO schemes presented here; the $B_y$ variable is shown.**

| Order 2 | L$_1$ Error | L$_1$ Accuracy | L$_{inf}$ Error | L$_{inf}$ Accuracy |
|---|---|---|---|---|
| 16$^3$ | 7.37622E-02 |  | 1.35489E-01 |  |
| 32$^3$ | 2.02438E-02 | 1.87 | 4.21736E-02 | 1.68 |
| 64$^3$ | 6.18236E-03 | 1.71 | 1.44576E-02 | 1.54 |
| $128^3$ | 1.68161E-03 | 1.88 | 5.87110E-03 | 1.30 |
| Order 3 | | | | |
| $16^3$ | 2.27828E-02 | | 3.54511E-02 | |
| $32^3$ | 3.02372E-03 | 2.91 | 4.72780E-03 | 2.91 |
| $64^3$ | 3.81555E-04 | 2.99 | 5.99318E-04 | 2.98 |
| $128^3$ | 4.78018E-05 | 3.00 | 7.50724E-05 | 3.00 |
| Order 5 | | | | |
| $16^3$ | 9.70335E-04 | | 1.53615E-03 | |
| $32^3$ | 3.17519E-05 | 4.93 | 4.96451E-05 | 4.95 |
| $64^3$ | 1.00029E-06 | 4.99 | 1.57239E-06 | 4.98 |
| $128^3$ | 3.13601E-08 | 5.00 | 4.92623E-08 | 5.00 |
| Order 7 | | | | |
| $16^3$ | 4.05552E-05 | | 6.42364E-05 | |
| $32^3$ | 2.71903E-07 | 7.22 | 4.26464E-07 | 7.23 |
| $48^3$ | 1.60355E-08 | 6.98 | 2.51501E-08 | 6.98 |
| $64^3$ | 2.14420E-09 | 6.99 | 3.37040E-09 | 6.99 |
| Order 9 | | | | |
| $16^3$ | 1.88289E-05 | | 2.94354E-05 | |
| $32^3$ | 3.37518E-09 | 12.45 | 5.27892E-09 | 12.45 |
| $48^3$ | 6.33200E-11 | 9.81 | 9.94403E-11 | 9.80 |
| $64^3$ | 3.48335E-12 | 10.08 | 5.50682E-12 | 10.06 |

## 6.2) Accuracy Analysis for Magnetohydrodynamics (MHD)

In this Sub-section we consider the two-dimensional Magnetohydrodynamics (MHD) system to study the accuracy of the presented schemes. The MHD system for an ideal fluid can be written in a conservation form as

$$\frac{\partial}{\partial t}\begin{pmatrix}\rho \\ \rho v_x \\ \rho v_y \\ \rho v_z \\ \varepsilon \\ B_x \\ B_y \\ B_z\end{pmatrix} + \frac{\partial}{\partial x}\begin{pmatrix}\rho v_x \\ \rho v_x^2 + p + \mathbf{B}^2/8\pi - B_x^2/4\pi \\ \rho v_x v_y - B_x B_y/4\pi \\ \rho v_x v_z - B_x B_z/4\pi \\ (\varepsilon + p + \mathbf{B}^2/8\pi)v_x - B_x(\mathbf{v}\cdot\mathbf{B})/4\pi \\ 0 \\ (v_x B_y - v_y B_x) \\ -(v_z B_x - v_x B_z)\end{pmatrix} + \frac{\partial}{\partial y}\begin{pmatrix}\rho v_y \\ \rho v_y v_x - B_y B_x/4\pi \\ \rho v_y^2 + p + \mathbf{B}^2/8\pi - B_y^2/4\pi \\ \rho v_y v_z - B_y B_z/4\pi \\ (\varepsilon + p + \mathbf{B}^2/8\pi)v_y - B_y(\mathbf{v}\cdot\mathbf{B})/4\pi \\ -(v_x B_y - v_y B_x) \\ 0 \\ (v_y B_z - v_z B_y)\end{pmatrix} = 0$$

(43)

where $\varepsilon = \rho \mathbf{v}^2/2 + p/(\gamma-1) + \mathbf{B}^2/8\pi$ is the total energy density. The density is denoted by $\rho$; the pressure is written as "$p$"; the velocity components are given by $v_x$, $v_y$, $v_z$; and the magnetic field components by $B_x$, $B_y$, $B_z$. "$\gamma$" denotes the ratio of specific heats. The usual challenge in MHD simulations has been that when the velocities or magnetic fields become too large, the pressure can become zero or negative.

In Balsara [12] we presented a genuinely two-dimensional vortex problem for the above MHD system. The problem consists of a smoothly varying and dynamically stable vortex that moves diagonally. In Balsara [12] the explicit expressions for the set up have been provided, therefore we do not describe the setup details here. In Table 2 we show the accuracy analysis for the $B_y$ variable. To minimize the effect of small jumps in the velocity field at the periodic boundaries, we double the computational domain and stopping time for the fifth, seventh, and ninth-order schemes. We observe that the presented schemes are able to reach their design accuracy.

**Table 2 shows the accuracy of the 2D Magnetohydrodynamics vortex problem for the divergence-preserving AFD-WENO schemes presented here; the $B_y$ variable is shown.**

| Order 2 | L₁ Error | L₁ Accuracy | L_inf Error | L_inf Accuracy |
|---|---|---|---|---|
| $32^2$ | 1.38248E-02 |  | 2.14786E-01 |  |
| $64^2$ | 3.82861E-03 | 1.85 | 6.56915E-02 | 1.71 |
| $128^2$ | 9.49995E-04 | 2.01 | 1.74273E-02 | 1.91 |
| $256^2$ | 2.61112E-04 | 1.86 | 4.82106E-03 | 1.85 |
| Order 3 |  |  |  |  |
| $32^2$ | 5.39616E-03 |  | 7.64776E-02 |  |
| $64^2$ | 8.63601E-04 | 2.64 | 1.14915E-02 | 2.73 |
| $128^2$ | 1.11966E-04 | 2.95 | 1.42516E-03 | 3.01 |
| $256^2$ | 1.41120E-05 | 2.99 | 1.77298E-04 | 3.01 |
| Order 5 |  |  |  |  |
| $32^2$ | 3.69666E-03 |  | 2.16003E-01 |  |
| $64^2$ | 4.52850E-04 | 3.03 | 2.65982E-02 | 3.02 |
| $128^2$ | 1.78924E-05 | 4.66 | 1.12614E-03 | 4.56 |
| $256^2$ | 5.76284E-07 | 4.96 | 3.61606E-05 | 4.96 |
| Order 7 |  |  |  |  |
| $32^2$ | 1.97898E-03 |  | 1.09211E-01 |  |
| $64^2$ | 7.34885E-05 | 4.75 | 4.46131E-03 | 4.61 |
| $128^2$ | 7.19814E-07 | 6.67 | 4.64895E-05 | 6.58 |
| $256^2$ | 5.89200E-09 | 6.93 | 3.76116E-07 | 6.95 |

| Order 9 | | | | |
|---|---|---|---|---|
| $32^2$ | 2.12589E-03 | | 1.29615E-01 | |
| $64^2$ | 1.38312E-05 | 7.26 | 8.22383E-04 | 7.30 |
| $128^2$ | 3.66605E-08 | 8.56 | 2.44937E-06 | 8.39 |
| $256^2$ | 7.70028E-11 | 8.90 | 5.06964E-09 | 8.92 |

### 6.3) Accuracy Analysis for Relativistic Magnetohydrodynamics (RMHD)

In this Sub-section we consider the two-dimensional Relativistic-Magnetohydrodynamics (RMHD) system to study the accuracy of the presented schemes. The RMHD system in two dimensions can be written in a conservation form as

$$\frac{\partial}{\partial t}\begin{pmatrix} \rho W \\ m_x \\ m_y \\ m_z \\ \varepsilon \\ B_x \\ B_y \\ B_z \end{pmatrix} + \frac{\partial}{\partial x}\begin{pmatrix} \rho W v_x \\ m_x v_x + p + \frac{|b|^2}{2} - B_x b_x / W \\ m_y v_x - B_x b_y / W \\ m_z v_x - B_x b_z / W \\ m_x \\ 0 \\ (v_x B_y - v_y B_x) \\ -(v_z B_x - v_x B_z) \end{pmatrix} + \frac{\partial}{\partial y}\begin{pmatrix} \rho W v_y \\ m_x v_y - B_y b_x / W \\ m_y v_y + p + \frac{|b|^2}{2} - B_y b_y / W \\ m_z v_y - B_y b_z / W \\ m_y \\ -(v_x B_y - v_y B_x) \\ 0 \\ (v_y B_z - v_z B_y) \end{pmatrix} = 0$$

(44)

Here, $\rho$ is the density, $W$ is the Lorentz factor, set $\{\varepsilon, m_x, m_y, m_z\}$ forms the four-momentum density, $\mathbf{v} = (v_x, v_y, v_z)^\top$ is the velocity vector, $p$ is the pressure and $B_x, B_y, B_z$ are the components of the magnetic field. The speed of light is assumed to be unity; therefore, the Lorentz factor is defined by $W = 1/\sqrt{1-\mathbf{v}^2}$. The conserved quantities $m_i$ and $\varepsilon$ are given by

$$m_i = (\rho h W^2 + \mathbf{B}^2)v_i - (\mathbf{v}\cdot\mathbf{B})B_i \quad ; \quad \varepsilon = \rho h W^2 - p + \frac{\mathbf{B}^2}{2} + \frac{\mathbf{v}^2 \mathbf{B}^2 - (\mathbf{v}\cdot\mathbf{B})^2}{2}$$

(45)

where $h = 1 + \frac{\gamma}{\gamma-1}\frac{p}{\rho}$ is the specific enthalpy and $\gamma$ denotes the polytropic index. In the above, Latin indices range from 1 to 3. The components and magnitude of the covariant magnetic field are defined as

$$b_\mu = W\left(\mathbf{v}\cdot\mathbf{B}, \frac{B_i}{W^2} + (\mathbf{v}\cdot\mathbf{B})v_i\right) \quad ; \quad |b|^2 = \frac{\mathbf{B}^2}{W^2} + (\mathbf{v}\cdot\mathbf{B})^2 \tag{46}$$

where the Greek indices range from 0 to 3. Cai *et al.* [44] have presented a provably convergent and robust Newton-Raphson method to extract the primitives from the conservative variables. We use the same strategy here. We define the vector of conservative variables as $\mathbf{U} = (\rho W, m_x, m_y, m_z, \varepsilon, B_x, B_y, B_z)^\top$. Balsara and Kim [26] suggested that the better way to reconstruct the state variables (given in vector $\mathbf{U}$) is to reconstruct in the special vector of primitive variables given by $\mathbf{V} = (\rho, v_x W, v_y W, v_z W, p, B_x, B_y, B_z)^\top$. We adopt the same idea and reconstruct $\mathbf{V}$ to obtain the reconstructed primitives and then use the change of variable matrix $\frac{\partial \mathbf{U}}{\partial \mathbf{V}}$ to obtain the derivative terms of the form $\frac{\partial \mathbf{U}}{\partial x}$ using the formula $\frac{\partial \mathbf{U}}{\partial x} \cong \frac{\partial \mathbf{U}}{\partial \mathbf{V}}\frac{\partial \mathbf{V}}{\partial x}$.

Balsara and Kim [26] have presented a genuinely two-dimensional vortex problem for the RMHD system. The problem consists of a smoothly varying and dynamically stable vortex that moves diagonally in the square computational domain. The detailed derivation and setup expression are given in Balsara and Kim [26]; therefore, we do not describe the setup details here. In Table 3 we show the accuracy analysis for the $B_y$ variable. To minimize the effect of small jumps in the velocity field at the periodic boundaries, we double the computational domain and stopping time for the seventh, and ninth-order schemes. From Table 3 we see that the presented schemes are able to reach their design accuracy for the RMHD system.

**Table 3 shows the accuracy of the 2D Relativistic-Magnetohydrodynamics vortex problem for the divergence-preserving AFD-WENO schemes presented here; the $B_y$ variable is shown.**

| Order 3 | L$_1$ Error | L$_1$ Accuracy | L$_{inf}$ Error | L$_{inf}$ Accuracy |
|---|---|---|---|---|
| 64$^2$ | 4.04260E-03 | | 1.28301E-01 | |
| 128$^2$ | 8.49900E-04 | 2.25 | 3.81457E-02 | 1.75 |
| $256^2$ | 1.23860E-04 | 2.78 | 7.30376E-03 | 2.38 |
| $512^2$ | 1.63154E-05 | 2.92 | 1.04766E-03 | 2.80 |
| Order 5 | | | | |
| $64^2$ | 7.01418E-04 | | 3.05950E-02 | |
| $128^2$ | 5.30257E-05 | 3.73 | 3.34089E-03 | 3.19 |
| $256^2$ | 2.11697E-06 | 4.65 | 1.88630E-04 | 4.15 |
| $512^2$ | 9.57948E-08 | 4.47 | 6.46521E-06 | 4.87 |
| Order 7 | | | | |
| $64^2$ | 1.79695E-03 | | 2.30325E-01 | |
| $128^2$ | 2.57514E-04 | 2.80 | 3.10841E-02 | 2.89 |
| $256^2$ | 3.68427E-06 | 6.13 | 9.95726E-04 | 4.96 |
| $512^2$ | 4.75228E-08 | 6.28 | 1.65906E-05 | 5.91 |
| Order 9 | | | | |
| $64^2$ | 2.82170E-03 | | 2.26298E-01 | |
| $128^2$ | 3.10540E-04 | 3.18 | 3.40935E-02 | 2.73 |
| $256^2$ | 1.55905E-06 | 7.64 | 3.14428E-04 | 6.76 |
| $512^2$ | 3.71425E-09 | 8.71 | 1.52784E-06 | 7.69 |

**7) Test problems for Hyperbolic PDE Systems with a Divergence-Based Involution Constraint**

We present several two-dimensional test problems for the CED, MHD and RMHD systems. For each of the problems, we take a CFL number of 0.4 unless stated otherwise.

**7.1) Test Problems Involving CED**

In this sub-section, we present a set of two test problems for the CED system. The first problem shows the refraction of a compact electromagnetic beam by a dielectric slab. The problem was described in detail in Balsara *et al.* [29], therefore we do not repeat the description here. The problem is set up on a rectangular $xy$ – domain that spans $[-5,8]\times[-2.5,7]\mu m$. For the simulation shown, we use a fifth-order accurate scheme using a 1300×900 zone mesh. The permittivity increases in a tapered fashion from $\epsilon_0$ for $x<0$ to 2.25 $\epsilon_0$ for $x>0$. A compact electromagnetic beam is incident on the slab at an angle of incidence given by $45^o$. Figs. 10a, 10b, 10c show, respectively, $B_z$, $D_x$, and $D_y$ at the initial time $t=0$. Figs. 10d, 10e, and 10f show the same profiles at the final time of $4.0\times10^{-14}$ s. The surface of the dielectric slab is shown by a vertical dashed black line. The inclined dashed black lines show the angles of incidence, refraction, and reflection, and these black lines are over-plotted on the field components to guide the eye. According to Snell's law, the angle of refraction should be $28.12^o$ since the refractive index of the dielectric slab is 1.5. We see that the code obtains the correct angle of refraction. We also observe that some of the radiation is reflected from the surface of the slab. The presence of a reflected wave is consistent with the Fresnel conditions for transmission and reflection of radiation at dielectric interfaces.

The second problem shows the total internal reflection of a compact electromagnetic beam by a dielectric slab. This problem was also described in detail in Balsara *et al.* [29], therefore we do not repeat the description here. This problem is set up on a rectangular $xy$ – domain that spans $[-6,1]\times[-2.5,6]\mu m$. We use a seventh-order accurate scheme using a 700×850 zone mesh. Here, $\epsilon$ is chosen such that it has a value of 4.0 $\epsilon_0$ for $x \leq 0$ and tapers rapidly to the ambient value of $\epsilon_0$ for $x > 0$. This value of permittivity for $x < 0$ implies a refractive index of 2 for the dielectric slab. For this mesh, the taper width that is applied to the variation in the permittivity is 0.25 times a zone width. The wavelength in the dielectric medium corresponds to about 30 zones. For such a slab, the critical angle for total internal reflection is $30^o$. The angle of incidence of the incident radiation is $45^o$, with the result that the incident radiation will undergo total internal reflection. Figs. 11a, 11b, and 11c show, respectively, $B_z$, $D_x$, and $D_y$ at the time $t = 0$. The same profiles at a final time of $t = 5\times10^{-14}$ are shown in Figs. 11d, 11e, and 11f. The surface of the dielectric slab is shown by a vertical black line. The inclined black lines for the incident and reflected rays are over-plotted on the field components to guide our eye. We see from the figures that the radiation undergoes total internal reflection, as expected. The results are consistent with the results reported in Balsara *et al.* [26], showing the occurrence of total internal reflection.

**7.2) Test Problems Involving MHD**

In this sub-section, we present various two-dimensional test problems for the MHD system. We begin by examining the Orszag-Tang problem that was first introduced in Orszag and Tang [71]. The Orszag-Tang problem is widely used as a test model for the emergence of MHD turbulence. The problem is initialized on a square domain that spans [0,2]×[0,2]. The ghost cells are filled with the periodic boundaries. The adiabatic constant is set as $\gamma = 5/3$. The density was initialized to $\gamma^2$ all over the domain. The pressure is uniformly set to $\gamma$. The velocity field is set up as follows:

$$\mathbf{v} = -\sin(\pi y)\hat{x} + \sin(\pi x)\hat{y}$$

The magnetic field is initialized in a divergence-free manner using the following vector potential:

$$A_z = -\frac{\sqrt{4\pi}}{2\pi}\left(\cos(2\pi x) + 2\cos(\pi y)\right)$$

We run simulations until a final time of $t = 1$ on a grid of 256×256 zones. We use the ninth-order scheme, and show results for the density variable, pressure variable, magnitude of the velocity vector and magnitude of the magnetic field vector in Fig. 11. We observe a close resemblance between our obtained results and those presented in Balsara [12]. The third, fifth and seventh order schemes also perform well on this problem, therefore they are not shown here.

The second test problem is the two-dimensional rotor problem that was first presented in Balsara & Spicer [8]; see also Balsara [12]. The computational domain spans the domain [-0.5, 0.5]×[-0.5,0.5]. A dense and rapidly spinning cylinder is set up in the center of an initially stationary, light ambient fluid. The ambient fluid is initially at rest. A uniform magnetic field initially threads the two fluids. Its value is set to 2.5 units, and it initially points in the x-direction. The total pressure in the fluid is set to unity, i.e. $p=1$. The density in the ambient fluid is uniformly set to unity, while the constant density in the rotor is 10 units out to a radius of 0.1. A linear taper is applied to the density between a radius of 0.1 and 0.13 so that the density in the rotor decreases linearly to the value of the density in the ambient fluid. Six zones are used for the taper to join the density of the two fluids. That number should be kept fixed if the resolution is increased or decreased. The initial angular velocity of the rotor is uniform out to a radius of 0.1. At this radius, the toroidal velocity has a value of one unit. The toroidal velocity decreases linearly from one unit to zero between a radius of 0.1 and 0.13 so that it joins the velocity of the ambient fluid at a radius of 0.13. The simulations were run until a final time of t=0.29 using a mesh consisting of 256×256 zones. We use the seventh-order accurate scheme for this test problem. Fig. 13 shows the results for the density variable, pressure variable, magnitude of the velocity vector and magnitude of the magnetic field vector. We see that the results in Fig. 13 closely match the profiles presented in Balsara & Spicer [8] and Balsara [12]. The third, fifth and ninth order schemes also perform well on this problem, therefore they are not shown here.

Next, we consider two variations of the blast problem (BLAST-I and BLAST-II). The first type is described in Balsara & Spicer [8], while the second case is sourced from Wu & Shu [85]. Both types are examples of stringent problems for the MHD system because they correspond to very low values of plasma β. As a result, the problems require a Physically Constraint Preserving (PCP) variant of the presented schemes (see Bhoriya *et al* [41] and Section 5 for the detailed strategy). The first problem is run on a three-dimensional grid with a domain $[-0.5, 0.5] \times [-0.5, 0.5] \times [-0.5, 0.5]$ using a $300^3$ zone mesh. The second problem was configured on a unit square grid that consists of 400×400 zones and spans the domain [-0.5,0.5]×[-0.5,0.5]. Initially, the density for both types is uniformly set to unity and the velocity is set to zero. The two problems differ in terms of magnetic and pressure strengths. For the first problem type (BLAST-I), the pressure is uniformly set to 0.1 except within a central spherical region of radius 0.1, where it is elevated to 1000. Additionally, a magnetic field with a magnitude of 100 is initialized along the *x*-direction that makes the problem challenging for numerical schemes. While this problem was originally presented in two-dimensions in Balsara & Spicer [8], here we present the three-dimensional variant of the same problem. We simulate the problem until a final time of *t* = 0.01 using the fifth-order scheme and show results for the density profile, pressure profile, magnitude of the velocity vector and magnitude of the magnetic field vector in Fig. 14. There is a good similarity between the obtained results, and those presented in Balsara & Spicer [8]. For the second problem type (BLAST-II), the pressure within a central circle of radius 0.1 is elevated to 10000, making initial jump larger in the pressure variable. Additionally, a much stronger magnetic field with a magnitude of 1000 is initialized along the *x*-direction that makes the problem extremely challenging for numerical schemes. Fig. 15 displays the numerical results at time *t*=0.001 obtained by the fifth-order scheme. The obtained results are consistent with those reported in Wu & Shu [85], highlighting the effectiveness of the PCP property of the presented schemes.

Next, we consider the astrophysical jet problem (Mach 800) from Balsara [19]. Following Wu & Shu [85], a magnetic field is added to this problem to simulate the MHD jet flows. The presence of a magnetic field makes this test more extreme. As a result, the problems require a Physically Constraint Preserving (PCP) variant of the presented schemes (see Bhoriya *et al* [41] and Section 5 for the detailed strategy). We consider three variations of the jet problem based on the strength of the magnetic field ($B_y$- component). All variants of the problem are configured on a 2D grid consisting of 400×600 zones and spans the domain [-0.5,0.5]×[0.0,1.5]. Initially, the domain is filled with a static medium with uniform density of $0.1\gamma$ and uniform pressure of unity ($p=1$), where the adiabatic index $\gamma$ is set to 1.4. The magnetic field is given by $(B_x, B_y, B_z) = (0, B_a, 0)$. At the bottom boundary $(y=0)$, a dense jet is injected through the inlet part $|x|<0.05$ with the states $(\rho, v_x, v_y, v_z, p, B_x, B_y, B_z) = (\gamma, 0, 800, 0, 1, 0, B_a, 0)$. Outflow boundaries are employed at all the other boundaries. We consider the below three test cases (based on the value of $B_a$)

Jet-I: Moderately magnetized case: $B_a = \sqrt{200}$ (corresponding plasma-beta, $\beta_a = 10^{-2}$)

Jet-II: Strongly magnetized case: $B_a = \sqrt{2000}$ (corresponding $\beta_a = 10^{-3}$)

Jet-III: Extremely strongly magnetized case: $B_a = \sqrt{20000}$ (corresponding $\beta_a = 10^{-4}$)

All three variants were run to a stopping time of $t=0.002$. Fig. 16a, 16b, and 16c show the resulting density profiles on logarithmic scales for the Jet-I, Jet-II and Jet-III problems, respectively. Fig. 17a, 17b, and 17c show the resulting pressure profiles on logarithmic scales for the Jet-I, Jet-II and Jet-III problems, respectively. A seventh-order accurate scheme is used for all the runs. From Figs. 16 and 17, we observe that the cocoon heads in the beam are well captured for these extreme cases; hence showing the robust performance of the proposed methods.

**7.3) Test Problems Involving RMHD**

In this sub-section, we present various two-dimensional test problems for the RMHD system. We begin by examining the blast problem for the RMHD system. The non-relativistic version of this test problem was first presented in Balsara and Spicer [8] and was extended to RMHD system by Komissarov [66]. Several variants of this test problem have been presented by Mignone and Bodo [70] by using a slightly modified version of the moderately magnetized case in Komissarov [66]. In this subsection, we use the setup given in Mignone and Bodo [70]; see also Balsara and Kim [26].

The blast test problem is set up on a two-dimensional square domain that spans $[-6,6]\times[-6,6]$. Within a radius of 0.8, the explosion zone has a density of $10^{-2}$ and a pressure of 1. Outside a radius of 1 unit, the ambient medium has a density of $10^{-4}$ and a pressure of $5\times10^{-4}$. A linear taper is applied to the density and pressure from a radius of 0.8 to 1. Accordingly, both the density and pressure linearly decrease with increasing radius in that range of radii. The magnetic field is initialized in the x-direction and has a magnitude of 0.1. A polytropic index of $\gamma = 4/3$ is used in this problem. The simulation is run to a final time of 4 using a fifth-order accurate scheme on a grid consisting of $512\times512$ zones. Obtained results have been presented in Fig. 18. Figure 18a shows the logarithm of the density, Fig. 18b shows the logarithm of the pressure, and Fig. 18c shows the magnetic pressure at time t=4. We see a close resemblance of the obtained results with the results presented in Mignone and Bodo [70] and Balsara and Kim [26].

Next, we consider the Orszag-Tang problem that was first introduced in Orszag and Tang [71] for the non-relativistic magnetohydrodynamics and later was extended to relativistic magnetohydrodynamics in van der Holst *et. al* [60]; also see Wu and Shu [86]. Our setup is same as in Wu and Shu [86]. The problem is initialized on a square domain that spans $[0, 2\pi]\times[0, 2\pi]$. The ghost cells are filled with the periodic boundaries. The adiabatic constant is set as $\gamma = 4/3$. The density is initialized to 1 all over the domain. The pressure is uniformly set to 10. The velocity field is set up as follows:

$$\mathbf{v} = -A\sin(y)\hat{x} + A\sin(x)\hat{y}$$

where $A = 0.99/\sqrt{2}$. The magnetic field is initialized in a divergence-free manner using the following vector potential:

$$A_z = -\left(\frac{\cos(2x)}{2} + \cos(y)\right)$$

We run two simulations, one until a final time of $t = 2.818127$ and another until a time of $t=6.8558$. We use a grid of 512×512 zones for both the runs. We use the seventh order scheme and show the obtained results in Fig. 19. Figures 19a and 19b show the logarithm of the density and the logarithm of the magnetic pressure at time t=2.818127; and Figures 19c and 19d show the logarithm of the density and the logarithm of the magnetic pressure at time t=6.8558. For both the runs, we observe a close resemblance between our obtained results and those presented in Wu and Shu [86]. The third, fifth and ninth order schemes also perform well on this problem, therefore they are not shown here.

Next, we consider the Shock-cloud interaction problem from Wu and Shu [86]. The problem describes the dynamics of a high-density cloud when it is hit by a strong shock wave. The problem is initialized on a two-dimensional domain that spans $[-0.2,1.2]\times[0,1]$. We use the inflow boundary conditions at the left boundary and outflow boundary conditions at the remaining boundaries. The inflow boundary states are same as the left shock states. A right moving shock at $x=0.05$ is initialized using the following states

$$(\rho, \mathbf{v}, \mathbf{B}, p) = \begin{cases} (3.86859, 0.68,\ 0,\ 0,\ 0,\ 0.84981,\ -0.84981,\ 1.25115) & \text{if } x < 0.05 \\ (1,\ 0,\ 0,\ 0,\ 0,\ 0.16106,\ 0.16106,\ 0.05) & \text{otherwise} \end{cases}$$

A cloud centered at (0.25,0.5) with radius 0.15 is initialized in the domain. The cloud has the same states to the ambient fluid except for a higher density of 30 units. The adiabatic constant is set as $\gamma = 5/3$. The simulation is run to a final time of 1.2 using a ninth-order accurate scheme on a grid that consists of $560\times 400$ zones. Results have been presented in Fig. 20. Fig. 20a shows the obtained logarithm of the density and Fig. 20b shows the logarithm of the magnetic pressure at time t=1.2. We see a close resemblance of the obtained results with the results presented in Wu and Shu [86]. The third, fifth and seventh order schemes also perform well on this problem, therefore they are not shown here.

**7.4) Speed Comparisons Between Finite Difference and Finite Volume WENO Schemes**

It is always interesting to ask whether there is a substantial speed advantage in the finite difference WENO schemes (like the one presented here) and the finite volume WENO schemes. An example of such a finite volume WENO scheme would be the one for MHD in Balsara *et al*. [39]. Since that is a recently published work, we take that as a point of reference and catalogue the speed differences within the context of numerical MHD. When considering speed, it is also important to recognize that finite difference (FD) methods do not suffer much from the curse of dimensionality, whereas finite volume (FV) methods do suffer a lot from the curse of dimensionality. As a result, a 2D FD scheme may have a certain speed advantage over a 2D FV scheme; however, a 3D FD scheme might have an even greater speed advantage over a 3D FV scheme. Even so, a fair discussion should not just stop there. We have to realize that a FV method may use a tensor product array of calls to the Riemann solver at the 2D boundary of a 3D zone; however, it is possible to use methods from van der Vegt, and van der Ven [83] [84] to reduce the number of calls to the Riemann solver; thereby leading to some efficiencies. Analogous FD methods do not have the ability to exploit those efficiencies. Similarly, FV methods can make efficient use of ADER methods (Dumbser *et al*. [54], [55], Balsara *et al*. [14], [15]) to improve the efficiency of the timestepping. FD methods, because they do not have a volumetric reconstruction, cannot indeed draw on ADER approaches for timestepping. In particular, Balsara *et al*. [15] found that ADER methods offer a factor of two speed improvement relative to SSP-RK methods.

We mention all these factors in the previous paragraph because every direct comparison between speeds of FD and FV codes of the same order of accuracy is always very nuanced. Firstly, as explained before, it depends on whether the problem is 2D or 3D. Secondly, it depends on whether efficient tricks were used for flux evaluation or not. Thirdly, it depends on whether ADER timestepping was used or not. Fourthly, it is important to appreciate that the stencil size increases very dramatically with order and dimensionality for FV WENO schemes, whereas this is not the case for FD WENO schemes. The comparison we present here is between an FD WENO scheme and a FV WENO scheme that does not use the efficiencies of van der Vegt, and van der Ven [83] [84] to reduce the number of calls to the Riemann solver. Furthermore, to keep the comparison fair, we use the same SSP-RK timestepping for both FD and FV schemes; but the FV scheme could have used ADER methods to improve its speed by at least a factor of two. We present both 2D speed comparisons and 3D speed comparisons so that the effects of dimensionality can be factored in. Table 4 shows the results. We see that FD WENO schemes do not suffer much from the curse of dimensionality as we go from 2D to 3D. We also see that in 2D, the FD methods have a 7 to 14 fold speed advantage over FV methods and that advantage becomes more pronounced at higher order. In 3D, the FD methods have a 5 to 12 fold speed advantage over FV methods and, as before, that advantage becomes more pronounced at higher order. We also see that as we go from 2D to 3D, FD methods do not have as pronounced a degradation in speed as FV methods. In making Table 4 we used identical algorithmic elements between FD and FV methods just to keep the comparison fair. However, it should also be noted in a spirit of fairness that FV methods do have access to tricks, as shown in the previous paragraph, that would enable them to close the gap.

**Table 4 was made on a single core of an Intel(R) Xeon(R) Gold 6248R CPU @ 3.00GHz with a GNU Fortran compiler (gcc v9.4). It shows the speeds, as measured in zones updated per second, by FD WENO schemes and FV WENO schemes at 3$^{rd}$ and 5$^{th}$ orders. 2D problems used meshes with 400$^2$ zones, 3D problems used meshes with 128$^3$ zones.**

| Dimension (zones) | Order | FD WENO Zones/sec | FV WENO Zones/sec | Speed ratio |
|---|---|---|---|---|
| 2D (400 × 400) | 3$^{rd}$ | 40000 | 6038 | 6.6 |
|  | 5$^{th}$ | 18824 | 1280 | 14.7 |
| 3D (128 × 128 × 128) | 3$^{rd}$ | 22075 | 4671 | 4.7 |
|  | 5$^{th}$ | 16777 | 1391 | 12.1 |

## 8) Conclusions

AFD-WENO schemes have seen increasing prominence in recent years (Balsara *et al*. [36], [37]) because they can handle hyperbolic systems in conservation form as well as hyperbolic systems with non-conservative products. Moreover, they do so while preserving conservation when parts or all of the PDE system is in conservation form. They have also been extended to retain a physical constraint preserving (PCP) property, Bhoriya *et al*. [41]. Well-balancing can also be achieved within the context of AFD-WENO schemes, as shown by Xu and Shu [87]. The free stream preserving property can also be guaranteed for AFD-WENO schemes, as shown by Jiang *et al*. [65]. It is, therefore, natural to ask for AFD-WENO schemes that preserve the discrete divergence for PDE systems that have an involution property that preserves the divergence. Such involution preserving PDE systems encompass some very important application areas such as computational electrodynamics, magnetohydrodynamics and relativistic magnetohydrodynamics.

By studying large classes of such PDEs, we find that for every vector field whose divergence has to be preserved on the mesh, we have an update equation that looks like eqn. (1). Such an update equation forms a sub-system of the overall PDE and a divergence-free or divergence-preserving update for such a system inevitably requires a Yee-style collocation of variables of the form shown in Fig. 1. If multiple vector fields have a divergence-preserving update, as is the case for computational electrodynamics, then the PDE system will have multiple copies of equations that look like eqn. (1), and all of those vector fields will have collocations of the type shown in Fig. 1. This brings us to the realization that a common solution strategy should be found for all these classes of involution-constrained PDEs. In all such cases, we find that stabilization of the numerical PDE requires a two-dimensional upwinding, as shown in Figs. 3 and 4. The associated math, which is based on the multidimensional Riemann solver, was initially developed in Balsara [17], [18], [22]. In Section 3, we specialize that math for PDEs that have a sub-system of the type shown in eqn. (1). In that same Section, eqns. (12) and (20) give explicit formulae for such a multidimensional LLF and HLL Riemann solvers for involution-constrained PDEs that have a divergence-preserving constraint. Section 4 then provides a step-by-step plan for implementing the divergence-preserving AFD-WENO scheme. Section 5 shows the extra steps that have to be taken to achieve an efficient and time-explicit scheme that has the PCP property for divergence-preserving PDEs.

Section 6 shows the versatility of our method because it shows that several rather different PDE systems with a divergence-based involution constraint can all be treated using the same formalism. For all these different involution-constrained PDE systems we show that the method meets its design accuracy. Section 7 shows several stringent test problems associated with the same PDE systems. Some of these test problems are so stringent that their simulation requires enforcement of a PCP condition. When solving such problems, we are introduced to another advantage of AFD-WENO schemes which is that it is just as convenient (and acceptable) to interpolate the primitive variables as it is to interpolate the conserved variables. This advantage is not shared by the classical FD-WENO method, where the right- and left-going fluxes can only be reconstructed in a dimension-by-dimension finite volume sense. Especially for relativistic MHD we find that only the primitive variables first documented in Balsara and Kim [26] will discriminate between flow velocities that are very close to the speed of light. (Of course, the primitive variables are indeed projected into the space of eigenvectors that are developed for primitive variables when carrying out 1D interpolation.) As a result, Section 5 of this paper, which describes PCP methods for involution-constrained PDEs, is very valuable when designing AFD-WENO schemes for very stringent MHD and relativistic MHD flows.

**Acknowledgements**

DSB acknowledges support via NSF grants NSF-AST-2009776, NSF-AST-2434532, NASA grant NASA-2020-1241 and NASA grant 80NSSC22K0628. CWS acknowledges support via NSF grant DMS-2309249 and NSF-AST-2434532.

**Ethical Statement**
**i. Compliance with Ethical Standards** : This manuscript complies with all ethical standards for scientific publishing.
**ii. (in case of Funding) Funding** : The funding has been acknowledged. DSB acknowledges support via NSF grants NSF-AST-2009776, NSF-AST-2434532, NASA grant NASA-2020-1241 and NASA grant 80NSSC22K0628. CWS acknowledges support via NSF grant DMS-2309249 and NSF-AST-2434532.
**iii. Conflict of Interest** : On behalf of all authors, the corresponding author states that there is no conflict of interest.
**iv. Ethical approval** : N/A
**v. Informed consent** : N/A

**vi. Data Statement** : All data that was used in the generation of the figures has been stored and available for later use.

## Appendix A: Efficient Finite-Difference 2D WENO-AO interpolations.

We describe a very efficient finite-difference WENO-AO interpolation in 2D at third, fifth and seventh orders. While the math for ninth order gives rise to very large expressions that are too long to document in a paper, the reader who understands this procedure can use a computer algebra system to easily make the extension to ninth and higher orders. We use this 2D WENO interpolation to obtain electric field at the edge-centered locations. We use the tensor product of following 1D Legendre polynomials, that span the interval $[-1/2, 1/2]$, for the basis of interpolated polynomial in 2D.

$$L_0(x) = 1 \; ; \; L_1(x) = x \; ; \; L_2(x) = x^2 - \frac{1}{12} \; ; \; L_3(x) = x^3 - \frac{3}{20}x \; ;$$

$$L_4(x) = x^4 - \frac{3}{14}x^2 + \frac{3}{560} \; ; \; L_5(x) = x^5 - \frac{5}{18}x^3 + \frac{5}{336}x \; ;$$

$$L_6(x) = x^6 - \frac{15}{44}x^4 + \frac{5}{176}x^2 - \frac{5}{14784} \; ; \quad (A.1)$$

$$L_7(x) = x^7 - \frac{21}{52}x^5 + \frac{105}{2288}x^3 - \frac{35}{27456}x \; ;$$

$$L_8(x) = x^8 - \frac{7}{15}x^6 + \frac{7}{104}x^4 - \frac{7}{2288}x^2 + \frac{7}{329472}.$$

We let "$r$" denote the order of accuracy of the interpolation; for example, an interpolation that is only based on $L_0(x) = 1$, $L_1(x) = x$ and $L_2(x) = x^2 - 1/12$ corresponds to $r = 3$.

### A.1) r=3 Finite-Difference WENO Interpolation in 2D

The 2D interpolation polynomial at third order is given by

$$u(x,y) = u_{00} + u_x L_1(x) + u_y L_1(y) + u_{xx} L_2(x) + u_{yy} L_2(y) + u_{xy} L_1(x) L_1(y). \quad (A.2)$$

Fig. 5 shows us five possible stencils that can each be used to evaluate the $u_{00}, u_x, u_y, u_{xx}, u_{yy}$ and $u_{xy}$ terms. The stencil $S_1^{r3}$ denotes the Right-Up biased third-order accurate stencil, $S_2^{r3}$ denotes the Left-Up biased third-order accurate stencil, $S_3^{r3}$ denotes the Left-Down biased third-order accurate stencil, $S_4^{r3}$ denotes the Right-Down biased third-order accurate stencil and $S_5^{r3}$ denotes the central third-order accurate stencil. For each stencil, we catalogue expressions for the unknown terms $u_{00}, u_x, u_y, u_{xx}, u_{yy}$ and $u_{xy}$ below.

For the stencil $S_1^{r_3}$ we obtain

$$u_{00} = \left(-2u_{0,1} + u_{0,2} + 26u_{0,0} - 2u_{1,0} + u_{2,0}\right)/24$$
$$u_x = \left(-3u_{0,0} + 4u_{1,0} - u_{2,0}\right)/2$$
$$u_y = \left(-3u_{0,0} + 4u_{0,1} - u_{0,2}\right)/2$$
$$u_{xx} = \left(u_{0,0} - 2u_{1,0} + u_{2,0}\right)/2 \qquad (A.3)$$
$$u_{yy} = \left(u_{0,0} - 2u_{0,1} + u_{0,2}\right)/2$$
$$u_{xy} = u_{1,1} + u_{0,0} - u_{0,1} - u_{1,0}$$

For the stencil $S_2^{r_3}$ we obtain

$$u_{00} = \left(-2u_{0,1} + u_{0,2} + 26u_{0,0} - 2u_{-1,0} + u_{-2,0}\right)/24$$
$$u_x = \left(3u_{0,0} - 4u_{-1,0} + u_{-2,0}\right)/2$$
$$u_y = \left(-3u_{0,0} + 4u_{0,1} - u_{0,2}\right)/2$$
$$u_{xx} = \left(u_{0,0} - 2u_{-1,0} + u_{-2,0}\right)/2 \qquad (A.4)$$
$$u_{yy} = \left(u_{0,0} - 2u_{0,1} + u_{0,2}\right)/2$$
$$u_{xy} = u_{0,1} + u_{-1,0} - u_{-1,1} - u_{0,0}$$

For the stencil $S_3^{r_3}$ we obtain

$$u_{00} = \left(-2u_{0,-1} + u_{0,-2} + 26u_{0,0} - 2u_{-1,0} + u_{-2,0}\right)/24$$
$$u_x = \left(3u_{0,0} - 4u_{-1,0} + u_{-2,0}\right)/2$$
$$u_y = \left(3u_{0,0} - 4u_{0,-1} + u_{0,-2}\right)/2$$
$$u_{xx} = \left(u_{0,0} - 2u_{-1,0} + u_{-2,0}\right)/2 \qquad (A.5)$$
$$u_{yy} = \left(u_{0,0} - 2u_{0,-1} + u_{0,-2}\right)/2$$
$$u_{xy} = -u_{0,-1} - u_{-1,0} + u_{-1,-1} + u_{0,0}$$

For the stencil $S_4^{r_3}$ we obtain

$$u_{00} = \left(-2u_{0,-1} + u_{0,-2} + 26u_{0,0} - 2u_{1,0} + u_{2,0}\right)/24$$
$$u_x = \left(-3u_{0,0} + 4u_{1,0} - u_{2,0}\right)/2$$
$$u_y = \left(3u_{0,0} - 4u_{0,-1} + u_{0,-2}\right)/2$$
$$u_{xx} = \left(u_{0,0} - 2u_{1,0} + u_{2,0}\right)/2 \tag{A.6}$$
$$u_{yy} = \left(u_{0,0} - 2u_{0,-1} + u_{0,-2}\right)/2$$
$$u_{xy} = u_{0,-1} - u_{0,0} + u_{1,0} - u_{1,-1}$$

For the stencil $S_5^{r_3}$ we obtain

$$u_{00} = \left(u_{0,-1} + u_{0,1} + u_{-1,0} + u_{1,0} + 20u_{0,0}\right)/24$$
$$u_x = \left(u_{1,0} - u_{-1,0}\right)/2$$
$$u_y = \left(u_{0,1} - u_{0,-1}\right)/2$$
$$u_{xx} = \left(u_{1,0} - 2u_{0,0} + u_{-1,0}\right)/2 \tag{A.7}$$
$$u_{yy} = \left(u_{0,1} - 2u_{0,0} + u_{0,-1}\right)/2$$
$$u_{xy} = u_{1,1} + u_{-1,-1} - u_{-1,1} - u_{1,-1}$$

The smoothness indicator (denoted by $\beta^{r_3}$) for the interpolated polynomial in eqn. (A.1), written as a sum of perfect squares, is given by

$$\beta^{r_3} = u_x^2 + u_y^2 + \frac{13}{3}u_{xx}^2 + \frac{13}{3}u_{yy}^2 + \frac{7}{6}u_{xy}^2.$$

The smoothness indicators for the stencil $S_1^{r_3}$ is denoted by $\beta_1^{r_3}$, for the stencil $S_2^{r_3}$ is denoted by $\beta_2^{r_3}$, for the stencil $S_3^{r_3}$ is denoted by $\beta_3^{r_3}$, for the stencil $S_4^{r_3}$ is denoted by $\beta_4^{r_3}$ and for the stencil $S_5^{r_3}$ is denoted by $\beta_5^{r_3}$.

In designing a WENO scheme at third order, we wish to make a non-linearly hybridized interpolation using stencils $S_1^{r_3}, S_2^{r_3}, S_3^{r_3}, S_4^{r_3}$ and $S_5^{r_3}$. The interpolation is described by one parameter $\gamma_{Lo} \in (0,1)$. We set $\gamma_{Lo} = 0.85$. The linear weights for the stencils $S_1^{r_3}, S_2^{r_3}, S_3^{r_3}, S_4^{r_3}$ are given by

$$\gamma_1^{r_3} = (1-\gamma_{Lo})/4 \quad ; \quad \gamma_2^{r_3} = (1-\gamma_{Lo})/4 \quad ; \quad \gamma_3^{r_3} = (1-\gamma_{Lo})/4 \quad ; \quad \gamma_4^{r_3} = (1-\gamma_{Lo})/4. \tag{A.8}$$

The linear weight for the stencil $S_5^{r_3}$ is given by $\gamma_5^{r_3} = \gamma_{Lo}$. Notice that for the linear weights we have, $\gamma_1^{r_3} + \gamma_2^{r_3} + \gamma_3^{r_3} + \gamma_4^{r_3} + \gamma_5^{r_3} = 1$. We see that when the central stencil $S_5^{r_3}$ is smooth we want most of our interpolation to come from the central stencil because it is the most stable choice when the solution is smooth. However, when a suitable comparison of the smoothness indicators shows that the central stencil is non-smooth, we want most (or all) of our interpolation to be weighted towards either the Right-Up biased stencil or the Left-Up biased stencil or the Left-Down biased stencil or the Right-Down biased stencil, based on the smoothest of the interpolated polynomial.

To avoid loss of accuracy at critical points (Borges *et al.*, 2008) we use the smoothness indicators to define

$$\tau = \frac{1}{4}\left(\left|\beta_5^{r_3} - \beta_1^{r_3}\right| + \left|\beta_5^{r_3} - \beta_2^{r_3}\right| + \left|\beta_5^{r_3} - \beta_3^{r_3}\right| + \left|\beta_5^{r_3} - \beta_4^{r_3}\right|\right) . \tag{A.9}$$

The unnormalized non-linear weights are given by

$$w_1^{r_3} = \gamma_1^{r_3}\left(1 + \tau^2 \big/ \left(\beta_1^{r_3} + \varepsilon\right)^2\right) \; ; \; w_2^{r_3} = \gamma_2^{r_3}\left(1 + \tau^2 \big/ \left(\beta_2^{r_3} + \varepsilon\right)^2\right) \; ;$$
$$w_3^{r_3} = \gamma_3^{r_3}\left(1 + \tau^2 \big/ \left(\beta_3^{r_3} + \varepsilon\right)^2\right) \; ; \; w_4^{r_3} = \gamma_4^{r_3}\left(1 + \tau^2 \big/ \left(\beta_4^{r_3} + \varepsilon\right)^2\right) \; ; \tag{A.10}$$
$$w_5^{r_3} = \gamma_5^{r_3}\left(1 + \tau^2 \big/ \left(\beta_5^{r_3} + \varepsilon\right)^2\right)$$

where $\epsilon = 10^{-12}$ is a small number to avoid the divisions by zero. The normalization of the non-linear weights is given by

$$\overline{w}_1^{r_3} = w_1^{r_3} \big/ \left(w_1^{r_3} + w_2^{r_3} + w_3^{r_3} + w_4^{r_3} + w_5^{r_3}\right) \; ; \; \overline{w}_2^{r_3} = w_2^{r_3} \big/ \left(w_1^{r_3} + w_2^{r_3} + w_3^{r_3} + w_4^{r_3} + w_5^{r_3}\right) \; ;$$
$$\overline{w}_3^{r_3} = w_3^{r_3} \big/ \left(w_1^{r_3} + w_2^{r_3} + w_3^{r_3} + w_4^{r_3} + w_5^{r_3}\right) \; ; \; \overline{w}_4^{r_3} = w_4^{r_3} \big/ \left(w_1^{r_3} + w_2^{r_3} + w_3^{r_3} + w_4^{r_3} + w_5^{r_3}\right) \; ; \tag{A.11}$$
$$\overline{w}_5^{r_3} = w_5^{r_3} \big/ \left(w_1^{r_3} + w_2^{r_3} + w_3^{r_3} + w_4^{r_3} + w_5^{r_3}\right).$$

Let $P_i^{r_3}(x,y)$ be the interpolated polynomial corresponding to stencil $S_i^{r_3}$. The non-linearly hybridized third order accurate WENO interpolation is given by

$$P^{(3)}(x,y) = \overline{w}_1^{r_3} P_1^{r_3}(x,y) + \overline{w}_2^{r_3} P_2^{r_3}(x,y) + \overline{w}_3^{r_3} P_3^{r_3}(x,y) + \overline{w}_4^{r_3} P_4^{r_3}(x,y) + \overline{w}_5^{r_3} P_5^{r_3}(x,y) .$$
$$\tag{A.12}$$

This completes our description of the third order 2D WENO interpolation.

## A.2) WENO-AO(5,3) Interpolation in 2D

The 2D WENO-AO(5,3) interpolation consists of a non-linear hybridization between a large, centered, fifth order stencil denoted by $S_1^{r5}$ and the five smaller stencils described in Appendix A.1). The larger central fifth order accurate stencil is shown in Fig. 6a. The fifth order accurate 2D polynomial is given

$$
\begin{aligned}
u(x,y) = &\, u_{00} + u_x L_1(x) + u_y L_1(y) + u_{xx} L_2(x) + u_{yy} L_2(y) + u_{xy} L_1(x) L_1(y) \\
&+ u_{xxx} L_3(x) + u_{yyy} L_3(y) + u_{xxy} L_2(x) L_1(y) + u_{xyy} L_1(x) L_2(y) \\
&+ u_{xxxx} L_4(x) + u_{yyyy} L_4(y) + u_{xxxy} L_3(x) L_1(y) + u_{xyyy} L_1(x) L_3(y) + u_{xxyy} L_2(x) L_2(y)
\end{aligned}
\tag{A.13}
$$

The large central fifth order accurate stencil gives

$$u_{00} = (4636 u_{0,0} + 288 u_{0,-1} - 17 u_{0,-2} + 288 u_{0,1} - 17 u_{0,2} + 288 u_{-1,0} + 10 u_{-1,-1} + 10 u_{-1,1}$$
$$- 17 u_{-2,0} + 288 u_{1,0} + 10 u_{1,-1} + 10 u_{1,1} - 17 u_{2,0}) / 5760$$

$$u_x = (-144 u_{-1,0} - 5 u_{-1,-1} - 5 u_{-1,1} + 17 u_{-2,0} + 144 u_{1,0} + 5 u_{1,-1} + 5 u_{1,1} - 17 u_{2,0}) / 240$$

$$u_y = (-144 u_{0,-1} - 5 u_{-1,-1} - 5 u_{1,-1} + 17 u_{0,-2} + 144 u_{0,1} + 5 u_{-1,1} + 5 u_{1,1} - 17 u_{0,2}) / 240$$

$$u_{xx} = (-374 u_{0,0} - 14 u_{0,-1} - 14 u_{0,1} + 198 u_{-1,0} + 7 u_{-1,-1} + 7 u_{-1,1} - 11 u_{-2,0} + 198 u_{1,0}$$
$$+ 7 u_{1,-1} + 7 u_{1,1} - 11 u_{2,0}) / 336$$

$$u_{yy} = (-374 u_{0,0} - 14 u_{-1,0} - 14 u_{1,0} + 198 u_{0,-1} + 7 u_{-1,-1} + 7 u_{1,-1} - 11 u_{0,-2} + 198 u_{0,1}$$
$$+ 7 u_{-1,1} + 7 u_{1,1} - 11 u_{0,2}) / 336$$

$$u_{xy} = (188 u_{-1,-1} - 17 u_{-1,-2} - 188 u_{-1,1} + 17 u_{-1,2} - 17 u_{-2,-1} + 17 u_{-2,1} - 188 u_{1,-1} + 17 u_{1,-2}$$
$$+ 188 u_{1,1} - 17 u_{1,2} + 17 u_{2,-1} - 17 u_{2,1}) / 480$$

$$u_{xxx} = (2 u_{-1,0} - u_{-2,0} - 2 u_{1,0} + u_{2,0}) / 12$$

$$u_{yyy} = (2 u_{0,-1} - u_{0,-2} - 2 u_{0,1} + u_{0,2}) / 12$$

$$u_{xxy} = (2 u_{0,-1} - 2 u_{0,1} - u_{-1,-1} + u_{-1,1} - u_{1,-1} + u_{1,1}) / 4$$

$$u_{xyy} = (2 u_{-1,0} - u_{-1,-1} - u_{-1,1} - 2 u_{1,0} + u_{1,-1} + u_{1,1}) / 4$$

$$u_{xxxx} = (6 u_{0,0} - 4 u_{-1,0} + u_{-2,0} - 4 u_{1,0} + u_{2,0}) / 24$$

$$u_{yyyy} = (6 u_{0,0} - 4 u_{0,-1} + u_{0,-2} - 4 u_{0,1} + u_{0,2}) / 24$$

$$u_{xxxy} = (-2 u_{-1,-1} + 2 u_{-1,1} + u_{-2,-1} - u_{-2,1} + 2 u_{1,-1} - 2 u_{1,1} - u_{2,-1} + u_{2,1}) / 24$$

$$u_{xyyy} = (-2 u_{-1,-1} + u_{-1,-2} + 2 u_{-1,1} - u_{-1,2} + 2 u_{1,-1} - u_{1,-2} - 2 u_{1,1} + u_{1,2}) / 24$$

$$u_{xxyy} = (4 u_{0,0} - 2 u_{0,-1} - 2 u_{0,1} - 2 u_{-1,0} + u_{-1,-1} + u_{-1,1} - 2 u_{1,0} + u_{1,-1} + u_{1,1}) / 4$$

The smoothness indicator (denoted by $\beta^{r5}$) for the interpolated polynomial in eqn. (A.13), written as a sum of perfect squares, is given by

$$\beta^{r_5} = \left(u_x + \frac{u_{xxx}}{10}\right)^2 + \frac{13}{3}\left(u_{xx} + \frac{123}{455}u_{xxxx}\right)^2 + \frac{781}{20}u_{xxx}^2 + \frac{1421461}{2275}u_{xxxx}^2$$
$$+ \left(u_y + \frac{u_{yyy}}{10}\right)^2 + \frac{13}{3}\left(u_{yy} + \frac{123}{455}u_{yyyy}\right)^2 + \frac{781}{20}u_{yyy}^2 + \frac{1421461}{2275}u_{yyyy}^2$$
$$+ \frac{7}{6}\left(u_{xy} + \frac{13}{140}u_{xxxy} + \frac{13}{140}u_{xyyy}\right)^2 + \frac{47}{10}u_{xxy}^2 + \frac{47}{10}u_{xyy}^2$$
$$+ \frac{88841}{2100}u_{xxxy}^2 + \frac{88841}{2100}u_{xyyy}^2 + \frac{1}{16800}\left(u_{xxxy} - u_{xyyy}\right)^2 + \frac{5083}{270}u_{xxyy}^2$$

In designing a WENO scheme at fifth order, we wish to make a non-linearly hybridized interpolation using the third order stencils $S_1^{r_3}, S_2^{r_3}, S_3^{r_3}, S_4^{r_3}, S_5^{r_3}$ and a fifth order stencil $S_1^{r_5}$. The interpolation is described by two parameters $\gamma_{Lo}, \gamma_{Hi} \in (0,1)$. We set $\gamma_{Lo} = \gamma_{Hi} = 0.85$. The linear weights for the stencils $S_1^{r_3}, S_2^{r_3}, S_3^{r_3}, S_4^{r_3}$ are given by

$$\gamma_1^{r_3} = (1-\gamma_{Hi})(1-\gamma_{Lo})/4 \quad ; \quad \gamma_2^{r_3} = (1-\gamma_{Hi})(1-\gamma_{Lo})/4 \quad ;$$
$$\gamma_3^{r_3} = (1-\gamma_{Hi})(1-\gamma_{Lo})/4 \quad ; \quad \gamma_4^{r_3} = (1-\gamma_{Hi})(1-\gamma_{Lo})/4 \quad .$$

The linear weight for the stencil $S_5^{r_3}$ is given by $\gamma_5^{r_3} = (1-\gamma_{Hi})\gamma_{Lo}$. The linear weight for the stencil $S_1^{r_5}$ is given by $\gamma_1^{r_5} = \gamma_{Hi}$. Notice that for the linear weights we have, $\gamma_1^{r_3} + \gamma_2^{r_3} + \gamma_3^{r_3} + \gamma_4^{r_3} + \gamma_5^{r_3} + \gamma_1^{r_5} = 1$. We see that when the fifth order central stencil $S_1^{r_5}$ is smooth we want most of our interpolation to come from the central stencil because it is the most stable choice when the solution is smooth. However, when a suitable comparison of the smoothness indicators shows that the fifth order central stencil is non-smooth, we want most (or all) of our interpolation to be weighted towards the third order stencils. We denote the fifth order central smoothness indicator by $\beta_1^{r_5}$.

To avoid loss of accuracy at critical points (Borges et al., 2008) we use the smoothness indicators to define

$$\tau = \frac{1}{5}\left(\left|\beta_1^{r_5} - \beta_1^{r_3}\right| + \left|\beta_1^{r_5} - \beta_2^{r_3}\right| + \left|\beta_1^{r_5} - \beta_3^{r_3}\right| + \left|\beta_1^{r_5} - \beta_4^{r_3}\right| + \left|\beta_1^{r_5} - \beta_5^{r_3}\right|\right).$$

The unnormalized non-linear weights are given by

$$w_1^{r3} = \gamma_1^{r3}\left(1+\tau^2/\left(\beta_1^{r3}+\varepsilon\right)^2\right) \quad ; \quad w_2^{r3} = \gamma_2^{r3}\left(1+\tau^2/\left(\beta_2^{r3}+\varepsilon\right)^2\right) \quad ;$$

$$w_3^{r3} = \gamma_3^{r3}\left(1+\tau^2/\left(\beta_3^{r3}+\varepsilon\right)^2\right) \quad ; \quad w_4^{r3} = \gamma_4^{r3}\left(1+\tau^2/\left(\beta_4^{r3}+\varepsilon\right)^2\right) \quad ;$$

$$w_5^{r3} = \gamma_5^{r3}\left(1+\tau^2/\left(\beta_5^{r3}+\varepsilon\right)^2\right) \quad : \quad w_1^{r5} = \gamma_1^{r5}\left(1+\tau^2/\left(\beta_1^{r5}+\varepsilon\right)^2\right)$$

where $\epsilon = 10^{-12}$ is a small number to avoid the divisions by zero. The normalization of the non-linear weights is given by

$$\overline{w}_1^{r3} = w_1^{r3}/\left(w_1^{r3}+w_2^{r3}+w_3^{r3}+w_4^{r3}+w_5^{r3}+w_1^{r5}\right) \quad ; \quad \overline{w}_2^{r3} = w_2^{r3}/\left(w_1^{r3}+w_2^{r3}+w_3^{r3}+w_4^{r3}+w_5^{r3}+w_1^{r5}\right) \quad ;$$

$$\overline{w}_3^{r3} = w_3^{r3}/\left(w_1^{r3}+w_2^{r3}+w_3^{r3}+w_4^{r3}+w_5^{r3}+w_1^{r5}\right) \quad ; \quad \overline{w}_4^{r3} = w_4^{r3}/\left(w_1^{r3}+w_2^{r3}+w_3^{r3}+w_4^{r3}+w_5^{r3}+w_1^{r5}\right) \quad ;$$

$$\overline{w}_5^{r3} = w_5^{r3}/\left(w_1^{r3}+w_2^{r3}+w_3^{r3}+w_4^{r3}+w_5^{r3}+w_1^{r5}\right) \quad ; \quad \overline{w}_1^{r5} = w_1^{r5}/\left(w_1^{r3}+w_2^{r3}+w_3^{r3}+w_4^{r3}+w_5^{r3}+w_1^{r5}\right).$$

Let $P_i^{r3}(x,y)$ be the interpolated polynomial corresponding to stencil $S_i^{r3}$ and $P_1^{r5}(x,y)$ be the interpolated polynomial corresponding to stencil $S_1^{r5}$. We denote the interpolated polynomial for WENO-AO(5,3) as $P^{AO(5,3)}(x,y)$. Our task in this paragraph is to describe the construction of the order-preserving, non-linearly hybridized, fifth order polynomial $P^{AO(5,3)}(x,y)$. The non-linear weights should be combined in such a way that when all the smoothness indicators seem to have almost similar values then only the higher order scheme is obtained. The non-linearly hybridized fifth order accurate interpolation is given by

$$P^{AO(5,3)}(x) = \frac{\overline{w}_1^{r5}}{\gamma_1^{r5}}\begin{pmatrix} P_1^{r5}(x,y) - \gamma_1^{r3}P_1^{r3}(x,y) - \gamma_2^{r3}P_2^{r3}(x,y) - \gamma_3^{r3}P_3^{r3}(x,y) \\ -\gamma_4^{r3}P_4^{r3}(x,y) - \gamma_5^{r3}P_5^{r3}(x,y) \end{pmatrix}$$
$$+ \overline{w}_1^{r3}P_1^{r3}(x,y) + \overline{w}_2^{r3}P_2^{r3}(x,y) + \overline{w}_3^{r3}P_3^{r3}(x,y) + \overline{w}_4^{r3}P_4^{r3}(x,y) + \overline{w}_5^{r3}P_5^{r3}(x,y)$$

This completes our description of the fifth order 2D WENO-AO(5,3) interpolation.

### A.3) WENO-AO(7,5,3) Interpolation in 2D

The 2D WENO-AO(7,5,3) interpolation consists of a non-linear hybridization between a large, centered, seventh order stencil denoted by $S_1^{r7}$, an intermediate central fifth order stencil denoted by $S_1^{r5}$ and the five smaller stencils described in Appendix A.1). The larger central seventh order accurate stencil is shown in Fig. 6b and The intermediate fifth order central stencil is shown in Fig. 6a. The seventh order accurate 2D polynomial is given

$$\begin{aligned}u(x,y) &= u_{00} + u_x L_1(x) + u_y L_1(y) + u_{xx} L_2(x) + u_{yy} L_2(y) + u_{xy} L_1(x)L_1(y) \\ &+ u_{xxx} L_3(x) + u_{yyy} L_3(y) + u_{xxy} L_2(x)L_1(y) + u_{xyy} L_1(x)L_2(y) \\ &+ u_{xxxx} L_4(x) + u_{yyyy} L_4(y) + u_{xxxy} L_3(x)L_1(y) + u_{xyyy} L_1(x)L_3(y) + u_{xxyy} L_2(x)L_2(y) \\ &+ u_{xxxxx} L_5(x) + u_{yyyyy} L_5(y) + u_{xxxxy} L_4(x)L_1(y) + u_{xyyyy} L_1(x)L_4(y) \\ &+ u_{xxxyy} L_3(x)L_2(y) + u_{xxyyy} L_2(x)L_3(y) \\ &+ u_{xxxxxx} L_6(x) + u_{yyyyyy} L_6(y) + u_{xxxxxy} L_5(x)L_1(y) + u_{xyyyyy} L_1(x)L_5(y) \\ &+ u_{xxxxyy} L_4(x)L_2(y) + u_{xxyyyy} L_2(x)L_4(y) + u_{xxxyyy} L_3(x)L_3(y)\end{aligned} \quad (A.14)$$

The large central stencil gives

$$\begin{aligned}u_{00} =& (767024 u_{0,0} + 52223 u_{0,-1} - 4820 u_{0,-2} + 367 u_{0,-3} + 52223 u_{0,1} - 4820 u_{0,2} + 367 u_{0,3} \\ &+ 52223 u_{-1,0} + 2632 u_{-1,-1} - 119 u_{-1,-2} + 2632 u_{-1,1} - 119 u_{-1,2} - 4820 u_{-2,0} - 119 u_{-2,-1} \\ &- 119 u_{-2,1} + 367 u_{-3,0} + 52223 u_{1,0} + 2632 u_{1,-1} - 119 u_{1,-2} + 2632 u_{1,1} - 119 u_{1,2} - 4820 u_{2,0} \\ &- 119 u_{2,-1} - 119 u_{2,1} + 367 u_{3,0}) / 967680 \\ u_x =& (-52223 u_{-1,0} - 2632 u_{-1,-1} + 119 u_{-1,-2} - 2632 u_{-1,1} + 119 u_{-1,2} + 9640 u_{-2,0} + 238 u_{-2,-1} \\ &+ 238 u_{-2,1} - 1101 u_{-3,0} + 52223 u_{1,0} + 2632 u_{1,-1} - 119 u_{1,-2} + 2632 u_{1,1} - 119 u_{1,2} - 9640 u_{2,0} \\ &- 238 u_{2,-1} - 238 u_{2,1} + 1101 u_{3,0}) / 80640 \\ u_y =& (-52223 u_{0,-1} + 9640 u_{0,-2} - 1101 u_{0,-3} + 52223 u_{0,1} - 9640 u_{0,2} + 1101 u_{0,3} - 2632 u_{-1,-1} \\ &+ 238 u_{-1,-2} + 2632 u_{-1,1} - 238 u_{-1,2} + 119 u_{-2,-1} - 119 u_{-2,1} - 2632 u_{1,-1} + 238 u_{1,-2} \\ &+ 2632 u_{1,1} - 238 u_{1,2} + 119 u_{2,-1} - 119 u_{2,1}) / 80640\end{aligned}$$

$$u_{xx} = (-93672u_{0,0} - 4972u_{0,-1} + 238u_{0,-2} - 4972u_{0,1} + 238u_{0,2} + 50921u_{-1,0} + 2596u_{-1,-1}$$
$$- 119u_{-1,-2} + 2596u_{-1,1} - 119u_{-1,2} - 4418u_{-2,0} - 110u_{-2,-1} - 110u_{-2,1} + 333u_{-3,0}$$
$$+ 50921u_{1,0} + 2596u_{1,-1} - 119u_{1,-2} + 2596u_{1,1} - 119u_{1,2} - 4418u_{2,0} - 110u_{2,-1} - 110u_{2,1}$$
$$+ 333u_{3,0})/80640$$

$$u_{yy} = (-93672u_{0,0} + 50921u_{0,-1} - 4418u_{0,-2} + 333u_{0,-3} + 50921u_{0,1} - 4418u_{0,2} + 333u_{0,3}$$
$$- 4972u_{-1,0} + 2596u_{-1,-1} - 110u_{-1,-2} + 2596u_{-1,1} - 110u_{-1,2} + 238u_{-2,0} - 119u_{-2,-1}$$
$$- 119u_{-2,1} - 4972u_{1,0} + 2596u_{1,-1} - 110u_{1,-2} + 2596u_{1,1} - 110u_{1,2} + 238u_{2,0} - 119u_{2,-1}$$
$$- 119u_{2,1})/80640$$

$$u_{xy} = (387074u_{-1,-1} - 58672u_{-1,-2} + 5505u_{-1,-3} - 387074u_{-1,1} + 58672u_{-1,2} - 5505u_{-1,3}$$
$$- 58672u_{-2,-1} + 4046u_{-2,-2} + 58672u_{-2,1} - 4046u_{-2,2} + 5505u_{-3,-1} - 5505u_{-3,1}$$
$$- 387074u_{1,-1} + 58672u_{1,-2} - 5505u_{1,-3} + 387074u_{1,1} - 58672u_{1,2} + 5505u_{1,3}$$
$$+ 58672u_{2,-1} - 4046u_{2,-2} - 58672u_{2,1} + 4046u_{2,2} - 5505u_{3,-1} + 5505u_{3,1})/806400$$

$$u_{xxx} = (217u_{-1,0} + 6u_{-1,-1} + 6u_{-1,1} - 134u_{-2,0} - 3u_{-2,-1} - 3u_{-2,1} + 17u_{-3,0} - 217u_{1,0} - 6u_{1,-1}$$
$$- 6u_{1,1} + 134u_{2,0} + 3u_{2,-1} + 3u_{2,1} - 17u_{3,0})/864$$

$$u_{yyy} = (217u_{0,-1} - 134u_{0,-2} + 17u_{0,-3} - 217u_{0,1} + 134u_{0,2} - 17u_{0,3} + 6u_{-1,-1} - 3u_{-1,-2} - 6u_{-1,1}$$
$$+ 3u_{-1,2} + 6u_{1,-1} - 3u_{1,-2} - 6u_{1,1} + 3u_{1,2})/864$$

$$u_{xxy} = (2486u_{0,-1} - 238u_{0,-2} - 2486u_{0,1} + 238u_{0,2} - 1298u_{-1,-1} + 119u_{-1,-2} + 1298u_{-1,1}$$
$$- 119u_{-1,2} + 55u_{-2,-1} - 55u_{-2,1} - 1298u_{1,-1} + 119u_{1,-2} + 1298u_{1,1} - 119u_{1,2} + 55u_{2,-1}$$
$$- 55u_{2,1})/3360$$

$$u_{xyy} = (2486u_{-1,0} - 1298u_{-1,-1} + 55u_{-1,-2} - 1298u_{-1,1} + 55u_{-1,2} - 238u_{-2,0} + 119u_{-2,-1}$$
$$+ 119u_{-2,1} - 2486u_{1,0} + 1298u_{1,-1} - 55u_{1,-2} + 1298u_{1,1} - 55u_{1,2} + 238u_{2,0} - 119u_{2,-1}$$
$$- 119u_{2,1})/3360$$

$$u_{xxxx} = (2272u_{0,0} + 66u_{0,-1} + 66u_{0,1} - 1583u_{-1,0} - 44u_{-1,-1} - 44u_{-1,1} + 488u_{-2,0} + 11u_{-2,-1}$$
$$+ 11u_{-2,1} - 41u_{-3,0} - 1583u_{1,0} - 44u_{1,-1} - 44u_{1,1} + 488u_{2,0} + 11u_{2,-1} + 11u_{2,1} - 41u_{3,0}) / 6336$$

$$u_{yyyy} = (2272u_{0,0} - 1583u_{0,-1} + 488u_{0,-2} - 41u_{0,-3} - 1583u_{0,1} + 488u_{0,2} - 41u_{0,3} + 66u_{-1,0}$$
$$- 44u_{-1,-1} + 11u_{-1,-2} - 44u_{-1,1} + 11u_{-1,2} + 66u_{1,0} - 44u_{1,-1} + 11u_{1,-2} - 44u_{1,1} + 11u_{1,2}) / 6336$$

$$u_{xxxy} = (-1349u_{-1,-1} + 102u_{-1,-2} + 1349u_{-1,1} - 102u_{-1,2} + 802u_{-2,-1} - 51u_{-2,-2} - 802u_{-2,1}$$
$$+ 51u_{-2,2} - 85u_{-3,-1} + 85u_{-3,1} + 1349u_{1,-1} - 102u_{1,-2} - 1349u_{1,1} + 102u_{1,2} - 802u_{2,-1}$$
$$+ 51u_{2,-2} + 802u_{2,1} - 51u_{2,2} + 85u_{3,-1} - 85u_{3,1}) / 8640$$

$$u_{xyyy} = (-1349u_{-1,-1} + 802u_{-1,-2} - 85u_{-1,-3} + 1349u_{-1,1} - 802u_{-1,2} + 85u_{-1,3} + 102u_{-2,-1} - 51u_{-2,-2}$$
$$- 102u_{-2,1} + 51u_{-2,2} + 1349u_{1,-1} - 802u_{1,-2} + 85u_{1,-3} - 1349u_{1,1} + 802u_{1,2} - 85u_{1,3} - 102u_{2,-1}$$
$$+ 51u_{2,-2} + 102u_{2,1} - 51u_{2,2}) / 8640$$

$$u_{xxyy} = (936u_{0,0} - 490u_{0,-1} + 22u_{0,-2} - 490u_{0,1} + 22u_{0,2} - 490u_{-1,0} + 256u_{-1,-1} - 11u_{-1,-2} + 256u_{-1,1}$$
$$- 11u_{-1,2} + 22u_{-2,0} - 11u_{-2,-1} - 11u_{-2,1} - 490u_{1,0} + 256u_{1,-1} - 11u_{1,-2} + 256u_{1,1} - 11u_{1,2} + 22u_{2,0}$$
$$- 11u_{2,-1} - 11u_{2,1}) / 672$$

$$u_{xxxxx} = (-5u_{-1,0} + 4u_{-2,0} - u_{-3,0} + 5u_{1,0} - 4u_{2,0} + u_{3,0}) / 240$$
$$u_{yyyyy} = (-5u_{0,-1} + 4u_{0,-2} - u_{0,-3} + 5u_{0,1} - 4u_{0,2} + u_{0,3}) / 240$$
$$u_{xxxxy} = (-6u_{0,-1} + 6u_{0,1} + 4u_{-1,-1} - 4u_{-1,1} - u_{-2,-1} + u_{-2,1} + 4u_{1,-1} - 4u_{1,1} - u_{2,-1} + u_{2,1}) / 48$$
$$u_{xyyyy} = (-6u_{-1,0} + 4u_{-1,-1} - u_{-1,-2} + 4u_{-1,1} - u_{-1,2} + 6u_{1,0} - 4u_{1,-1} + u_{1,-2} - 4u_{1,1} + u_{1,2}) / 48$$
$$u_{xxxyy} = (-4u_{-1,0} + 2u_{-1,-1} + 2u_{-1,1} + 2u_{-2,0} - u_{-2,-1} - u_{-2,1} + 4u_{1,0} - 2u_{1,-1} - 2u_{1,1} - 2u_{2,0} + u_{2,-1} + u_{2,1}) / 24$$
$$u_{xxyyy} = (-4u_{0,-1} + 2u_{0,-2} + 4u_{0,1} - 2u_{0,2} + 2u_{-1,-1} - u_{-1,-2} - 2u_{-1,1} + u_{-1,2} + 2u_{1,-1} - u_{1,-2} - 2u_{1,1} + u_{1,2}) / 24$$

$$u_{xxxxxx} = (-20u_{0,0} + 15u_{-1,0} - 6u_{-2,0} + u_{-3,0} + 15u_{1,0} - 6u_{2,0} + u_{3,0}) / 720$$
$$u_{yyyyyy} = (-20u_{0,0} + 15u_{0,-1} - 6u_{0,-2} + u_{0,-3} + 15u_{0,1} - 6u_{0,2} + u_{0,3}) / 720$$
$$u_{xxxxxy} = (5u_{-1,-1} - 5u_{-1,1} - 4u_{-2,-1} + 4u_{-2,1} + u_{-3,-1} - u_{-3,1} - 5u_{1,-1} + 5u_{1,1} + 4u_{2,-1} - 4u_{2,1} - u_{3,-1} + u_{3,1}) / 480$$
$$u_{xyyyyy} = (5u_{-1,-1} - 4u_{-1,-2} + u_{-1,-3} - 5u_{-1,1} + 4u_{-1,2} - u_{-1,3} - 5u_{1,-1} + 4u_{1,-2} - u_{1,-3} + 5u_{1,1} - 4u_{1,2} + u_{1,3}) / 480$$
$$u_{xxxxyy} = (-12u_{0,0} + 6u_{0,-1} + 6u_{0,1} + 8u_{-1,0} - 4u_{-1,-1} - 4u_{-1,1} - 2u_{-2,0} + u_{-2,-1} + u_{-2,1} + 8u_{1,0} - 4u_{1,-1} - 4u_{1,1}$$
$$- 2u_{2,0} + u_{2,-1} + u_{2,1}) / 48$$
$$u_{xxyyyy} = (-12u_{0,0} + 8u_{0,-1} - 2u_{0,-2} + 8u_{0,1} - 2u_{0,2} + 6u_{-1,0} - 4u_{-1,-1} + u_{-1,-2} - 4u_{-1,1} + u_{-1,2} + 6u_{1,0} - 4u_{1,-1}$$
$$+ u_{1,-2} - 4u_{1,1} + u_{1,2}) / 48$$
$$u_{xxxyyy} = (4u_{-1,-1} - 2u_{-1,-2} - 4u_{-1,1} + 2u_{-1,2} - 2u_{-2,-1} + u_{-2,-2} + 2u_{-2,1} - u_{-2,2} - 4u_{1,-1} + 2u_{1,-2} + 4u_{1,1} - 2u_{1,2}$$
$$+ 2u_{2,-1} - u_{2,-2} - 2u_{2,1} + u_{2,2}) / 144$$

The smoothness indicator (denoted by $\beta^{r_7}$) for the interpolated polynomial in eqn. (A.14), written as a sum of perfect squares, is given by

$$\beta^{r_7} = \left(u_x + \frac{u_{xxx}}{10} + \frac{u_{xxxxx}}{126}\right)^2 + \left(u_y + \frac{u_{yyy}}{10} + \frac{u_{yyyyy}}{126}\right)^2$$

$$+ \frac{13}{3}\left(u_{xx} + \frac{123}{455}u_{xxxx} + \frac{85}{2002}u_{xxxxxx}\right)^2 + \frac{13}{3}\left(u_{yy} + \frac{123}{455}u_{yyyy} + \frac{85}{2002}u_{yyyyyy}\right)^2$$

$$+ \frac{7}{6}\left(u_{xy} + \frac{13}{140}u_{xxxy} + \frac{13}{140}u_{xyyy} + \frac{13}{1764}u_{xxxxxy} + \frac{13}{1764}u_{xyyyyy} + \frac{3}{350}u_{xxxyyy}\right)^2$$

$$+ \frac{781}{20}\left(u_{xxx} + \frac{26045}{49203}u_{xxxxxx}\right)^2 + \frac{781}{20}\left(u_{yyy} + \frac{26045}{49203}u_{yyyyy}\right)^2$$

$$+ \frac{47}{10}\left(u_{xxy} + \frac{533}{1974}u_{xxxxy} + \frac{781}{8460}u_{xxyyy}\right)^2 + \frac{47}{10}\left(u_{xyy} + \frac{533}{1974}u_{xyyyy} + \frac{781}{8460}u_{xxxyy}\right)^2$$

$$+ \frac{1421461}{2275}\left(u_{xxxx} + \frac{81596225}{93816426}u_{xxxxxx}\right)^2 + \frac{1421461}{2275}\left(u_{yyyy} + \frac{81596225}{93816426}u_{yyyyyy}\right)^2$$

$$+ \frac{88841}{2100}\left(u_{xxxy} + \frac{2100}{177682}\left(\frac{4740203}{105840}u_{xxxxxy} + \frac{109343}{14000}u_{xxxyyy} - \frac{1}{105840}u_{xyyyyy}\right)\right)^2$$

$$+ \frac{88841}{2100}\left(u_{xyyy} + \frac{2100}{177682}\left(\frac{4740203}{105840}u_{xyyyyy} + \frac{109343}{14000}u_{xxxyyy} - \frac{1}{105840}u_{xxxxxy}\right)\right)^2$$

$$+ \frac{1}{16800}(u_{xxxy} - u_{xyyy})^2 + \frac{5083}{270}\left(u_{xxyy} + \frac{270}{10166}\left(\frac{32021}{3150}u_{xxxxyy} + \frac{32021}{3150}u_{xxyyyy}\right)\right)^2$$

$$+ \frac{21520059541}{1377684}u_{xxxxx}^2 + \frac{21520059541}{1377684}u_{yyyyy}^2$$

$$+ \frac{2805965789}{4145400}\left(u_{xxxxy} + \frac{287}{16835794734}u_{xyyyy}\right)^2 + \frac{2805965789}{4145400}\left(u_{xyyyy} + \frac{287}{16835794734}u_{xxxxy}\right)^2$$

$$+ \frac{13478684678301203}{79549130118150}u_{xxxyy}^2 + \frac{13478684678301203}{79549130118150}u_{xxyyy}^2$$

$$+ \frac{15510384942580921}{27582029244}u_{xxxxxx}^2 + \frac{15510384942580921}{27582029244}u_{yyyyyy}^2$$

$$+ \frac{150025516743043}{8865621072}\left(u_{xxxxxy} + \frac{278487}{1200204133944344}u_{xxxyyy}\right)^2$$

$$+ \frac{150025516743043}{8865621072}\left(u_{xyyyyy} + \frac{278487}{1200204133944344}u_{xxxyyy}\right)^2$$

$$+ \frac{1461613}{315938496384}(u_{xxxxxy} + u_{xyyyyy})^2 + \frac{2278884634589}{840601125}u_{xxxxyy}^2 + \frac{2278884634589}{840601125}u_{xxyyyy}^2$$

$$+ \frac{1681}{747201000}(u_{xxxxyy} - u_{xxyyyy})^2 + \frac{18211835086345151119164258353}{11942261571939940114048000}u_{xxxyyy}^2$$

In designing a WENO scheme at seventh order, we wish to make a non-linearly hybridized interpolation using the third order stencils $S_1^{r3}, S_2^{r3}, S_3^{r3}, S_4^{r3}, S_5^{r3}$, an intermediate central fifth order stencil $S_1^{r5}$, and a large central seventh order stencil $S_1^{r7}$. We would like to have a scheme that produces a smoother transition from seventh to fifth order and then from fifth to third order. The interpolation is described by three parameters $\gamma_{Lo}, \gamma_{Avg}, \gamma_{Hi} \in (0,1)$. We set $\gamma_{Lo} = \gamma_{Avg} = \gamma_{Hi} = 0.85$. The linear weights for the stencils $S_1^{r3}, S_2^{r3}, S_3^{r3}, S_4^{r3}$ are given by

$$\gamma_1^{r3} = (1-\gamma_{Hi})(1-\gamma_{Avg})(1-\gamma_{Lo})/4 \;\; ; \;\; \gamma_2^{r3} = (1-\gamma_{Hi})(1-\gamma_{Avg})(1-\gamma_{Lo})/4 \;\; ;$$
$$\gamma_3^{r3} = (1-\gamma_{Hi})(1-\gamma_{Avg})(1-\gamma_{Lo})/4 \;\; ; \;\; \gamma_4^{r3} = (1-\gamma_{Hi})(1-\gamma_{Avg})(1-\gamma_{Lo})/4$$

The linear weight for the stencil $S_5^{r3}$ is given by $\gamma_5^{r3} = (1-\gamma_{Hi})(1-\gamma_{Avg})\gamma_{Lo}$. The linear weight for the stencil $S_1^{r5}$ is given by $\gamma_1^{r5} = (1-\gamma_{Hi})\gamma_{Avg}$. The linear weight for the stencil $S_1^{r7}$ is given by $\gamma_1^{r7} = \gamma_{Hi}$. Notice that for the linear weights we have, $\gamma_1^{r3} + \gamma_2^{r3} + \gamma_3^{r3} + \gamma_4^{r3} + \gamma_5^{r3} + \gamma_1^{r5} + \gamma_1^{r7} = 1$. We denote the seventh order central smoothness indicator by $\beta_1^{r7}$.

To avoid loss of accuracy at critical points (Borges *et al.*, 2008) we use the smoothness indicators to define

$$\tau = \frac{1}{6}\left(\left|\beta_1^{r7}-\beta_1^{r3}\right|+\left|\beta_1^{r7}-\beta_2^{r3}\right|+\left|\beta_1^{r7}-\beta_3^{r3}\right|+\left|\beta_1^{r7}-\beta_4^{r3}\right|+\left|\beta_1^{r7}-\beta_5^{r3}\right|+\left|\beta_1^{r7}-\beta_1^{r5}\right|\right).$$

The unnormalized non-linear weights are given by

$$w_1^{r3} = \gamma_1^{r3}\left(1+\tau^3/\left(\beta_1^{r3}+\varepsilon\right)^2\right) \;\; ; \;\; w_2^{r3} = \gamma_2^{r3}\left(1+\tau^3/\left(\beta_2^{r3}+\varepsilon\right)^2\right) \;\; ;$$
$$w_3^{r3} = \gamma_3^{r3}\left(1+\tau^3/\left(\beta_3^{r3}+\varepsilon\right)^2\right) \;\; ; \;\; w_4^{r3} = \gamma_4^{r3}\left(1+\tau^3/\left(\beta_4^{r3}+\varepsilon\right)^2\right) \;\; ;$$
$$w_5^{r3} = \gamma_5^{r3}\left(1+\tau^3/\left(\beta_5^{r3}+\varepsilon\right)^2\right) \;\; : \;\; w_1^{r5} = \gamma_1^{r5}\left(1+\tau^3/\left(\beta_1^{r5}+\varepsilon\right)^2\right) \;\; ;$$
$$w_1^{r7} = \gamma_1^{r7}\left(1+\tau^3/\left(\beta_1^{r7}+\varepsilon\right)^2\right)$$

where $\epsilon = 10^{-12}$ is a small number to avoid the divisions by zero. The normalization of the non-linear weights is given by

$$\overline{w}_1^{r_3} = w_1^{r_3} \Big/ \left(w_1^{r_3} + w_2^{r_3} + w_3^{r_3} + w_4^{r_3} + w_5^{r_3} + w_1^{r_5} + w_1^{r_7}\right) \ ; \quad \overline{w}_2^{r_3} = w_2^{r_3} \Big/ \left(w_1^{r_3} + w_2^{r_3} + w_3^{r_3} + w_4^{r_3} + w_5^{r_3} + w_1^{r_5} + w_1^{r_7}\right) \ ;$$

$$\overline{w}_3^{r_3} = w_3^{r_3} \Big/ \left(w_1^{r_3} + w_2^{r_3} + w_3^{r_3} + w_4^{r_3} + w_5^{r_3} + w_1^{r_5} + w_1^{r_7}\right) \ ; \quad \overline{w}_4^{r_3} = w_4^{r_3} \Big/ \left(w_1^{r_3} + w_2^{r_3} + w_3^{r_3} + w_4^{r_3} + w_5^{r_3} + w_1^{r_5} + w_1^{r_7}\right) \ ;$$

$$\overline{w}_5^{r_3} = w_5^{r_3} \Big/ \left(w_1^{r_3} + w_2^{r_3} + w_3^{r_3} + w_4^{r_3} + w_5^{r_3} + w_1^{r_5} + w_1^{r_7}\right) \ ; \quad \overline{w}_1^{r_5} = w_1^{r_5} \Big/ \left(w_1^{r_3} + w_2^{r_3} + w_3^{r_3} + w_4^{r_3} + w_5^{r_3} + w_1^{r_5} + w_1^{r_7}\right) \ ;$$

$$\overline{w}_1^{r_7} = w_1^{r_7} \Big/ \left(w_1^{r_3} + w_2^{r_3} + w_3^{r_3} + w_4^{r_3} + w_5^{r_3} + w_1^{r_5} + w_1^{r_7}\right).$$

Let $P_i^{r_3}(x,y)$ be the interpolated polynomial corresponding to stencil $S_i^{r_3}$, $P_1^{r_5}(x,y)$ be the interpolated polynomial corresponding to stencil $S_1^{r_5}$ and $P_1^{r_7}(x,y)$ be the interpolated polynomial corresponding to stencil $S_1^{r_7}$. We denote the interpolated polynomial for WENO-AO(7,5,3) as $P^{AO(7,5,3)}(x,y)$. Our task in this paragraph is to describe the construction of the order-preserving, non-linearly hybridized, seventh order polynomial $P^{AO(7,5,3)}(x,y)$. The non-linear weights should be combined in such a way that when all the smoothness indicators seem to have almost similar values then only the higher order scheme is obtained. The non-linearly hybridized seventh order accurate interpolation is given by

$$P^{AO(7,5,3)}(x) = \frac{\overline{w}_1^{r_7}}{\gamma_1^{r_7}} \begin{pmatrix} P_1^{r_7}(x,y) - \gamma_1^{r_5} P_1^{r_5}(x,y) - \gamma_1^{r_3} P_1^{r_3}(x,y) - \gamma_2^{r_3} P_2^{r_3}(x,y) \\ - \gamma_3^{r_3} P_3^{r_3}(x,y) - \gamma_4^{r_3} P_4^{r_3}(x,y) - \gamma_5^{r_3} P_5^{r_3}(x,y) \end{pmatrix}$$
$$+ \overline{w}_1^{r_5} P_1^{r_5}(x,y) + \overline{w}_1^{r_3} P_1^{r_3}(x,y) + \overline{w}_2^{r_3} P_2^{r_3}(x,y) + \overline{w}_3^{r_3} P_3^{r_3}(x,y)$$
$$+ \overline{w}_4^{r_3} P_4^{r_3}(x,y) + \overline{w}_5^{r_3} P_5^{r_3}(x,y)$$

This completes our description of the seventh order 2D WENO-AO(7,5,3) interpolation.

### A.4) WENO-AO(9,5,3) Interpolation in 2D

We will not describe all elements of a ninth order scheme because the expressions become much too large. However, we provide the most important ingredients that most readers will find hard to discover on their own. This includes the polynomial in which the interpolation is done and the smoothness indicator expressed as a sum of perfect squares. The 2D WENO-AO(9,5,3) interpolation consists of a non-linear hybridization between a large, centered, ninth order stencil denoted by $S_1^{r_9}$, an intermediate central fifth order stencil denoted by $S_1^{r_5}$ and the five smaller stencils described in Appendix A.1). The larger central ninth order accurate stencil is shown in Fig. 6c and The intermediate fifth order central stencil is shown in Fig. 6a. The ninth order accurate 2D polynomial is given

$$\begin{aligned}
u(x,y) = &\ u_{00} + u_x L_1(x) + u_y L_1(y) + u_{xx} L_2(x) + u_{yy} L_2(y) + u_{xy} L_1(x) L_1(y) \\
&+ u_{xxx} L_3(x) + u_{yyy} L_3(y) + u_{xxy} L_2(x) L_1(y) + u_{xyy} L_1(x) L_2(y) \\
&+ u_{xxxx} L_4(x) + u_{yyyy} L_4(y) + u_{xxxy} L_3(x) L_1(y) + u_{xyyy} L_1(x) L_3(y) + u_{xxyy} L_2(x) L_2(y) \\
&+ u_{xxxxx} L_5(x) + u_{yyyyy} L_5(y) + u_{xxxxy} L_4(x) L_1(y) + u_{xyyyy} L_1(x) L_4(y) \\
&+ u_{xxxyy} L_3(x) L_2(y) + u_{xxyyy} L_2(x) L_3(y) \\
&+ u_{xxxxxx} L_6(x) + u_{yyyyyy} L_6(y) + u_{xxxxxy} L_5(x) L_1(y) + u_{xyyyyy} L_1(x) L_5(y) \\
&+ u_{xxxxyy} L_4(x) L_2(y) + u_{xxyyyy} L_2(x) L_4(y) + u_{xxxyyy} L_3(x) L_3(y) \\
&+ u_{xxxxxxx} L_7(x) + u_{yyyyyyy} L_7(y) + u_{xxxxxxy} L_6(x) L_1(y) + u_{xyyyyyy} L_1(x) L_6(y) \\
&+ u_{xxxxxyy} L_5(x) L_2(y) + u_{xxyyyyy} L_2(x) L_5(y) + u_{xxxxyyy} L_4(x) L_3(y) + u_{xxxyyyy} L_3(x) L_4(y) \\
&+ u_{xxxxxxxx} L_8(x) + u_{yyyyyyyy} L_8(y) + u_{xxxxxxxy} L_7(x) L_1(y) + u_{xyyyyyyy} L_1(x) L_7(y) \\
&+ u_{xxxxxxyy} L_6(x) L_2(y) + u_{xxyyyyyy} L_2(x) L_6(y) + u_{xxxxxyyy} L_5(x) L_3(y) + u_{xxxyyyyy} L_3(x) L_5(y) \\
&+ u_{xxxxyyyy} L_4(x) L_4(y)
\end{aligned}$$

(A.15)

For the above ninth order polynomial, the smoothness indicator expressed as a sum of perfect squares is given by:-

$$\beta^{r_9} = \left(u_x + \frac{u_{xxx}}{10} + \frac{u_{xxxxx}}{126} + \frac{u_{xxxxxxx}}{1716}\right)^2 + \left(u_y + \frac{u_{yyy}}{10} + \frac{u_{yyyyy}}{126} + \frac{u_{yyyyyyy}}{1716}\right)^2$$

$$+ \frac{13}{3}\left(u_{xx} + \frac{123}{455}u_{xxxx} + \frac{85}{2002}u_{xxxxxx} + \frac{29}{5577}u_{xxxxxxxx}\right)^2 + \frac{13}{3}\left(u_{yy} + \frac{123}{455}u_{yyyy} + \frac{85}{2002}u_{yyyyyy} + \frac{29}{5577}u_{yyyyyyyy}\right)^2$$

$$+ \frac{7}{6}\left(\begin{array}{l} u_{xy} + \frac{13}{140}u_{xxxy} + \frac{13}{140}u_{xyyy} + \frac{13}{1764}u_{xxxxxy} + \frac{13}{1764}u_{xyyyyy} + \frac{3}{350}u_{xxxyyy} \\ + \frac{3}{7}\left(\frac{u_{xxxxxxy}}{792} + \frac{u_{xyyyyyyy}}{792} + \frac{u_{xxxxxyyy}}{630} + \frac{u_{xxxyyyyy}}{630}\right) \end{array}\right)^2$$

$$+ \frac{781}{20}\left(u_{xxx} + \frac{26045}{49203}u_{xxxxx} + \frac{8395}{60918}u_{xxxxxxx}\right)^2 + \frac{781}{20}\left(u_{yyy} + \frac{26045}{49203}u_{yyyyy} + \frac{8395}{60918}u_{yyyyyyy}\right)^2$$

$$+ \frac{47}{10}\left(u_{xxy} + \frac{533}{1974}u_{xxxxy} + \frac{781}{8460}u_{xxyyy} + \frac{10}{94}\left(\frac{1105}{2772}u_{xxxxxxy} + \frac{41}{175}u_{xxxxyyy} + \frac{781}{11340}u_{xxyyyyy}\right)\right)^2$$

$$+ \frac{47}{10}\left(u_{xyy} + \frac{533}{1974}u_{xyyyy} + \frac{781}{8460}u_{xxxyy} + \frac{10}{94}\left(\frac{1105}{2772}u_{xyyyyyy} + \frac{41}{175}u_{xxxyyyy} + \frac{781}{11340}u_{xxxxxyy}\right)\right)^2$$

$$+ \frac{1421461}{2275}\left(u_{xxxx} + \frac{2275}{2842922}\left(\frac{3263849}{3003}u_{xxxxxx} + \frac{35339362}{83655}u_{xxxxxxxx}\right)\right)^2$$

$$+ \frac{1421461}{2275}\left(u_{yyyy} + \frac{2275}{2842922}\left(\frac{3263849}{3003}u_{yyyyyy} + \frac{35339362}{83655}u_{yyyyyyyy}\right)\right)^2$$

$$+ \frac{88841}{2100}\left(u_{xxxy} + \frac{2100}{177682}\left(\begin{array}{l}\frac{4740203}{105840}u_{xxxxxy} + \frac{109343}{14000}u_{xxxyyy} - \frac{1}{105840}u_{xyyyyy} \\ + \frac{1292831}{110880}u_{xxxxxxxy} + \frac{364631}{88200}u_{xxxxxyyy} + \frac{109343}{176400}u_{xxxyyyyy} - \frac{1}{1441440}u_{xyyyyyyy}\end{array}\right)\right)^2$$

$$+ \frac{88841}{2100}\left(u_{xyyy} + \frac{2100}{177682}\left(\begin{array}{l}\frac{4740203}{105840}u_{xyyyyy} + \frac{109343}{14000}u_{xxxyyy} - \frac{1}{105840}u_{xxxxxy} \\ + \frac{1292831}{110880}u_{xyyyyyyy} + \frac{364631}{88200}u_{xxxyyyyy} + \frac{109343}{176400}u_{xxxxxyyy} - \frac{1}{1441440}u_{xxxxxxxy}\end{array}\right)\right)^2$$

$$+ \frac{1}{16800}(u_{xxxy} - u_{xyyy})^2 + \frac{5083}{270}\left(u_{xxyy} + \frac{270}{10166}\left(\begin{array}{l}\frac{32021}{3150}(u_{xxxxyy} + u_{xxyyyy}) + \frac{1207}{756}(u_{xxxxxxyy} + u_{xxyyyyyy}) \\ + \frac{3362}{1225}u_{xxxxyyyy}\end{array}\right)\right)^2$$

$$+\frac{21520059541}{1377684}\left(u_{xxxxx}+\frac{1377684}{43040119082}\left(\frac{34399031911}{852852}u_{xxxxxxx}\right)\right)^2$$

$$+\frac{21520059541}{1377684}\left(u_{yyyyy}+\frac{1377684}{43040119082}\left(\frac{34399031911}{852852}u_{yyyyyyy}\right)\right)^2$$

$$+\frac{2805965789}{4145400}\left(u_{xxxxy}+\frac{4145400}{5611931578}\left(\begin{array}{c}\frac{6442838623}{5471928}u_{xxxxxy}+\frac{1295061137}{10363500}u_{xxxxyyy}\\-\frac{41}{1776600}u_{xxyyy}-\frac{41}{22385160}u_{xyyyyy}\end{array}\right)\right)^2$$

$$+\frac{2805965789}{4145400}\left(u_{xyyyy}+\frac{4145400}{5611931578}\left(\begin{array}{c}\frac{6442838623}{5471928}u_{xyyyyyy}+\frac{1295061137}{10363500}u_{xxxyyyy}\\-\frac{41}{1776600}u_{xxxyy}-\frac{41}{22385160}u_{xxxxxy}\end{array}\right)\right)^2$$

$$+\frac{13478684678301203}{79549130118150}\left(u_{xxxyy}+\frac{39774565059075}{13478684678301203}\left(\begin{array}{c}\frac{102738105636515023}{572753736850680}u_{xxxxyy}\\+\frac{33964215589171433}{371229273884700}u_{xxxyyyy}\\+\frac{2302712395}{140006469007944}u_{xyyyyyy}\end{array}\right)\right)^2$$

$$+\frac{13478684678301203}{79549130118150}\left(u_{xxyyy}+\frac{39774565059075}{13478684678301203}\left(\begin{array}{c}\frac{102738105636515023}{572753736850680}u_{xyyyyy}\\+\frac{33964215589171433}{371229273884700}u_{xxxxyy}\\+\frac{2302712395}{140006469007944}u_{xxxxxy}\end{array}\right)\right)^2$$

$$+\frac{15510384942580921}{27582029244}\left(u_{xxxxxx}+\frac{5423630339859998294}{3024525063803279595}u_{xxxxxxxx}\right)^2$$

$$+\frac{15510384942580921}{27582029244}\left(u_{yyyyyy}+\frac{5423630339859998294}{3024525063803279595}u_{yyyyyyyy}\right)^2$$

$$+\frac{150025516743043}{8865621072}\left(u_{xxxxxy}+\frac{4432810536}{150025516743043}\left(\begin{array}{c}\frac{8545967960535563}{195580973952}u_{xxxxxxxy}\\+\frac{1206359027303963}{386147051136}u_{xxxxxyyy}\\+\frac{128619}{79681137536}u_{xyyyyyyy}\\-\frac{2813}{358206912}u_{xxxyyy}\\-\frac{364631}{1579692481920}u_{xxxyyyyy}\end{array}\right)\right)$$

$$+\frac{150025516743043}{8865621072}\left(u_{xyyyyy}+\frac{4432810536}{150025516743043}\left(\begin{array}{l}\frac{8545967960535563}{195580973952}u_{xyyyyyy}\\+\frac{1206359027303963}{386147051136}u_{xxxyyyyy}\\+\frac{128619}{79681137536}u_{xxxxxxy}\\-\frac{2813}{358206912}u_{xxxyyy}\\-\frac{364631}{1579692481920}u_{xxxxxyyy}\end{array}\right)\right)^2$$

$$+\frac{1461613}{315938496384}(u_{xxxxxy}+u_{xyyyyy})^2$$

$$+\frac{2278884634589}{840601125}\left(u_{xxxxyy}+\frac{840601125}{4557769269178}\left(\begin{array}{l}\frac{74618124571}{15823080}u_{xxxxxxy}+\frac{5742439135469}{3922805250}u_{xxxyyyy}\\-\frac{41}{58017960}u_{xyyyyyyy}\end{array}\right)\right)^2$$

$$+\frac{2278884634589}{840601125}\left(u_{xxyyyy}+\frac{840601125}{4557769269178}\left(\begin{array}{l}\frac{74618124571}{15823080}u_{xyyyyyy}+\frac{5742439135469}{3922805250}u_{xxxxyyy}\\-\frac{41}{58017960}u_{xxxxxxyy}\end{array}\right)\right)^2$$

$$+\frac{1681}{747201000}(u_{xxxxyy}-u_{xxyyyy})^2$$

$$+\frac{18211835086345151119164258353}{119422615719399401 1404800}\left(u_{xxxyyy}+\frac{5971130785969970057 02400}{18211835086345151119164258353}\left(\begin{array}{l}\frac{34477149017196090259}{4346983212186138201513472}(u_{xxxxxxy}+u_{xyyyyyyy})\\+\frac{85024449502692624249 8309371309}{52665373532255135902 9516800}(u_{xxxxxyy}+u_{xxyyyyy})\end{array}\right)\right)^2$$

$$+ \frac{12210527897166191835083}{443141066068272} u_{xxxxxx}^2 + \frac{12210527897166191835083}{443141066068272} u_{yyyyyy}^2$$

$$+ \frac{7290935622087218652389387744441}{119680885894973394082291} \left( u_{xxxxxy} + \frac{59840442947486697041456}{7290935622087218652389387744441} \left( \begin{array}{c} -\frac{1012830861343669}{213502365304291055520} u_{xxxyyy} \\ -\frac{92906143860066025}{125539390798923140645760} u_{xyyyyy} \end{array} \right) \right)^2$$

$$+ \frac{7290935622087218652389387744441}{119680885894973394082291} \left( u_{xyyyyyy} + \frac{59840442947486697041456}{7290935622087218652389387744441} \left( \begin{array}{c} -\frac{1012830861343669}{213502365304291055520} u_{xxxyyy} \\ -\frac{92906143860066025}{125539390798923140645760} u_{xxxxxy} \end{array} \right) \right)^2$$

$$+ \frac{10786969275341029621981757001451763059666}{15915747916235293957233413976659809} \times$$

$$\left( u_{xxxxyy} - \frac{30351858961788546932481157806}{26967423188352574054954392503629407649165} u_{xxxyyy} \left( \begin{array}{c} -\frac{1012830861343669}{213502365304291055520} u_{xxxyyy} \\ -\frac{92906143860066025}{125539390798923140645760} u_{xxxxyy} \end{array} \right) \right)^2$$

$$+ \frac{10786969275341029621981757001451763059666}{15915747916235293957233413976659809} \times$$

$$\left( u_{xxyyyy} - \frac{30351858961788546932481157806}{26967423188352574054954392503629407649165} u_{xxxxyy} \left( \begin{array}{c} -\frac{1012830861343669}{213502365304291055520} u_{xxxyyy} \\ -\frac{92906143860066025}{125539390798923140645760} u_{xxyyyy} \end{array} \right) \right)^2$$

$$+ \frac{2763606153533153038150205824278415130488282001719}{113263177391080811030808448515243512126493000} \left( u_{xxxxyyy}^2 + u_{xxxyyyy}^2 \right)$$

$$+ \frac{75509368098103789336083731407561}{428182013282630292264515} \left( u_{xxxxxxx}^2 + u_{yyyyyyy}^2 \right)$$

$$+ \frac{1997376815537062209802912385779467592212372365}{100349412338258061005083687187257696748484264781248} \left( u_{xxxxxxy} - u_{xyyyyyy} \right)^2$$

$$+\frac{12479128819356707044610023473398906994677171}{41805143600226707413387447323546996 0}\times$$

$$\left(u_{xxxxxxy}+\frac{20902571800113353706693723661773498 0}{12479128819356707044610023473398906994677171}\times\right.$$

$$\left.\left(\begin{array}{c}-\frac{7977829830562415501482175809036941341308898 9}{14736626986737197769977324691835046375651535387456 0}u_{xxxxxyyy}\\-\frac{2060961342250291550354757696295988240818 07593}{4912208995579065923325774897278348791883845129152 0}u_{xxxyyyyy}\end{array}\right)\right)^2$$

$$+\frac{12479128819356707044610023473398906994677171}{41805143600226707413387447323546996 0}\times$$

$$\left(u_{xyyyyyy}+\frac{20902571800113353706693723661773498 0}{12479128819356707044610023473398906994677171}\times\right.$$

$$\left.\left(\begin{array}{c}-\frac{7977829830562415501482175809036941341308898 9}{14736626986737197769977324691835046375651535387456 0}u_{xxxyyyyy}\\-\frac{2060961342250291550354757696295988240818 07593}{4912208995579065923325774897278348791883845129152 0}u_{xxxxxyyy}\end{array}\right)\right)^2$$

$$+\frac{26109152643905}{46695768186740593536}\left(u_{xxxxxyy}+u_{xxyyyyy}\right)^2$$

$$+\frac{71208597550871338859711633}{29184855116712870960}\left(u_{xxxxxyy}-\frac{26993580654975}{142417195101742677719423266}u_{xxxyyyy}\right)^2$$

$$+\frac{71208597550871338859711633}{29184855116712870960}\left(u_{xxyyyyyy}-\frac{26993580654975}{142417195101742677719423266}u_{xxxxyyy}\right)^2$$

$$+\frac{74444650667725915845814934912264574382747147172348831 83513}{12204108786784196192078089932082752522765677634633600}\left(u_{xxxxxyy}^2+u_{xxyyyyy}^2\right)^2$$

$$+\frac{25969660968505236829899732631451382865756595852879535 74573816311975259530963}{87899481081979187574464425121693808582659625310765620 7181523209141421581473044480 0}\times$$

$$\left(u_{xxxxxyyy}-u_{xxxyyyyy}\right)^2$$

$$+\frac{1117534769052383552198004388237023006772308002354 09}{28625517941962741122230963620029922384846800 0}u_{xxxxyyyy}^2$$

## Appendix B : A WENO-Stabilized Transcription Strategy from 3D Point Values to Volume Averages

In 3D, we non-linearly hybridize between nine smaller third order stencils and a larger fifth order stencil. The third order stencils give rise to interpolating polynomials that have the form

$$u(x,y,z) = u_0 + u_x x + u_y y + u_z z + u_{xx}(x^2 - 1/12) + u_{yy}(y^2 - 1/12) + u_{zz}(z^2 - 1/12)$$
$$+ u_{xy} xy + u_{yz} yz + u_{xz} xz \ .$$
(B.1)

and a smoothness indicator of the form

$$IS = u_x^2 + 13 u_{xx}^2 / 3 + u_y^2 + 13 u_{yy}^2 / 3 + u_z^2 + 13 u_{zz}^2 / 3 + 7 u_{xy}^2 / 6 + 7 u_{yz}^2 / 6 + 7 u_{xz}^2 / 6 \ .$$
(B.2)

The large fifth order stencil gives rise to an interpolating polynomial of the form

$$u(x,y,z) = u_0 + u_x x + u_y y + u_z z + u_{xx}(x^2 - 1/12) + u_{yy}(y^2 - 1/12) + u_{zz}(z^2 - 1/12)$$
$$+ u_{xy} xy + u_{yz} yz + u_{xz} xz + u_{xxx}(x^3 - 3x/20) + u_{yyy}(y^3 - 3y/20) + u_{zzz}(z^3 - 3z/20)$$
$$+ u_{xxy}(x^2 - 1/12)y + u_{xxz}(x^2 - 1/12)z + u_{xyy} x(y^2 - 1/12) + u_{yyz}(y^2 - 1/12)z$$
$$+ u_{xzz} x(z^2 - 1/12) + u_{yzz} y(z^2 - 1/12) + u_{xyz} xyz$$
$$+ u_{xxxx}\left(x^4 - 3x^2/14 + 3/560\right) + u_{yyyy}\left(y^4 - 3y^2/14 + 3/560\right) + u_{zzzz}\left(z^4 - 3z^2/14 + 3/560\right)$$
$$+ u_{xxxy}(x^3 - 3x/20)y + u_{xxxz}(x^3 - 3x/20)z + u_{xyyy} x(y^3 - 3y/20) + u_{yyyz}(y^3 - 3y/20)z$$
$$+ u_{xzzz} x(z^3 - 3z/20) + u_{yzzz} y(z^3 - 3z/20) + u_{xxyy}(x^2 - 1/12)(y^2 - 1/12)$$
$$+ u_{yyzz}(y^2 - 1/12)(z^2 - 1/12) + u_{xxzz}(x^2 - 1/12)(z^2 - 1/12) + u_{xxyz}(x^2 - 1/12)yz$$
$$+ u_{xyyz} x(y^2 - 1/12)z + u_{xyzz} xy(z^2 - 1/12) \ .$$
(B.3)

and a smoothness indicator of the form

$$IS = (u_x + u_{xxx}/10)^2 + 13(u_{xx} + 123 u_{xxxx}/455)^2/3 + 781 u_{xxx}^2/20 + 1421461 u_{xxxx}^2/2275 + (u_y + u_{yyy}/10)^2$$
$$+ 13(u_{yy} + 123 u_{yyyy}/455)^2/3 + 781 u_{yyy}^2/20 + 1421461 u_{yyyy}^2/2275 + (u_z + u_{zzz}/10)^2$$
$$+ 13(u_{zz} + 123 u_{zzzz}/455)^2/3 + 781 u_{zzz}^2/20 + 1421461 u_{zzzz}^2/2275 + 47 u_{xxy}^2/10 + 47 u_{xyy}^2/10$$
$$+ 47 u_{xxz}^2/10 + 47 u_{xzz}^2/10 + 47 u_{yyz}^2/10 + 47 u_{yzz}^2/10 + 61 u_{xyz}^2/48$$
$$+ 7(u_{xy} + 13 u_{xxxy}/140 + 13 u_{xyyy}/140)^2/6 + 7(u_{yz} + 13 u_{yyyz}/140 + 13 u_{yzzz}/140)^2/6$$
$$+ 7(u_{xz} + 13 u_{xxxz}/140 + 13 u_{xzzz}/140)^2/6 + 88841(u_{xxxy}^2 + u_{xyyy}^2)/2100 + (u_{xxxy} - u_{xyyy})^2/16800$$
$$+ 5083 u_{xxyy}^2/270 + 88841(u_{yyyz}^2 + u_{yzzz}^2)/2100 + (u_{yyyz} - u_{yzzz})^2/16800 + 5083 u_{yyzz}^2/270$$
$$+ 88841(u_{xzzz}^2 + u_{xxxz}^2)/2100 + (u_{xzzz} - u_{xxxz})^2/16800 + 5083 u_{xxzz}^2/270 + 10999 u_{xxyz}^2/2160$$
$$+ 10999 u_{xyyz}^2/2160 + 10999 u_{xyzz}^2/2160 \ .$$
(B.4)

Please note that this is not a finite volume style reconstruction but rather a pointwise interpolation. The eight smaller third order interpolating stencils are one-sided and they can be identified on the mesh. We shall not describe all of them, but we will provide all the coefficients for one of them below:-

Stencil-1 (Biased in 1$^{st}$ octant, *x, y, z* >0):

$$u_0 = (27u_{0,0,0} - 2u_{0,0,+1} + u_{0,0,+2} - 2u_{0,+1,0} + u_{0,+2,0} - 2u_{+1,0,0} + u_{+2,0,0})/24$$
$$u_x = (-3u_{0,0,0} + 4u_{+1,0,0} - u_{+2,0,0})/2$$
$$u_y = (-3u_{0,0,0} + 4u_{0,+1,0} - u_{0,+2,0})/2$$
$$u_z = (-3u_{0,0,0} + 4u_{0,0,+1} - u_{0,0,+2})/2$$
$$u_{xx} = (u_{0,0,0} - 2u_{+1,0,0} + u_{+2,0,0})/2$$
$$u_{yy} = (u_{0,0,0} - 2u_{0,+1,0} + u_{0,+2,0})/2$$
$$u_{zz} = (u_{0,0,0} - 2u_{0,0,+1} + u_{0,0,+2})/2$$
$$u_{xy} = u_{0,0,0} - u_{0,+1,0} - u_{+1,0,0} + u_{+1,+1,0}$$
$$u_{yz} = u_{0,0,0} - u_{0,0,+1} - u_{0,+1,0} + u_{0,+1,+1}$$
$$u_{xz} = u_{0,0,0} - u_{0,0,+1} - u_{+1,0,0} + u_{+1,0,+1}$$

These eight one-sided stencils have to be non-linearly hybridized with a smaller third order central stencil for which we provide all the coefficients below:-

$$u_0 = (18u_{0,0,0} + u_{0,0,-1} + u_{0,0,+1} + u_{0,-1,0} + u_{0,+1,0} + u_{-1,0,0} + u_{+1,0,0})/24$$
$$u_x = (-u_{-1,0,0} + u_{+1,0,0})/2 \;;\; u_y = (-u_{0,-1,0} + u_{0,+1,0})/2 \;;\; u_z = (-u_{0,0,-1} + u_{0,0,+1})/2$$
$$u_{xx} = (-2u_{0,0,0} + u_{-1,0,0} + u_{+1,0,0})/2$$
$$u_{yy} = (-2u_{0,0,0} + u_{0,-1,0} + u_{0,+1,0})/2$$
$$u_{zz} = (-2u_{0,0,0} + u_{0,0,-1} + u_{0,0,+1})/2$$
$$u_{xy} = (u_{-1,-1,0} - u_{-1,+1,0} - u_{+1,-1,0} + u_{+1,+1,0})/4$$
$$u_{yz} = (u_{0,-1,-1} - u_{0,-1,+1} - u_{0,+1,-1} + u_{0,+1,+1})/4$$
$$u_{xz} = (u_{-1,0,-1} - u_{-1,0,+1} - u_{+1,0,-1} + u_{+1,0,+1})/4$$

In addition, we will need a larger fifth order interpolating stencil against which we have to perform non-linear hybridization via a vis the nine smaller third order stencil. For the large fifth order stencil we provide all the coefficients below:-

$$u_0 = (4134u_{0,0,0} + 268u_{0,0,-1} - 17u_{0,0,-2} + 268u_{0,0,+1} - 17u_{0,0,+2} + 268u_{0,-1,0} + 10u_{0,-1,-1} + 10u_{0,-1,+1}$$
$$- 17u_{0,-2,0} + 268u_{0,+1,0} + 10u_{0,+1,-1} + 10u_{0,+1,+1} - 17u_{0,+2,0} + 268u_{-1,0,0} + 10u_{-1,0,-1} + 10u_{-1,0,+1}$$
$$+ 10u_{-1,-1,0} + 10u_{-1,+1,0} - 17u_{-2,0,0} + 268u_{+1,0,0} + 10u_{+1,0,-1} + 10u_{+1,0,+1} + 10u_{+1,-1,0} + 10u_{+1,+1,0}$$
$$- 17u_{+2,0,0}) / 5760$$

$$u_x = (-134u_{-1,0,0} - 5u_{-1,0,-1} - 5u_{-1,0,+1} - 5u_{-1,-1,0} - 5u_{-1,+1,0} + 17u_{-2,0,0} + 134u_{+1,0,0} + 5u_{+1,0,-1}$$
$$+ 5u_{+1,0,+1} + 5u_{+1,-1,0} + 5u_{+1,+1,0} - 17u_{+2,0,0}) / 240$$

$$u_y = (-134u_{0,-1,0} - 5u_{0,-1,-1} - 5u_{0,-1,+1} + 17u_{0,-2,0} + 134u_{0,+1,0} + 5u_{0,+1,-1} + 5u_{0,+1,+1} - 17u_{0,+2,0}$$
$$- 5u_{-1,-1,0} + 5u_{-1,+1,0} - 5u_{+1,-1,0} + 5u_{+1,+1,0}) / 240$$

$$u_z = (-134u_{0,0,-1} + 17u_{0,0,-2} + 134u_{0,0,+1} - 17u_{0,0,+2} - 5u_{0,-1,-1} + 5u_{0,-1,+1} - 5u_{0,+1,-1} + 5u_{0,+1,+1}$$
$$- 5u_{-1,0,-1} + 5u_{-1,0,+1} - 5u_{+1,0,-1} + 5u_{+1,0,+1}) / 240$$

$$u_{xx} = (-346u_{0,0,0} - 14u_{0,0,-1} - 14u_{0,0,+1} - 14u_{0,-1,0} - 14u_{0,+1,0} + 184u_{-1,0,0} + 7u_{-1,0,-1} + 7u_{-1,0,+1}$$
$$+ 7u_{-1,-1,0} + 7u_{-1,+1,0} - 11u_{-2,0,0} + 184u_{+1,0,0} + 7u_{+1,0,-1} + 7u_{+1,0,+1} + 7u_{+1,-1,0} + 7u_{+1,+1,0}$$
$$- 11u_{+2,0,0}) / 336$$

$$u_{yy} = (-346u_{0,0,0} - 14u_{0,0,-1} - 14u_{0,0,+1} + 184u_{0,-1,0} + 7u_{0,-1,-1} + 7u_{0,-1,+1} - 11u_{0,-2,0} + 184u_{0,+1,0}$$
$$+ 7u_{0,+1,-1} + 7u_{0,+1,+1} - 11u_{0,+2,0} - 14u_{-1,0,0} + 7u_{-1,-1,0} + 7u_{-1,+1,0} - 14u_{+1,0,0} + 7u_{+1,-1,0}$$
$$+ 7u_{+1,+1,0}) / 336$$

$$u_{zz} = (-346u_{0,0,0} + 184u_{0,0,-1} - 11u_{0,0,-2} + 184u_{0,0,+1} - 11u_{0,0,+2} - 14u_{0,-1,0} + 7u_{0,-1,-1} + 7u_{0,-1,+1}$$
$$- 14u_{0,+1,0} + 7u_{0,+1,-1} + 7u_{0,+1,+1} - 14u_{-1,0,0} + 7u_{-1,0,-1} + 7u_{-1,0,+1} - 14u_{+1,0,0} + 7u_{+1,0,-1}$$
$$+ 7u_{+1,0,+1}) / 336$$

$$u_{xy} = (386u_{-1,-1,0} - 10u_{-1,-1,-1} + 5u_{-1,-1,-2} - 10u_{-1,-1,+1} + 5u_{-1,-1,+2} - 34u_{-1,-2,0} - 386u_{-1,+1,0} + 10u_{-1,+1,-1}$$
$$- 5u_{-1,+1,-2} + 10u_{-1,+1,+1} - 5u_{-1,+1,+2} + 34u_{-1,+2,0} - 34u_{-2,-1,0} + 34u_{-2,+1,0} - 386u_{+1,-1,0} + 10u_{+1,-1,-1}$$
$$- 5u_{+1,-1,-2} + 10u_{+1,-1,+1} - 5u_{+1,-1,+2} + 34u_{+1,-2,0} + 386u_{+1,+1,0} - 10u_{+1,+1,-1} + 5u_{+1,+1,-2} - 10u_{+1,+1,+1}$$
$$+ 5u_{+1,+1,+2} - 34u_{+1,+2,0} + 34u_{+2,-1,0} - 34u_{+2,+1,0}) / 960$$

$$u_{yz} = (386u_{0,-1,-1} - 34u_{0,-1,-2} - 386u_{0,-1,+1} + 34u_{0,-1,+2} - 34u_{0,-2,-1} + 34u_{0,-2,+1} - 386u_{0,+1,-1} + 34u_{0,+1,-2}$$
$$+ 386u_{0,+1,+1} - 34u_{0,+1,+2} + 34u_{0,+2,-1} - 34u_{0,+2,+1} - 10u_{-1,-1,-1} + 10u_{-1,-1,+1} + 10u_{-1,+1,-1} - 10u_{-1,+1,+1}$$
$$+ 5u_{-2,-1,-1} - 5u_{-2,-1,+1} - 5u_{-2,+1,-1} + 5u_{-2,+1,+1} - 10u_{+1,-1,-1} + 10u_{+1,-1,+1} + 10u_{+1,+1,-1} - 10u_{+1,+1,+1}$$
$$+ 5u_{+2,-1,-1} - 5u_{+2,-1,+1} - 5u_{+2,+1,-1} + 5u_{+2,+1,+1}) / 960$$

$$u_{xz} = (386u_{-1,0,-1} - 34u_{-1,0,-2} - 386u_{-1,0,+1} + 34u_{-1,0,+2} - 10u_{-1,-1,-1} + 10u_{-1,-1,+1} + 5u_{-1,-2,-1} - 5u_{-1,-2,+1}$$
$$- 10u_{-1,+1,-1} + 10u_{-1,+1,+1} + 5u_{-1,+2,-1} - 5u_{-1,+2,+1} - 34u_{-2,0,-1} + 34u_{-2,0,+1} - 386u_{+1,0,-1} + 34u_{+1,0,-2}$$
$$+ 386u_{+1,0,+1} - 34u_{+1,0,+2} + 10u_{+1,-1,-1} - 10u_{+1,-1,+1} - 5u_{+1,-2,-1} + 5u_{+1,-2,+1} + 10u_{+1,+1,-1} - 10u_{+1,+1,+1}$$
$$- 5u_{+1,+2,-1} + 5u_{+1,+2,+1} + 34u_{+2,0,-1} - 34u_{+2,0,+1}) / 960$$

$$u_{xxx} = (2u_{-1,0,0} - u_{-2,0,0} - 2u_{+1,0,0} + u_{+2,0,0})/12$$
$$u_{yyy} = (2u_{0,-1,0} - u_{0,-2,0} - 2u_{0,+1,0} + u_{0,+2,0})/12$$
$$u_{zzz} = (2u_{0,0,-1} - u_{0,0,-2} - 2u_{0,0,+1} + u_{0,0,+2})/12$$
$$u_{xxy} = (2u_{0,-1,0} - 2u_{0,+1,0} - u_{-1,-1,0} + u_{-1,+1,0} - u_{+1,-1,0} + u_{+1,+1,0})/4$$
$$u_{xxz} = (2u_{0,0,-1} - 2u_{0,0,+1} - u_{-1,0,-1} + u_{-1,0,+1} - u_{+1,0,-1} + u_{+1,0,+1})/4$$
$$u_{xyy} = (2u_{-1,0,0} - u_{-1,-1,0} - u_{-1,+1,0} - 2u_{+1,0,0} + u_{+1,-1,0} + u_{+1,+1,0})/4$$
$$u_{yyz} = (2u_{0,0,-1} - 2u_{0,0,+1} - u_{0,-1,-1} + u_{0,-1,+1} - u_{0,+1,-1} + u_{0,+1,+1})/4$$
$$u_{xzz} = (2u_{-1,0,0} - u_{-1,0,-1} - u_{-1,0,+1} - 2u_{+1,0,0} + u_{+1,0,-1} + u_{+1,0,+1})/4$$
$$u_{yzz} = (2u_{0,-1,0} - u_{0,-1,-1} - u_{0,-1,+1} - 2u_{0,+1,0} + u_{0,+1,-1} + u_{0,+1,+1})/4$$
$$\begin{aligned}u_{xyz} = (&-8u_{-1,-1,-1} + u_{-1,-1,-2} + 8u_{-1,-1,+1} - u_{-1,-1,+2} + u_{-1,-2,-1} - u_{-1,-2,+1} + 8u_{-1,+1,-1} - u_{-1,+1,-2} - 8u_{-1,+1,+1} \\ &+ u_{-1,+1,+2} - u_{-1,+2,-1} + u_{-1,+2,+1} + u_{-2,-1,-1} - u_{-2,-1,+1} - u_{-2,+1,-1} + u_{-2,+1,+1} + 8u_{+1,-1,-1} - u_{+1,-1,-2} \\ &- 8u_{+1,-1,+1} + u_{+1,-1,+2} - u_{+1,-2,-1} + u_{+1,-2,+1} - 8u_{+1,+1,-1} + u_{+1,+1,-2} + 8u_{+1,+1,+1} - u_{+1,+1,+2} + u_{+1,+2,-1} \\ &- u_{+1,+2,+1} - u_{+2,-1,-1} + u_{+2,-1,+1} + u_{+2,+1,-1} - u_{+2,+1,+1})/16\end{aligned}$$

$$u_{xxxx} = (6u_{0,0,0} - 4u_{-1,0,0} + u_{-2,0,0} - 4u_{+1,0,0} + u_{+2,0,0})/24$$
$$u_{yyyy} = (6u_{0,0,0} - 4u_{0,-1,0} + u_{0,-2,0} - 4u_{0,+1,0} + u_{0,+2,0})/24$$
$$u_{zzzz} = (6u_{0,0,0} - 4u_{0,0,-1} + u_{0,0,-2} - 4u_{0,0,+1} + u_{0,0,+2})/24$$

$$u_{xxxy} = (-2u_{-1,-1,0} + 2u_{-1,+1,0} + u_{-2,-1,0} - u_{-2,+1,0} + 2u_{+1,-1,0} - 2u_{+1,+1,0} - u_{+2,-1,0} + u_{+2,+1,0})/24$$
$$u_{xxxz} = (-2u_{-1,0,-1} + 2u_{-1,0,+1} + u_{-2,0,-1} - u_{-2,0,+1} + 2u_{+1,0,-1} - 2u_{+1,0,+1} - u_{+2,0,-1} + u_{+2,0,+1})/24$$
$$u_{xyyy} = (-2u_{-1,-1,0} + u_{-1,-2,0} + 2u_{-1,+1,0} - u_{-1,+2,0} + 2u_{+1,-1,0} - u_{+1,-2,0} - 2u_{+1,+1,0} + u_{+1,+2,0})/24$$
$$u_{yyyz} = (-2u_{0,-1,-1} + 2u_{0,-1,+1} + u_{0,-2,-1} - u_{0,-2,+1} + 2u_{0,+1,-1} - 2u_{0,+1,+1} - u_{0,+2,-1} + u_{0,+2,+1})/24$$
$$u_{xzzz} = (-2u_{-1,0,-1} + u_{-1,0,-2} + 2u_{-1,0,+1} - u_{-1,0,+2} + 2u_{+1,0,-1} - u_{+1,0,-2} - 2u_{+1,0,+1} + u_{+1,0,+2})/24$$
$$u_{yzzz} = (-2u_{0,-1,-1} + u_{0,-1,-2} + 2u_{0,-1,+1} - u_{0,-1,+2} + 2u_{0,+1,-1} - u_{0,+1,-2} - 2u_{0,+1,+1} + u_{0,+1,+2})/24$$
$$u_{xxyy} = (4u_{0,0,0} - 2u_{0,-1,0} - 2u_{0,+1,0} - 2u_{-1,0,0} + u_{-1,-1,0} + u_{-1,+1,0} - 2u_{+1,0,0} + u_{+1,-1,0} + u_{+1,+1,0})/4$$
$$u_{yyzz} = (4u_{0,0,0} - 2u_{0,0,-1} - 2u_{0,0,+1} - 2u_{0,-1,0} + u_{0,-1,-1} + u_{0,-1,+1} - 2u_{0,+1,0} + u_{0,+1,-1} + u_{0,+1,+1})/4$$
$$u_{xxzz} = (4u_{0,0,0} - 2u_{0,0,-1} - 2u_{0,0,+1} - 2u_{-1,0,0} + u_{-1,0,-1} + u_{-1,0,+1} - 2u_{+1,0,0} + u_{+1,0,-1} + u_{+1,0,+1})/4$$

$$u_{xxyz} = (2u_{0,-1,-1} - 2u_{0,-1,+1} - 2u_{0,+1,-1} + 2u_{0,+1,+1} - 2u_{-1,-1,-1} + 2u_{-1,-1,+1} + 2u_{-1,+1,-1} - 2u_{-1,+1,+1}$$
$$+ u_{-2,-1,-1} - u_{-2,-1,+1} - u_{-2,+1,-1} + u_{-2,+1,+1} - 2u_{+1,-1,-1} + 2u_{+1,-1,+1} + 2u_{+1,+1,-1} - 2u_{+1,+1,+1} + u_{+2,-1,-1}$$
$$- u_{+2,-1,+1} - u_{+2,+1,-1} + u_{+2,+1,+1})/16$$

$$u_{xyyz} = (2u_{-1,0,-1} - 2u_{-1,0,+1} - 2u_{-1,-1,-1} + 2u_{-1,-1,+1} + u_{-1,-2,-1} - u_{-1,-2,+1} - 2u_{-1,+1,-1} + 2u_{-1,+1,+1}$$
$$+ u_{-1,+2,-1} - u_{-1,+2,+1} - 2u_{+1,0,-1} + 2u_{+1,0,+1} + 2u_{+1,-1,-1} - 2u_{+1,-1,+1} - u_{+1,-2,-1} + u_{+1,-2,+1} + 2u_{+1,+1,-1}$$
$$- 2u_{+1,+1,+1} - u_{+1,+2,-1} + u_{+1,+2,+1})/16$$

$$u_{xyzz} = (2u_{-1,-1,0} - 2u_{-1,-1,-1} + u_{-1,-1,-2} - 2u_{-1,-1,+1} + u_{-1,-1,+2} - 2u_{-1,+1,0} + 2u_{-1,+1,-1} - u_{-1,+1,-2} + 2u_{-1,+1,+1}$$
$$- u_{-1,+1,+2} - 2u_{+1,-1,0} + 2u_{+1,-1,-1} - u_{+1,-1,-2} + 2u_{+1,-1,+1} - u_{+1,-1,+2} + 2u_{+1,+1,0} - 2u_{+1,+1,-1} + u_{+1,+1,-2}$$
$$- 2u_{+1,+1,+1} + u_{+1,+1,+2})/16$$

The non-linear WENO hybridization follows the standard WENO-AO approach. This gives us a high order interpolating polynomial within the zone of interest where the WENO has retained all the higher order modes when the solution is smooth and higher order. If the solution is non-smooth then the WENO gives us the best lower order polynomial that is non-oscillatory. The variable $u_0$ then gives us the volume averaged value in the zone of interest.

**Appendix C : A WENO-Stabilized Transcription Strategy from 2D Point Values to Area Averages**

In 2D, we non-linearly hybridize between five smaller third order stencils and a larger fifth order stencil. The third order stencils give rise to interpolating polynomials that have the form

$$u(x,y) = u_0 + u_x x + u_y y + u_{xx}(x^2 - 1/12) + u_{yy}(y^2 - 1/12) + u_{xy} xy \ . \tag{C.1}$$

and a smoothness indicator of the form

$$IS = u_x^2 + 13 u_{xx}^2 / 3 + u_y^2 + 13 u_{yy}^2 / 3 + 7 u_{xy}^2 / 6 \ . \tag{C.2}$$

The large fifth order stencil gives rise to an interpolating polynomial of the form

$$\begin{aligned} u(x,y) &= u_0 + u_x x + u_y y + u_{xx}(x^2 - 1/12) + u_{yy}(y^2 - 1/12) + u_{xy} xy \\ &+ u_{xxx}(x^3 - 3x/20) + u_{yyy}(y^3 - 3y/20) + u_{xxy}(x^2 - 1/12) y + u_{xyy} x(y^2 - 1/12) \\ &+ u_{xxxx}(x^4 - 3x^2/14 + 3/560) + u_{yyyy}(y^4 - 3y^2/14 + 3/560) + u_{xxxy}(x^3 - 3x/20) y \\ &+ u_{xyyy} x(y^3 - 3y/20) + u_{xxyy}(x^2 - 1/12)(y^2 - 1/12) \ . \end{aligned} \tag{C.3}$$

and a smoothness indicator of the form

$$
\begin{aligned}
IS = & (u_x + u_{xxx}/10)^2 + \frac{13}{3}(u_{xx} + 123 u_{xxxx}/455)^2 + \frac{781}{20} u_{xxx}^2 \\
& + \frac{1421461}{2275} u_{xxxx}^2 + (u_y + u_{yyy}/10)^2 + \frac{13}{3}(u_{yy} + 123 u_{yyyy}/455)^2 \\
& + \frac{781}{20} u_{yyy}^2 + \frac{1421461}{2275} u_{yyyy}^2 + \frac{7}{6}(u_{xy} + 13 u_{xxxy}/140 + 13 u_{xyyy}/140)^2 \\
& + \frac{47}{10} u_{xxy}^2 + \frac{47}{10} u_{xyy}^2 + \frac{88841}{2100}(u_{xxxy}^2 + u_{xyyy}^2) + \frac{1}{16800}(u_{xxxy} - u_{xyyy})^2 + \frac{5083}{270} u_{xxyy}^2 \, .
\end{aligned}
\quad (C.4)
$$

Please note that this is not a finite volume style reconstruction but rather a pointwise interpolation. The four smaller third order interpolating stencils are one-sided and they can be identified on the mesh. We shall not describe all of them, but we will provide all the coefficients for one of them below:-

Stencil-1 (Biased in 1$^{st}$ quadrant, $x, y > 0$):

$$u_0 = (-2u_{0,+1} + u_{0,+2} + 26 u_{0,0} - 2 u_{+1,0} + u_{+2,0})/24$$
$$u_x = (-3u_{0,0} + 4 u_{+1,0} - u_{+2,0})/2 \; ; \; u_y = (4 u_{0,+1} - u_{0,+2} - 3 u_{0,0})/2$$
$$u_{xx} = (u_{0,0} - 2 u_{+1,0} + u_{+2,0})/2 \; ; \; u_{yy} = (-2 u_{0,+1} + u_{0,+2} + u_{0,0})/2 \; ; \; u_{xy} = -u_{0,+1} + u_{0,0} - u_{+1,0} + u_{+1,+1}$$

These four one-sided stencils have to be non-linearly hybridized with a smaller third order central stencil for which we provide all the coefficients below:-

$$u_0 = (20 u_{0,0} + u_{0,-1} + u_{0,+1} + u_{-1,0} + u_{+1,0})/24$$
$$u_x = (-u_{-1,0} + u_{+1,0})/2 \; ; \; u_y = (-u_{0,-1} + u_{0,+1})/2$$
$$u_{xx} = (-2 u_{0,0} + u_{-1,0} + u_{+1,0})/2 \; ; \; u_{yy} = (-2 u_{0,0} + u_{0,-1} + u_{0,+1})/2 \; ; \; u_{xy} = (u_{-1,-1} - u_{-1,+1} - u_{+1,-1} + u_{+1,+1})/4$$

In addition, we will need a larger fifth order interpolating stencil against which we have to perform non-linear hybridization via a vis the five smaller third order stencil. For the large fifth order stencil we provide all the coefficients below:-

$$u_0 = (4636u_{0,0} + 288u_{0,-1} - 17u_{0,-2} + 288u_{0,+1} - 17u_{0,+2} + 288u_{-1,0} + 10u_{-1,-1} + 10u_{-1,+1} - 17u_{-2,0}$$
$$+ 288u_{+1,0} + 10u_{+1,-1} + 10u_{+1,+1} - 17u_{+2,0})/5760$$

$$u_x = (-144u_{-1,0} - 5u_{-1,-1} - 5u_{-1,+1} + 17u_{-2,0} + 144u_{+1,0} + 5u_{+1,-1} + 5u_{+1,+1} - 17u_{+2,0})/240$$

$$u_y = (-144u_{0,-1} + 17u_{0,-2} + 144u_{0,+1} - 17u_{0,+2} - 5u_{-1,-1} + 5u_{-1,+1} - 5u_{+1,-1} + 5u_{+1,+1})/240$$

$$u_{xx} = (-374u_{0,0} - 14u_{0,-1} - 14u_{0,+1} + 198u_{-1,0} + 7u_{-1,-1} + 7u_{-1,+1} - 11u_{-2,0} + 198u_{+1,0} + 7u_{+1,-1} + 7u_{+1,+1}$$
$$- 11u_{+2,0})/336$$

$$u_{yy} = (-374u_{0,0} + 198u_{0,-1} - 11u_{0,-2} + 198u_{0,+1} - 11u_{0,+2} - 14u_{-1,0} + 7u_{-1,-1} + 7u_{-1,+1} - 14u_{+1,0} + 7u_{+1,-1}$$
$$+ 7u_{+1,+1})/336$$

$$u_{xy} = (188u_{-1,-1} - 17u_{-1,-2} - 188u_{-1,+1} + 17u_{-1,+2} - 17u_{-2,-1} + 17u_{-2,+1} - 188u_{+1,-1} + 17u_{+1,-2} + 188u_{+1,+1}$$
$$- 17u_{+1,+2} + 17u_{+2,-1} - 17u_{+2,+1})/480$$

$$u_{xxx} = (2u_{-1,0} - u_{-2,0} - 2u_{+1,0} + u_{+2,0})/12$$

$$u_{yyy} = (2u_{0,-1} - u_{0,-2} - 2u_{0,+1} + u_{0,+2})/12$$

$$u_{xxy} = (2u_{0,-1} - 2u_{0,+1} - u_{-1,-1} + u_{-1,+1} - u_{+1,-1} + u_{+1,+1})/4$$

$$u_{xyy} = (2u_{-1,0} - u_{-1,-1} - u_{-1,+1} - 2u_{+1,0} + u_{+1,-1} + u_{+1,+1})/4$$

$$u_{xxxx} = (6u_{0,0} - 4u_{-1,0} + u_{-2,0} - 4u_{+1,0} + u_{+2,0})/24$$

$$u_{yyyy} = (6u_{0,0} - 4u_{0,-1} + u_{0,-2} - 4u_{0,+1} + u_{0,+2})/24$$

$$u_{xxxy} = (-2u_{-1,-1} + 2u_{-1,+1} + u_{-2,-1} - u_{-2,+1} + 2u_{+1,-1} - 2u_{+1,+1} - u_{+2,-1} + u_{+2,+1})/24$$

$$u_{xyyy} = (-2u_{-1,-1} + u_{-1,-2} + 2u_{-1,+1} - u_{-1,+2} + 2u_{+1,-1} - u_{+1,-2} - 2u_{+1,+1} + u_{+1,+2})/24$$

$$u_{xxyy} = (4u_{0,0} - 2u_{0,-1} - 2u_{0,+1} - 2u_{-1,0} + u_{-1,-1} + u_{-1,+1} - 2u_{+1,0} + u_{+1,-1} + u_{+1,+1})/4$$

The non-linear WENO hybridization follows the standard WENO-AO approach. This gives us a high order interpolating polynomial within the zone of interest where the WENO has retained all the higher order modes when the solution is smooth and higher order. If the solution is non-smooth then the WENO gives us the best lower order polynomial that is non-oscillatory. The variable $u_0$ then gives us the volume averaged value in the zone of interest.

**Figure Captions**

*Fig. 1 shows us that the primal variables of the scheme, given by the normal components of the magnetic induction, are facially-collocated. The components of the primal magnetic field vector are shown by the thick blue arrows. The overbars on the magnetic field components indicate that these are facially averaged. They undergo an update from the induction equation. The edge-collocated electric fields, which are used for updating the facial magnetic induction components, are shown by the thin blue arrows close to the appropriate edge. They too have overbars to indicate that they are edge-averaged. The superscript "num" for the electric field components indicates that they are multidimensionally stabilized and, therefore, suitable for use in a numerical scheme.*

*Fig. 2 shows four zones in the xy-plane that come together at the z-edge of a three-dimensional mesh. Since the mesh is viewed from the top in plan view, the z-edge is shown by the black dot and the four abutting zones are shown as four squares. The four incoming states have subscripts given by "RU" for right-upper; "LU" for left-upper; "LD" for left-down and "RD" for right-down. Fig 2 shows the situation before the states start interacting via four one-dimensional and one multidimensional Riemann problems. The black arrows indicate that higher-order 2D interpolation can eventually be used to obtain the centered part of the electric field at the z-edge. The blue arrows denote the normal components of the magnetic field at the zone faces. Owing to the constraint, these field components (blue arrows) are continuous across zone faces. They can, therefore, be obtained from the higher order facial reconstruction of the normal component of the magnetic field within each face. This provides higher order values of the x- and y-components of the magnetic field at the z-edge that minimize the dissipation terms.*

*Fig. 3 shows the same situation as Fig. 2. However, it shows the situation after the four incoming states start interacting with each other. Four one-dimensional Riemann problems, shown by dashed lines, develop between the four pairs of incoming states. The resolved states from the one-dimensional Riemann problems are shown by a superscript with a single star. The shaded region depicts the strongly interacting state that arises when the four one-dimensional Riemann problems interact with one another. The strongly interacting state is shown by a superscript with a double star. We want to find the z-component of the electric field in the strongly interacting state. This*

*gives us the z-component of the electric field at the z-edge, which is shown by the dot in this two-dimensional projection.*

*Figs. 4a to 4d show the four cases that are fully supersonic in both directions. The axes, in similarity variables, are shown in blue. The wave model is shown in black.*

*Figs. 4e to 4h show the four cases that are fully supersonic in only one of the two directions. The axes, in similarity variables, are shown in blue. The wave model is shown in black.*

*Fig. 5 shows the four one-sided third order and one central third order stencil for the 2D-WENO interpolation. The $(i,j)^{th}$ computational cell is denoted by the shaded region. The zone-centered point values of the neighboring zones are denoted by "centered-dot".*

*Fig. 6 shows the higher order central stencils for the 2D-WENO interpolation. Fig. 6a shows the $5^{th}$ order accurate central stencil, Fig. 6b shows the $7^{th}$ order accurate central stencil and Fig. 6c shows the $9^{th}$ order accurate central stencil. The $(i,j)^{th}$ computational cell is denoted by the shaded region. The zone-centered point values of the neighboring zones are denoted by "centered-dot".*

*Fig. 7 shows part of the z-edge of the mesh. The points that make up the edge-centers of the mesh are shown by the solid black dots. The high order accurate pointwise values of $e_z^{num}$ have been evaluated at the edge-centers of the mesh by using the two-dimensional Riemann solver at each edge-center. We want the high order integral $\bar{E}_z^{num}$ for values of "z" in the range [-1/2,1/2]. This edge of interest is also identified in green. The figure also shows the one-dimensional stencils associated with the edge of interest for the third and fifth order WENO-AO interpolation schemes. We have three smaller third order stencils and a large fifth order stencil. We want a higher order line integral of $e_z^{num}$ along the green edge of interest in the figure. Once a high order, one-dimensional WENO interpolation polynomial has been evaluated using Legendre bases in the edge of interest, the leading term of that polynomial will give us the line integral with fifth order of accuracy. The interested reader should compare this figure with Fig. 1 from Balsara et al. [36] in*

*order to realize that the solution to the problem described here is already available from Section 3 of the previously cited paper.*

*Fig. 8a and 8b are analogous to Fig. 2 because they show four zones in the xy-plane that come together at the z-edge of a three-dimensional mesh. Fig. 8a shows the four zones that surround a z-edge. It shows how the "θ" variables that are evaluated at the xz- and yz-faces can be used to form an effective "θ" at the z-edge of the mesh. This effective "θ" at the z-edge can then be used to lower the order of the edge-centered z-component of the electric field that is used in the update of the facial magnetic fields in the xz- and yz-faces. Like Fig. 2, Fig. 8b shows the inputs that go into the evaluation of the z-component of the electric field. The only difference from Fig. 2 is that we now have the option of making a high order evaluation (which uses all the high order WENO reconstructions and interpolations as described in the text) which is superscripted with "HO"; and a low order (first order) evaluation which is superscripted with "LO".*

*Fig. 9 is analogous to Fig. 1 because it shows the components of the magnetic field in the faces of the mesh. The difference from Fig. 1 is that within each face we now have a high order component which is superscripted with "HO"; and a low order component which is superscripted with "LO". Both components in each face have been advanced in time using a forward Euler scheme with a timestep $\Delta t$. This temporal advance of the facial magnetic field components is made before a call to the AFD-WENO routine, and both the components within each face will be used for the PCP update.*

*Fig. 10) CED: Refraction of a compact electromagnetic beam by a dielectric slab. Figs. a), b), and c) show $B_z$, $D_x$, and $D_y$ at the initial time. Figs. d), e), and f) show the same at a final time of $4 \times 10^{-14}$s. The surface of the dielectric slab is identified by the dashed vertical black line. The oblique dashed black lines demarcate the angle of incidence, the angle of refraction, and the angle of reflection.*

*Fig. 11) CED: Total internal reflection of a compact electromagnetic beam by a dielectric slab. Figs. a), b), and c) show $B_z$, $D_x$, and $D_y$ at the initial time. Figs. d), e), and f) show the same at a final time of $5 \times 10^{-14}$s. The surface of the dielectric slab is identified by the vertical black line. The oblique black lines demarcate the angle of incidence and the angle of total internal reflection.*

*Fig. 12) MHD: Orszag-Tang Problem. Fig. 12a shows the density, Fig. 12b shows the pressure, Fig. 12c shows magnitude of the velocity and Fig. 12d shows magnitude of the magnetic field vector using the ninth-order accurate scheme.*

*Fig. 13) MHD: Rotor Problem. Fig. 13a shows the density, Fig. 13b shows the pressure, Fig. 13c shows magnitude of the velocity and Fig. 13d shows magnitude of the magnetic field vector using the seventh-order accurate scheme.*

*Fig. 14) MHD: Blast Problem-I (BLAST-I). Fig. 14a shows the density, Fig. 14b shows the pressure, Fig. 14c shows magnitude of the velocity and Fig. 14d shows magnitude of the magnetic field vector using the fifth-order accurate scheme. Here we show the three-dimensional variant of this problem.*

*Fig. 15) MHD: Blast Problem-II (BLAST-II). Fig. 15a shows the density, Fig. 15b shows the pressure, Fig. 15c shows magnitude of the velocity and Fig. 15d shows magnitude of the magnetic field vector using the fifth-order accurate scheme.*

*Fig. 16) MHD: Astrophysical Jet Problem (Jet-I,II and III). Figs. 16a,b,c show the resulting density profiles on logarithmic scales for the Jet-I, Jet-II and Jet-III problems, respectively. Seventh-order accurate schemes is used.*

*Fig. 17) MHD: Astrophysical Jet Problem (Jet-I,II and III). Figs. 17a,b,c show the resulting pressure profiles on logarithmic scales for the Jet-I, Jet-II and Jet-III problems, respectively. Seventh-order accurate schemes is used.*

*Fig. 18) RMHD: Blast Problem. Fig. 18a shows the logarithm of density, Fig. 18b shows the logarithm of pressure, and Fig. 18c shows the magnetic pressure using the fifth-order accurate scheme.*

*Fig. 19) RMHD: Orszag-Tang Problem. Fig. 19a, b show the logarithm of density and the logarithm of magnetic pressure at time t=2.818127; and Fig. 19c, d show the logarithm of density and the logarithm of magnetic pressure at time t=6.8558. The seventh-order accurate scheme is used.*

*Fig. 20) RMHD: Shock-Cloud interaction problem. Fig. 20a shows the logarithm of density and Fig. 20b shows the logarithm of magnetic pressure using the ninth-order accurate scheme.*

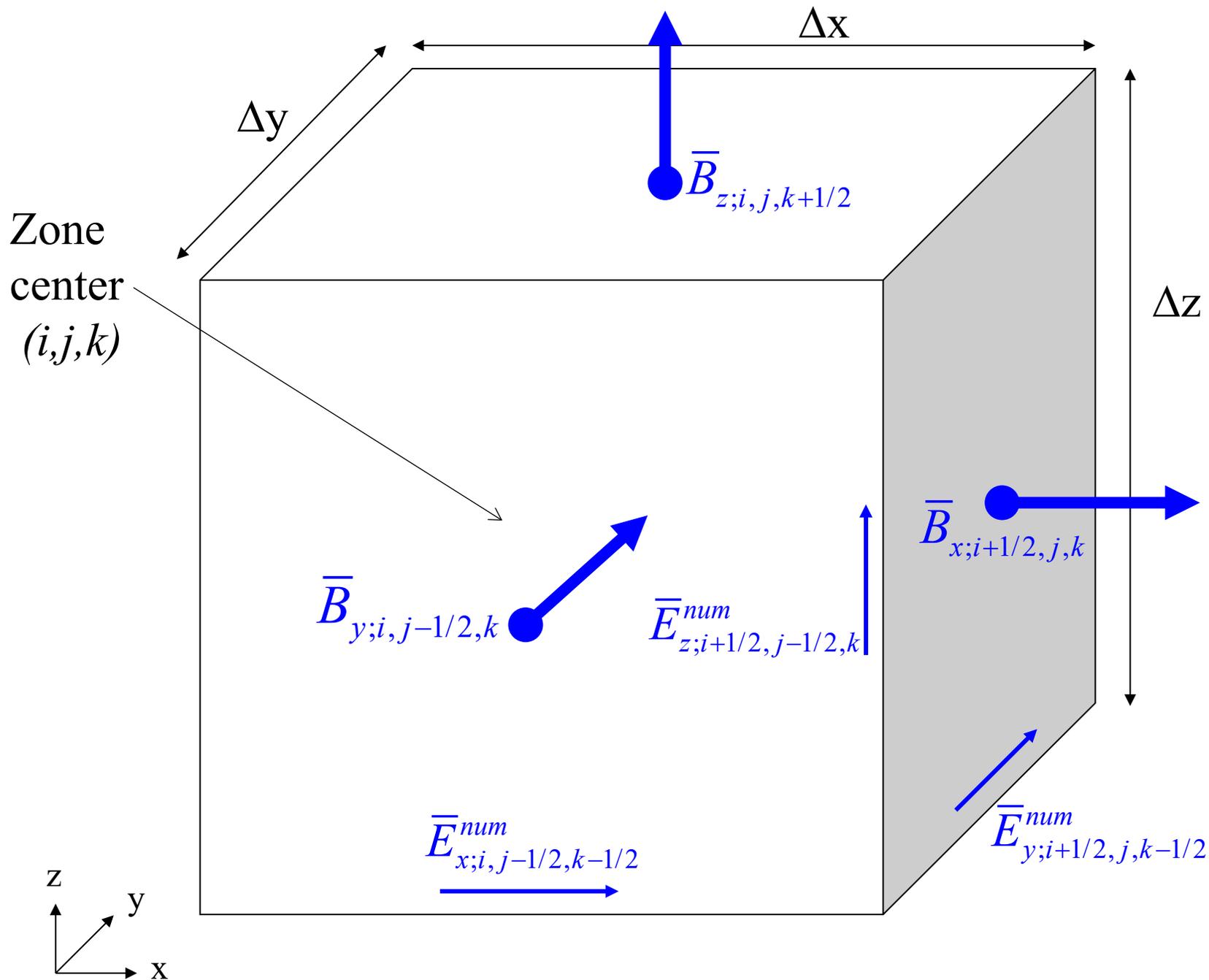

Fig. 1 shows us that the primal variables of the scheme, given by the normal components of the magnetic induction, are facially-collocated. The components of the primal magnetic field vector are shown by the thick blue arrows. The overbars on the magnetic field components indicate that these are facially averaged. They undergo an update from the induction equation. The edge-collocated electric fields, which are used for updating the facial magnetic induction components, are shown by the thin blue arrows close to the appropriate edge. They too have overbars to indicate that they are edge-averaged. The superscript "num" for the electric field components indicates that they are multidimensionally stabilized and, therefore, suitable for use in a numerical scheme.

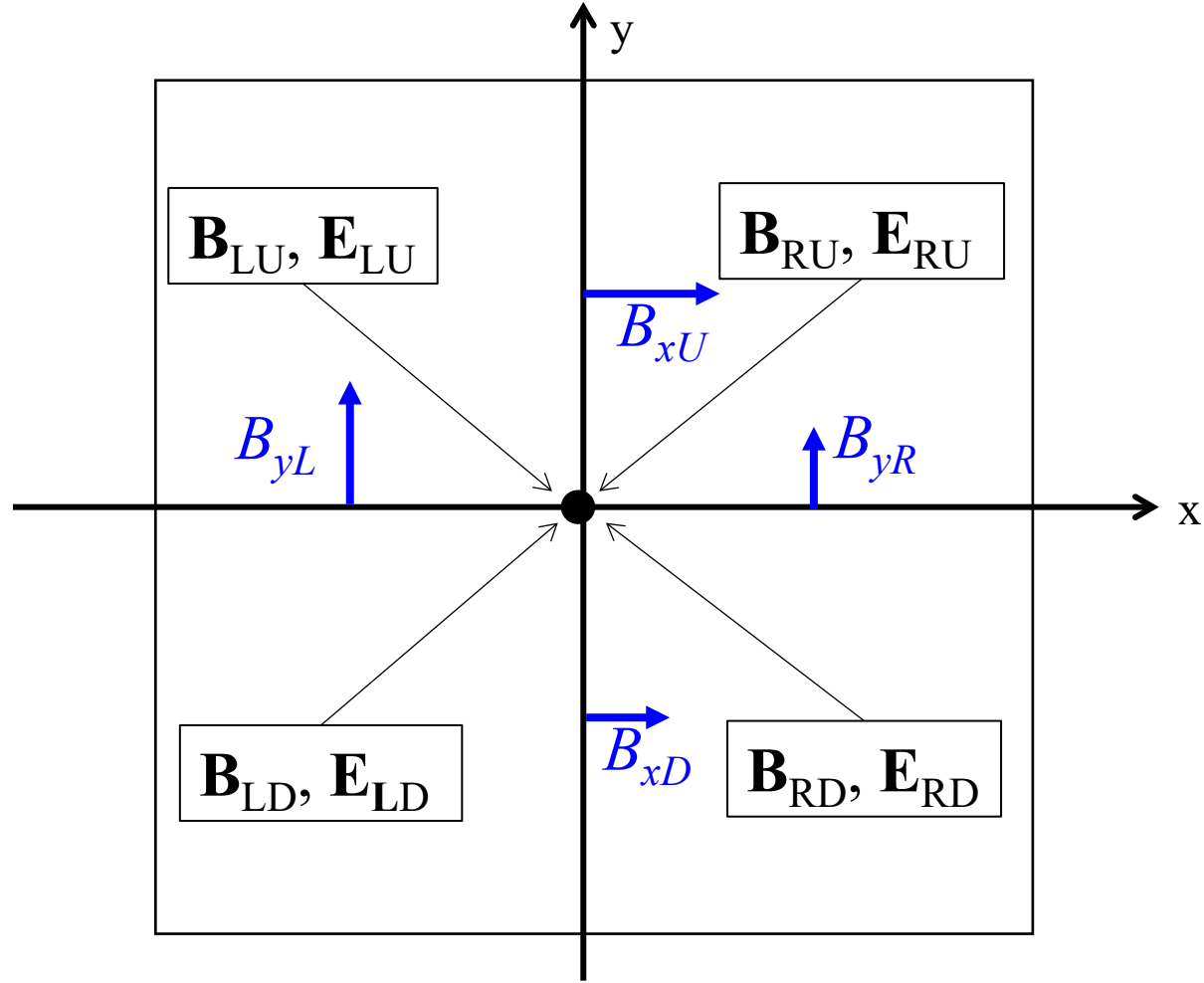

Fig. 2 shows four zones in the xy-plane that come together at the z-edge of a three-dimensional mesh. Since the mesh is viewed from the top in plan view, the z-edge is shown by the black dot and the four abutting zones are shown as four squares. The four incoming states have subscripts given by "RU" for right-upper; "LU" for left-upper; "LD" for left-down and "RD" for right-down. Fig 2 shows the situation before the states start interacting via four one-dimensional and one multidimensional Riemann problems. The black arrows indicate that higher-order 2D interpolation can eventually be used to obtain the centered part of the electric field at the z-edge. The blue arrows denote the normal components of the magnetic field at the zone faces. Owing to the constraint, these field components (blue arrows) are continuous across zone faces. They can, therefore, be obtained from the higher order facial reconstruction of the normal component of the magnetic field within each face. This provides higher order values of the x- and y-components of the magnetic field at the z-edge that minimize the dissipation terms.

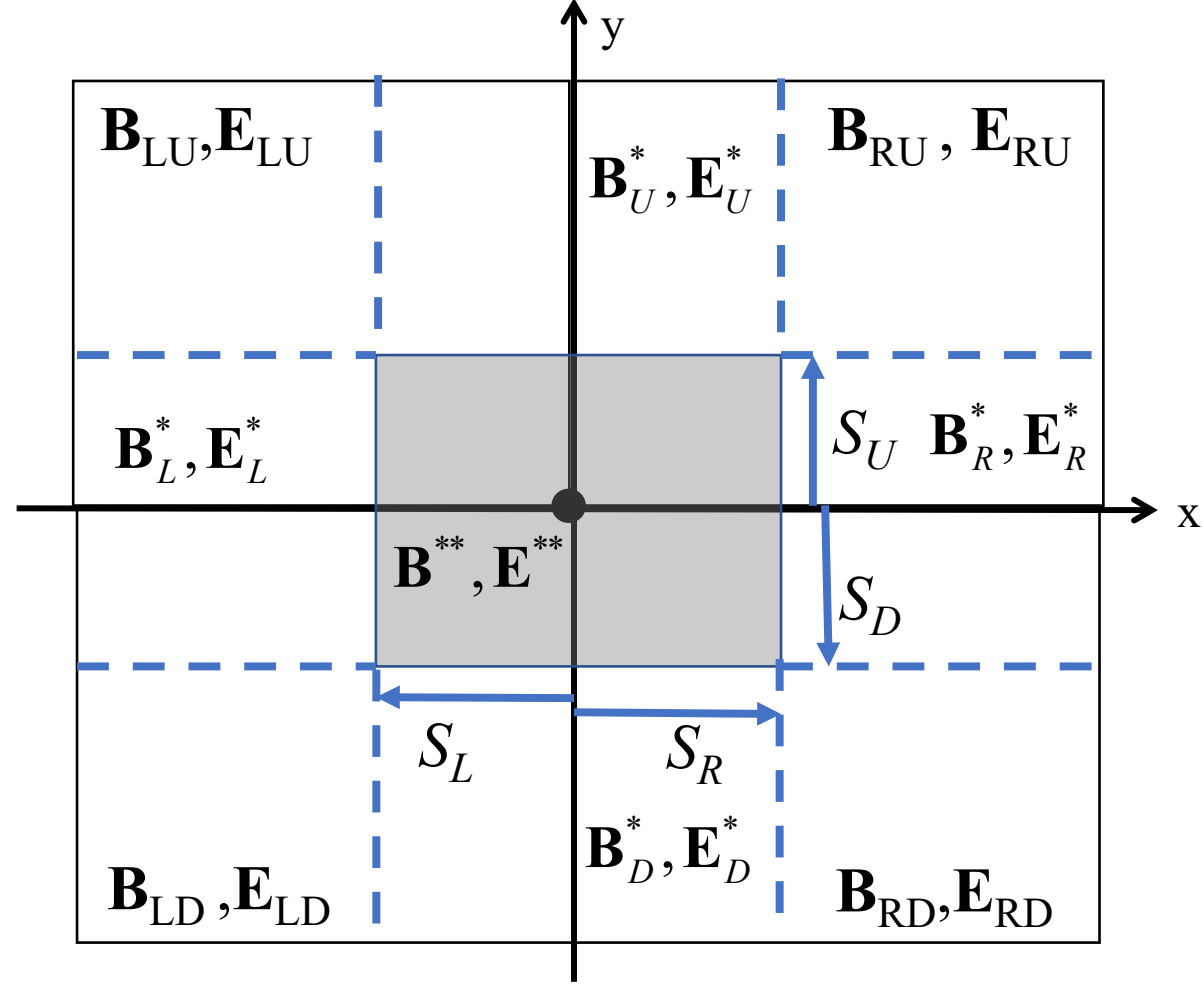

Fig. 3 shows the same situation as Fig. 2. However, it shows the situation after the four incoming states start interacting with each other. Four one-dimensional Riemann problems, shown by dashed lines, develop between the four pairs of incoming states. The resolved states from the one-dimensional Riemann problems are shown by a superscript with a single star. The shaded region depicts the strongly interacting state that arises when the four one-dimensional Riemann problems interact with one another. The strongly interacting state is shown by a superscript with a double star. We want to find the z-component of the electric field in the strongly interacting state. This gives us the z-component of the electric field at the z-edge, which is shown by the dot in this two-dimensional projection.

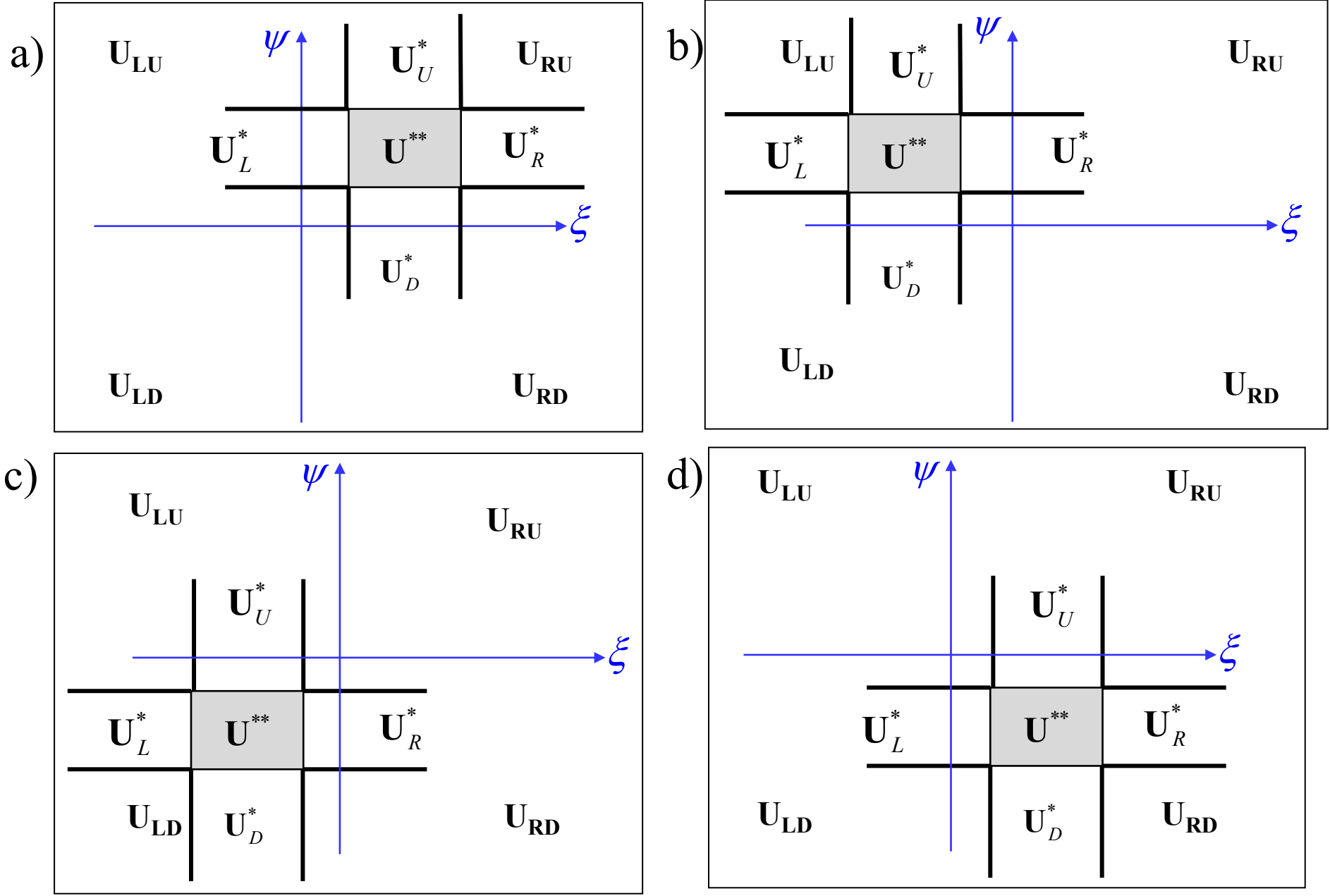

*Figs. 4a to 4d show the four cases that are fully supersonic in both directions. The axes, in similarity variables, are shown in blue. The wave model is shown in black.*

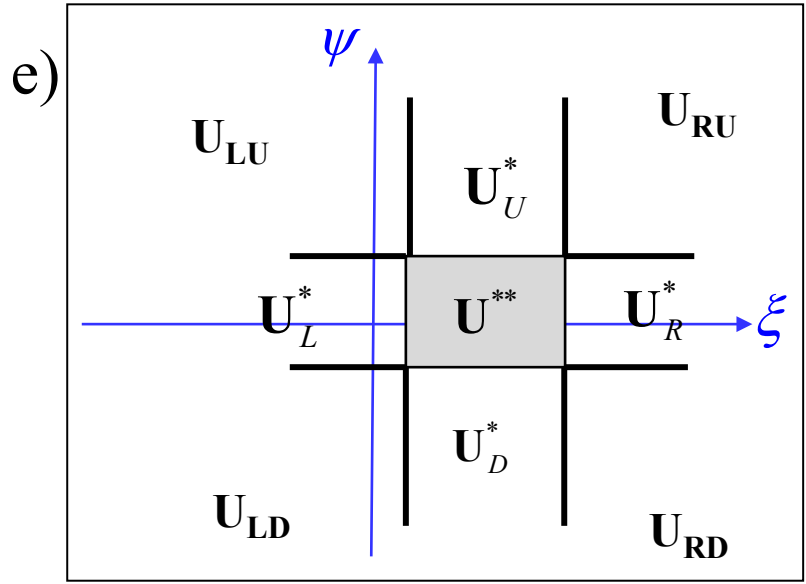
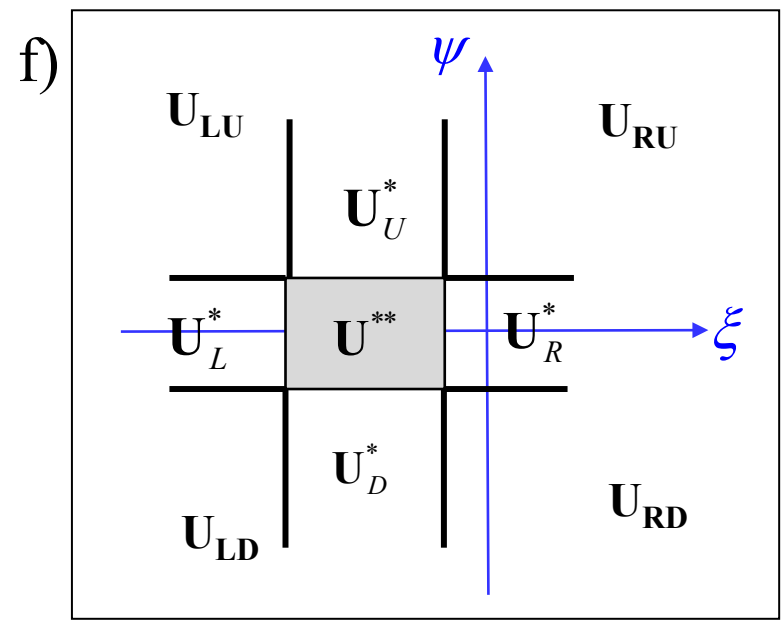
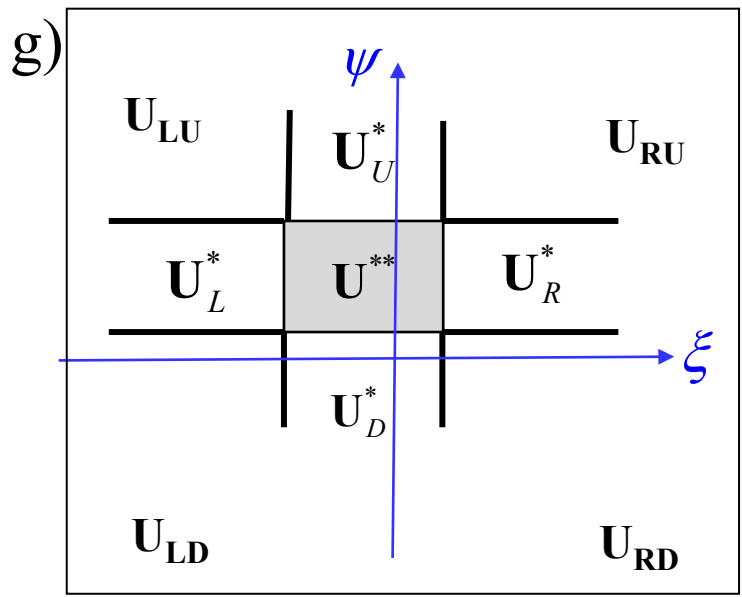
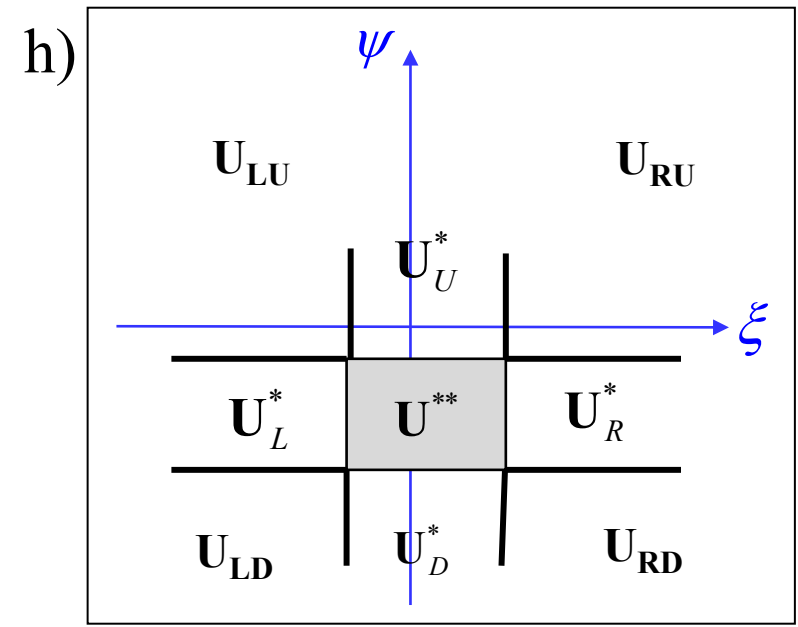

*Figs. 4e to 4h show the four cases that are fully supersonic in only one of the two directions. The axes, in similarity variables, are shown in blue. The wave model is shown in black.*

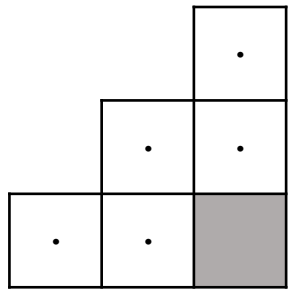 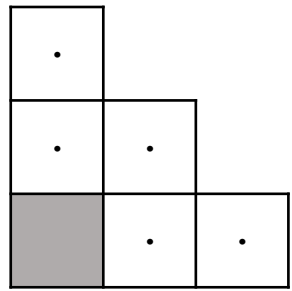 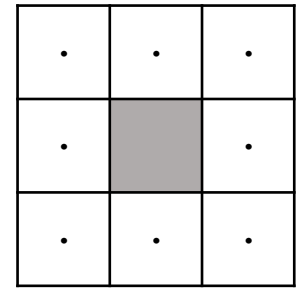

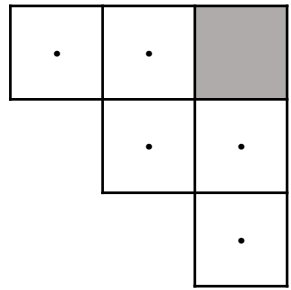 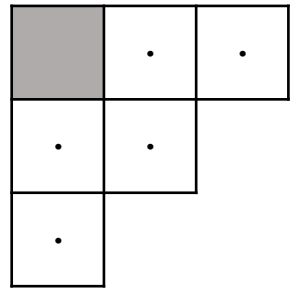

*Fig. 5 shows the four one-sided third order and one central third order stencil for the 2D-WENO interpolation. The (i,j)$^{th}$ computational cell is denoted by the shaded region. The zone-centered point values of the neighbouring zones are denoted by "centered-dot".*

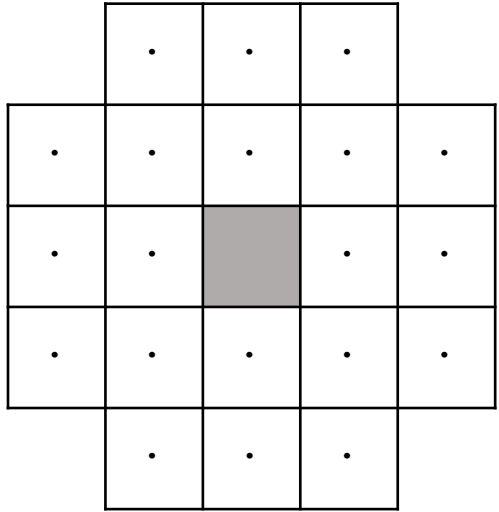 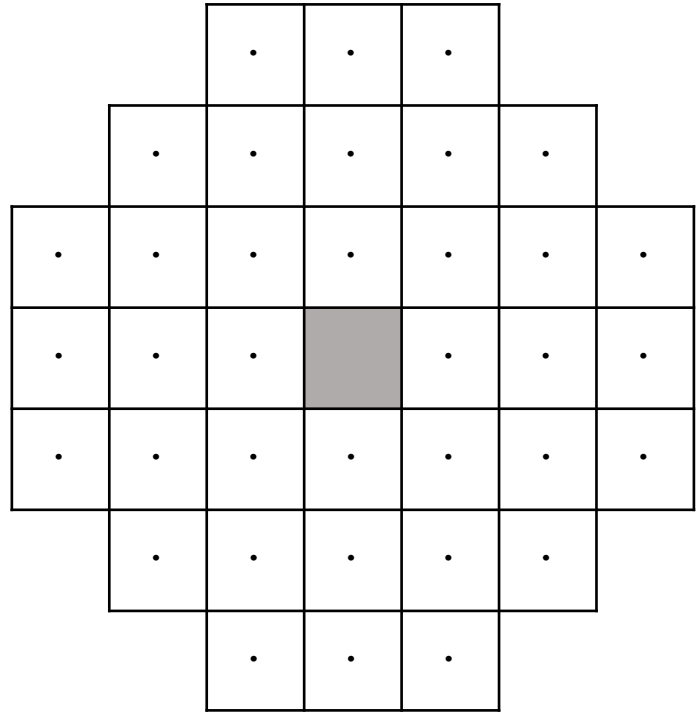 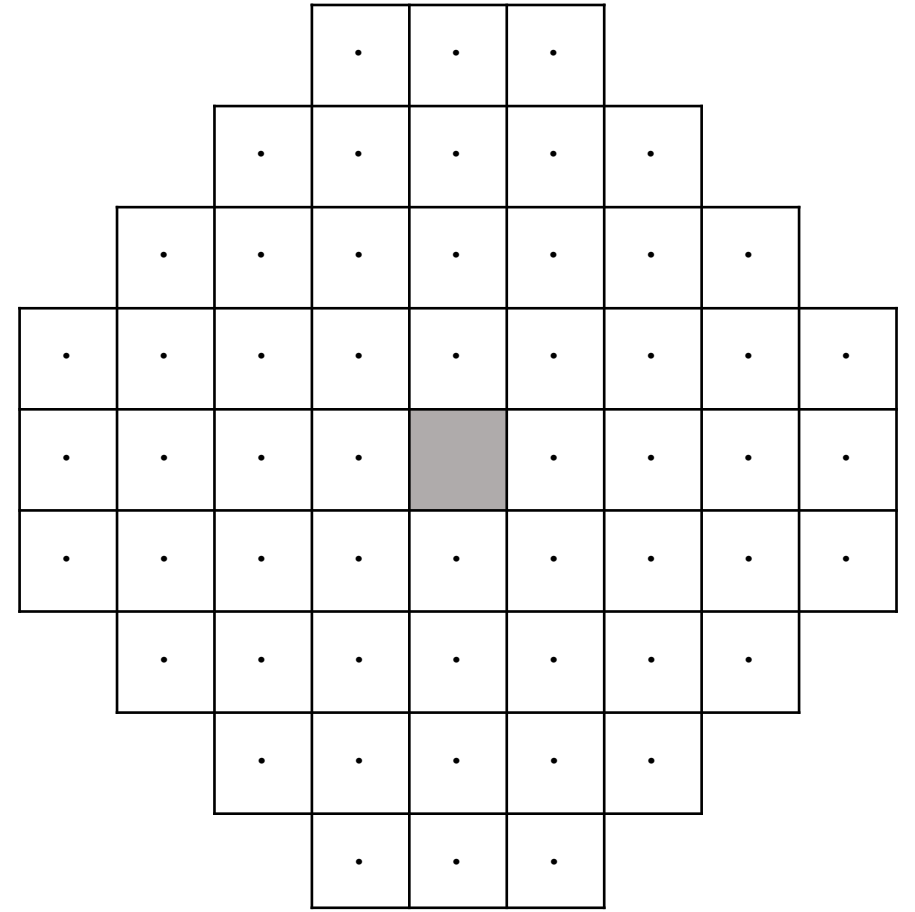

a) $S_1^{r_5}$      b) $S_1^{r_7}$      c) $S_1^{r_9}$

*Fig. 6 shows the higher order central stencils for the 2D-WENO interpolation. Fig. 6a shows the 5$^{th}$ order accurate central stencil, Fig. 6b shows the 7$^{th}$ order accurate central stencil and Fig. 6c shows the 9$^{th}$ order accurate central stencil. The (i,j)$^{th}$ computational cell is denoted by the shaded region. The zone-centered point values of the neighboring zones are denoted by "centered-dot".*

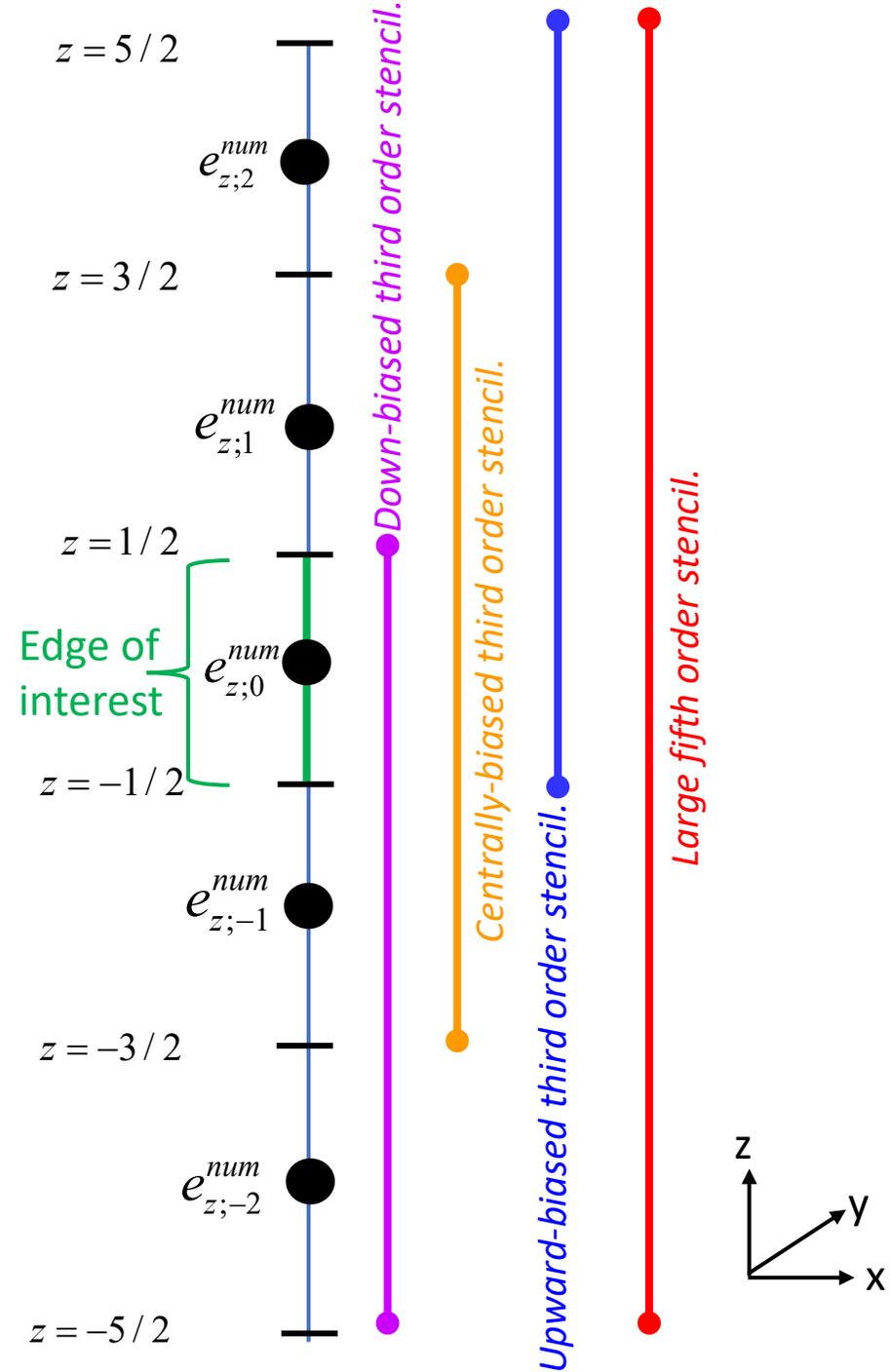

*Fig. 7 shows part of the z-edge of the mesh. The points that make up the edge-centers of the mesh are shown by the solid black dots. The high order accurate pointwise values of $e_z^{num}$ have been evaluated at the edge-centers of the mesh by using the two-dimensional Riemann solver at each edge-center. We want the high order integral $\bar{E}_z^{num}$ for values of "z" in the range [-1/2,1/2]. This edge of interest is also identified in green. The figure also shows the one-dimensional stencils associated with the edge of interest for the third and fifth order WENO-AO interpolation schemes. We have three smaller third order stencils and a large fifth order stencil. We want a higher order line integral of $e_z^{num}$ along the green edge of interest in the figure. Once a high order, one-dimensional WENO interpolation polynomial has been evaluated using Legendre bases in the edge of interest, the leading term of that polynomial will give us the line integral with fifth order of accuracy. The interested reader should compare this figure with Fig. 1 from Balsara et al. [32] in order to realize that the solution to the problem described here is already available from Section 3 of the previously cited paper.*

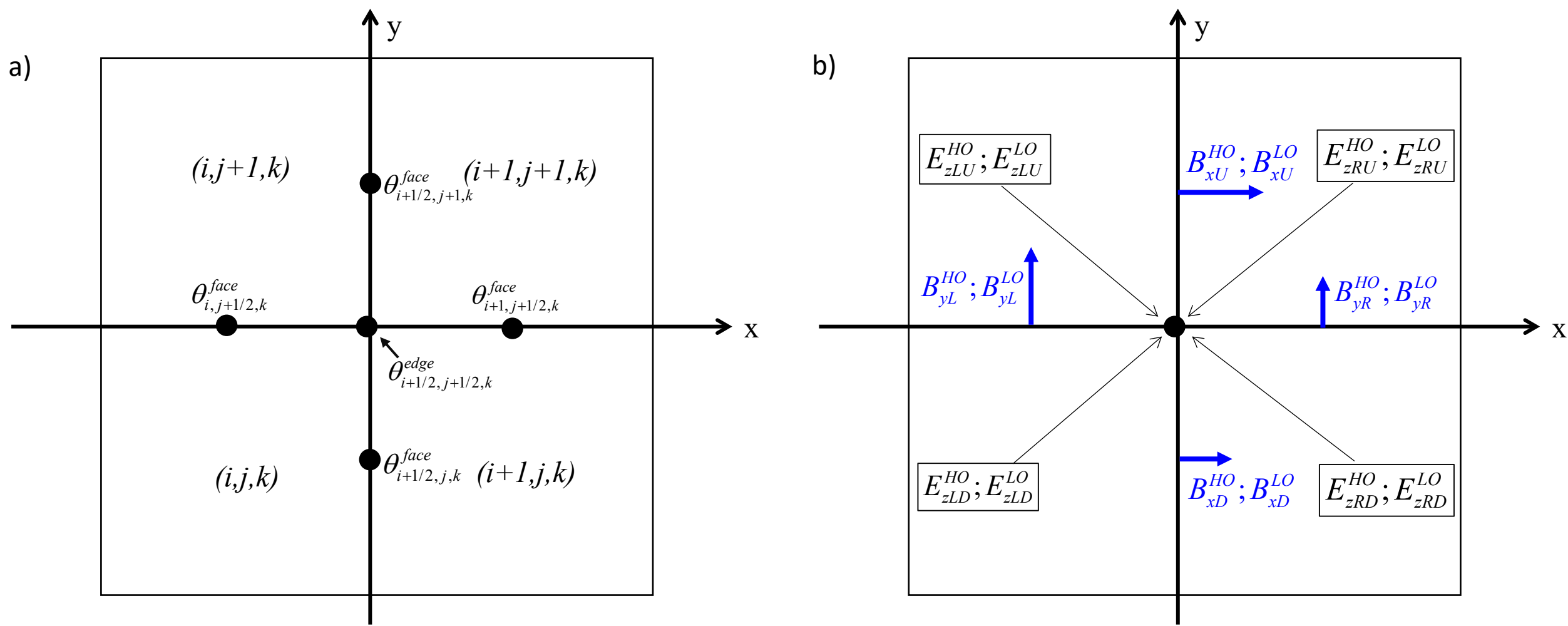

Fig. 8a and 8b are analogous to Fig. 2 because they show four zones in the xy-plane that come together at the z-edge of a three-dimensional mesh. Fig. 8a shows the four zones that surround a z-edge. It shows how the "θ" variables that are evaluated at the xz- and yz-faces can be used to form an effective "θ" at the z-edge of the mesh. This effective "θ" at the z-edge can then be used to lower the order of the edge-centered z-component of the electric field that is used in the update of the facial magnetic fields in the xz- and yz-faces. Like Fig. 2, Fig. 8b shows the inputs that go into the evaluation of the z-component of the electric field. The only difference from Fig. 2 is that we now have the option of making a high order evaluation (which uses all the high order WENO reconstructions and interpolations as described in the text) which is superscripted with "HO"; and a low order (first order) evaluation which is superscripted with "LO".

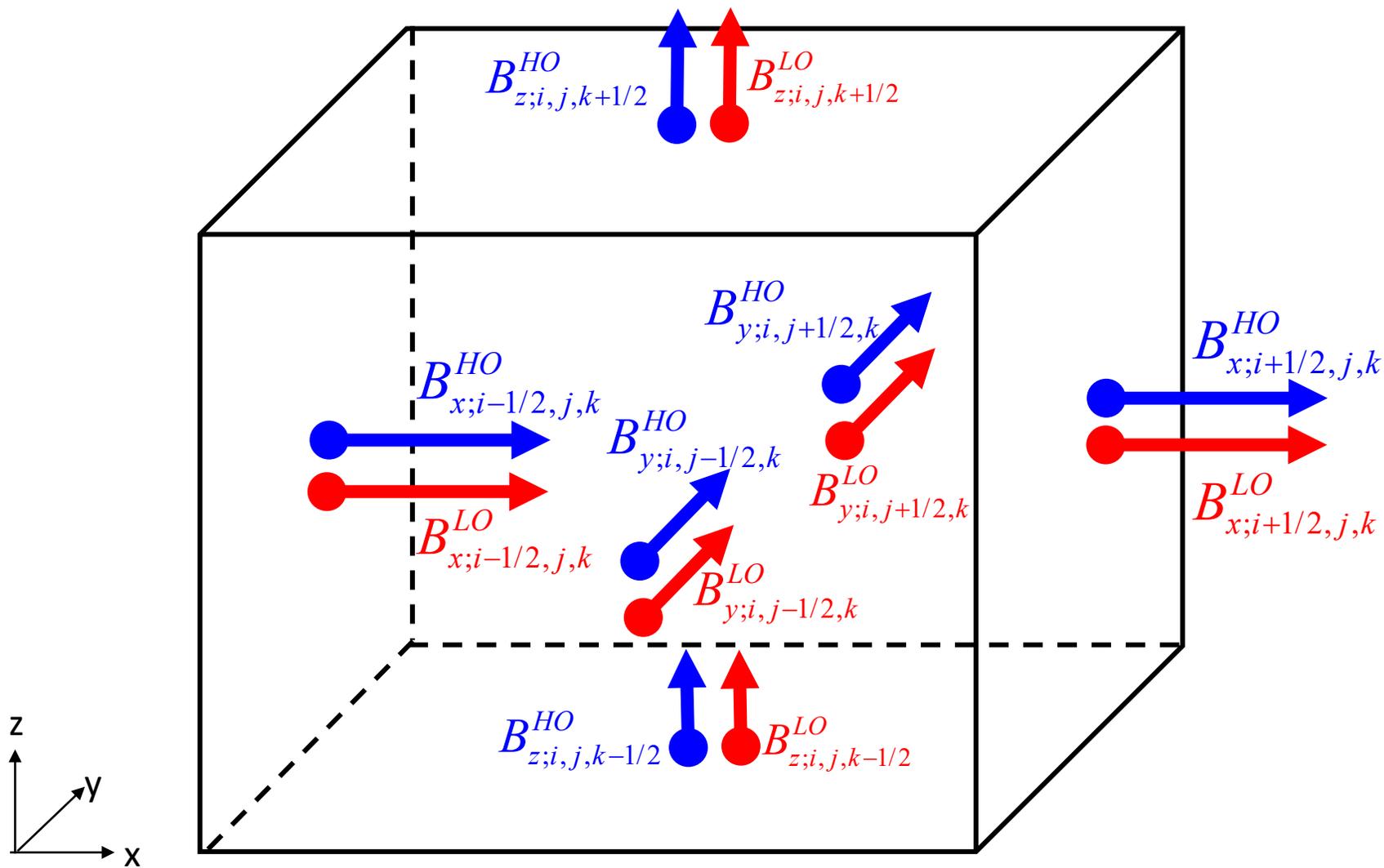

Fig. 9 is analogous to Fig. 1 because it shows the components of the magnetic field in the faces of the mesh. The difference from Fig. 1 is that within each face we now have a high order component which is superscripted with "HO"; and a low order component which is superscripted with "LO". Both components in each face have been advanced in time using a forward Euler scheme with a timestep $\Delta t$. This temporal advance of the facial magnetic field components is made before a call to the AFD-WENO routine, and both the components within each face will be used for the PCP update.

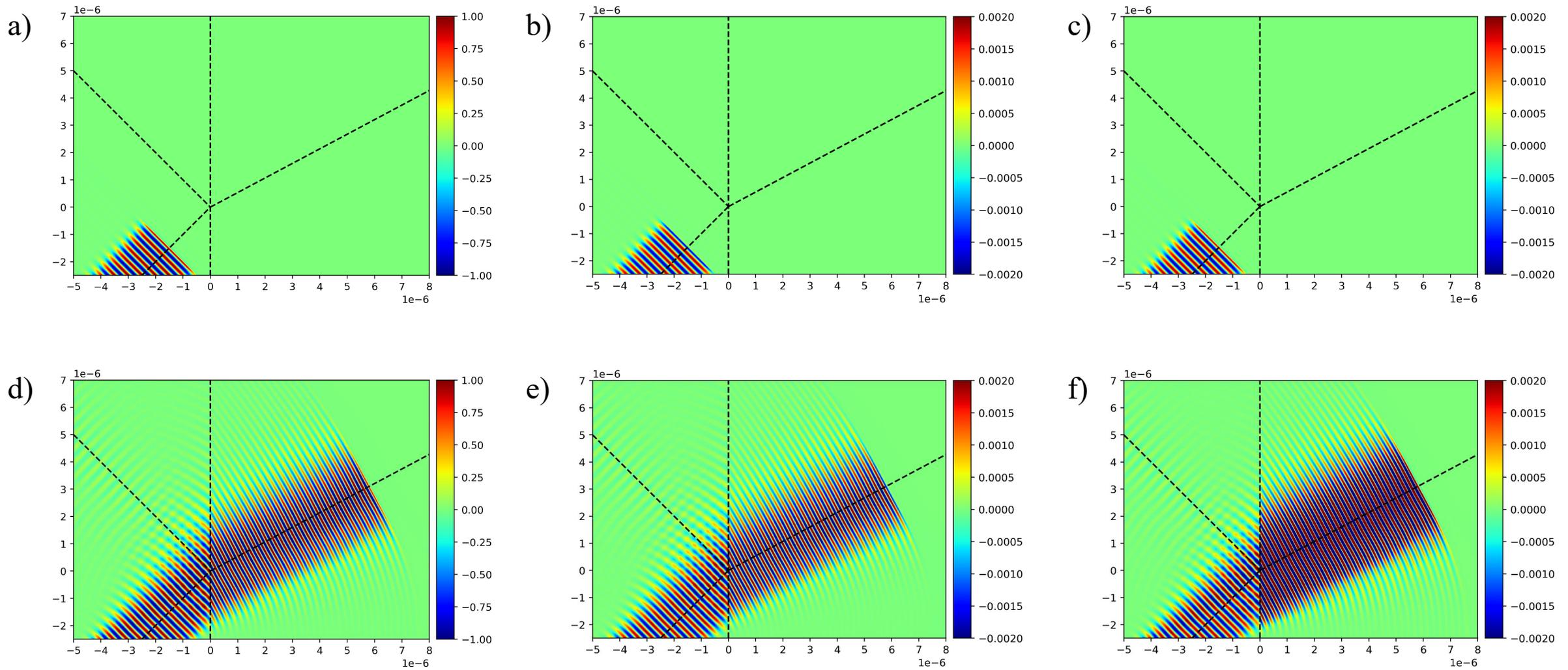

*Fig. 10) CED: Refraction of a compact electromagnetic beam by a dielectric slab. Figs. a), b), and c) show $B_z$, $D_x$, and $D_y$ at the initial time. Figs. d), e), and f) show the same at a final time of $4\times10^{-14}$s. The surface of the dielectric slab is identified by the dashed vertical black line. The oblique dashed black lines demarcate the angle of incidence, the angle of refraction, and the angle of reflection.*

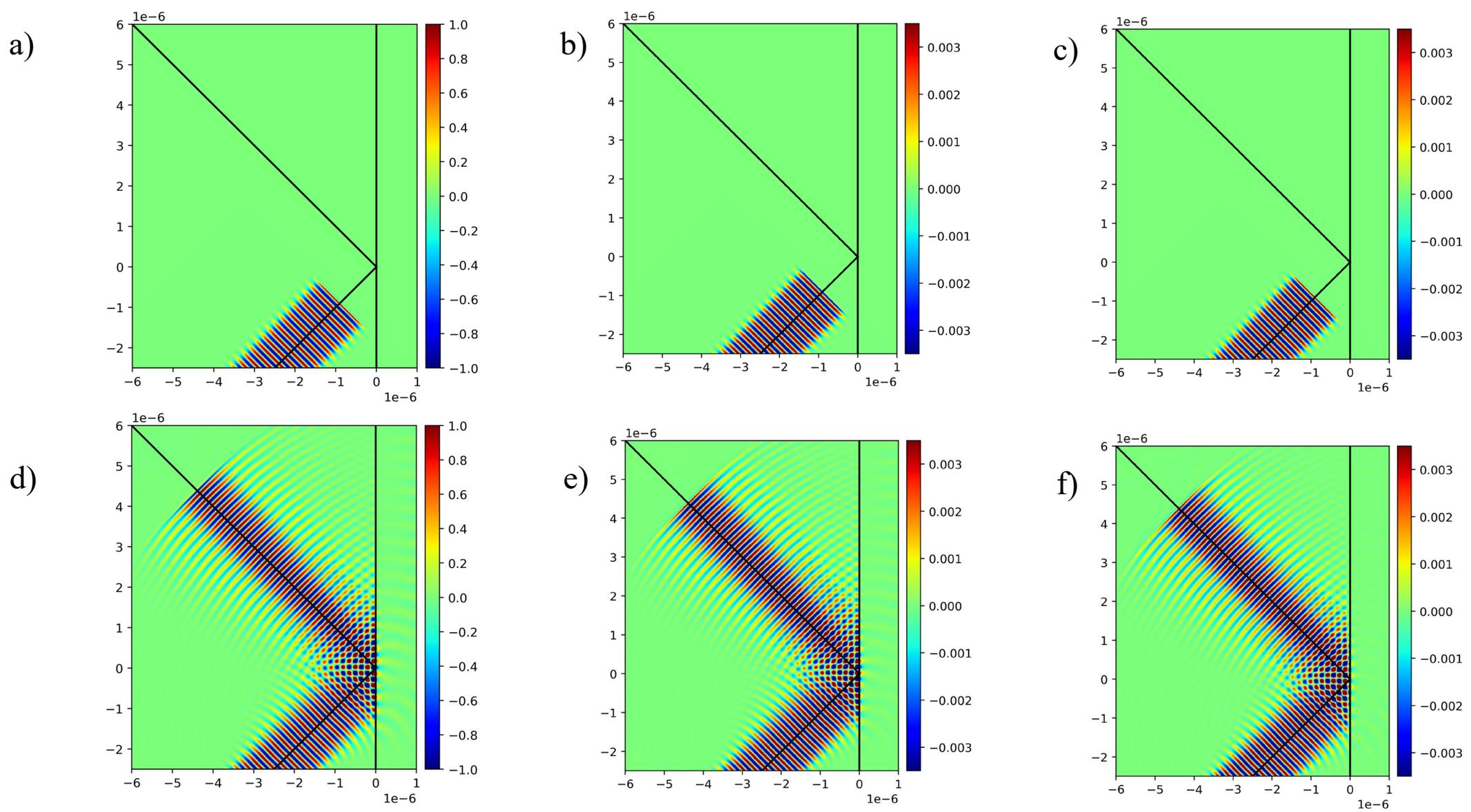

*Fig. 11) CED: Total internal reflection of a compact electromagnetic beam by a dielectric slab. Figs. a), b), and c) show $B_z$, $D_x$, and $D_y$ at the initial time. Figs. d), e), and f) show the same at a final time of $5 \times 10^{-14}$s. The surface of the dielectric slab is identified by the vertical black line. The oblique black lines demarcate the angle of incidence and the angle of total internal reflection.*

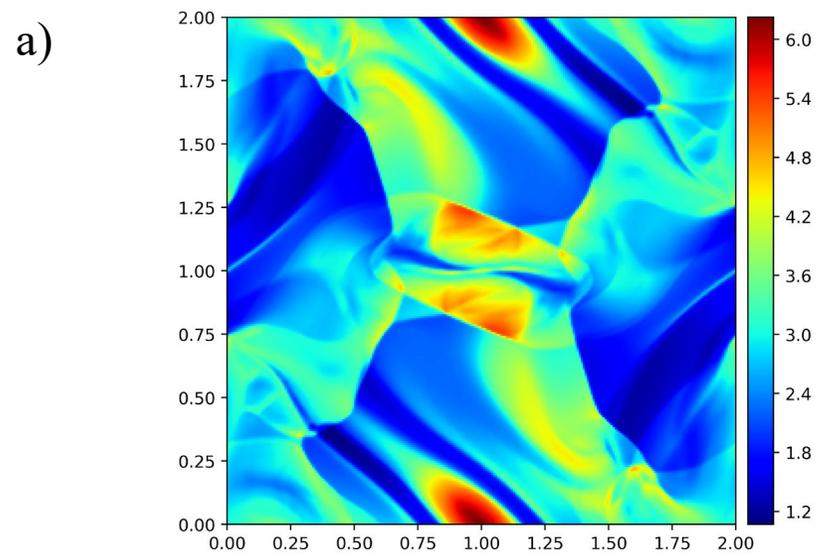 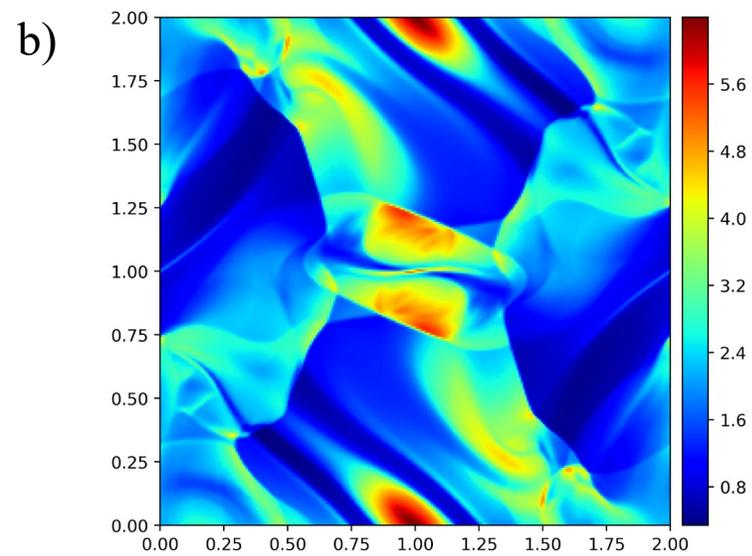
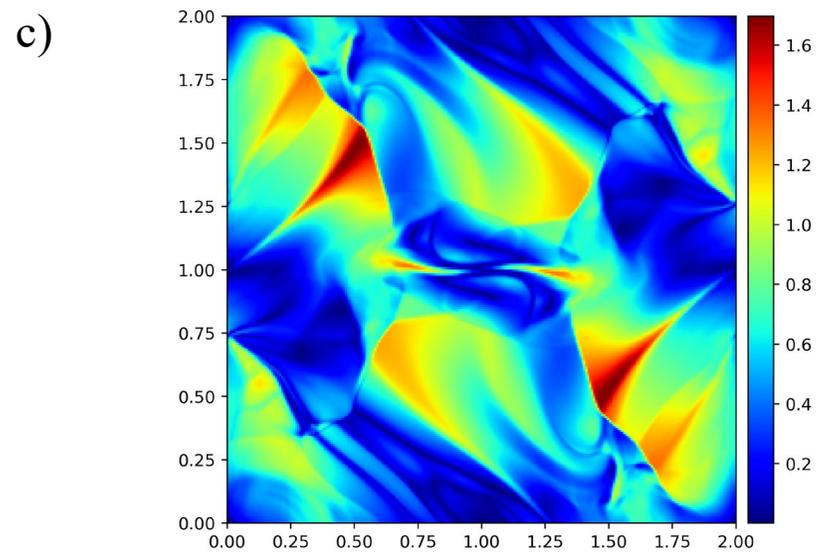 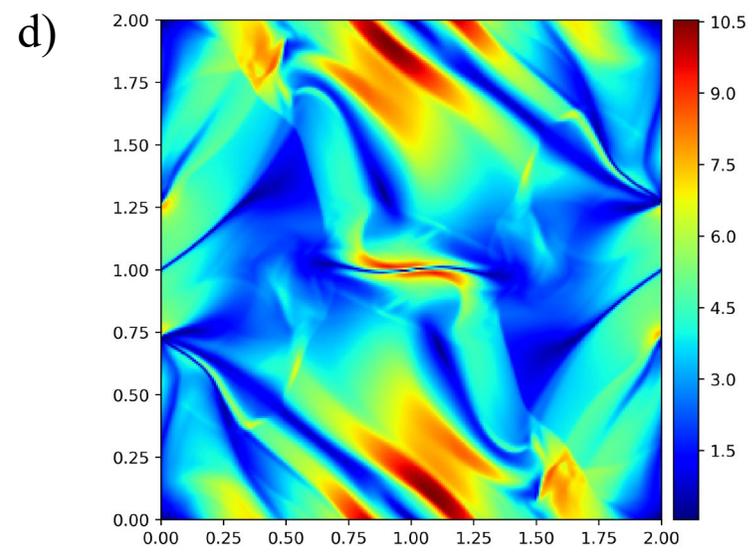

*Fig. 12) MHD: Orszag-Tang Problem. Fig. 12a shows the density, Fig. 12b shows the pressure, Fig. 12c shows magnitude of the velocity and Fig. 12d shows magnitude of the magnetic field vector using the ninth-order accurate scheme.*

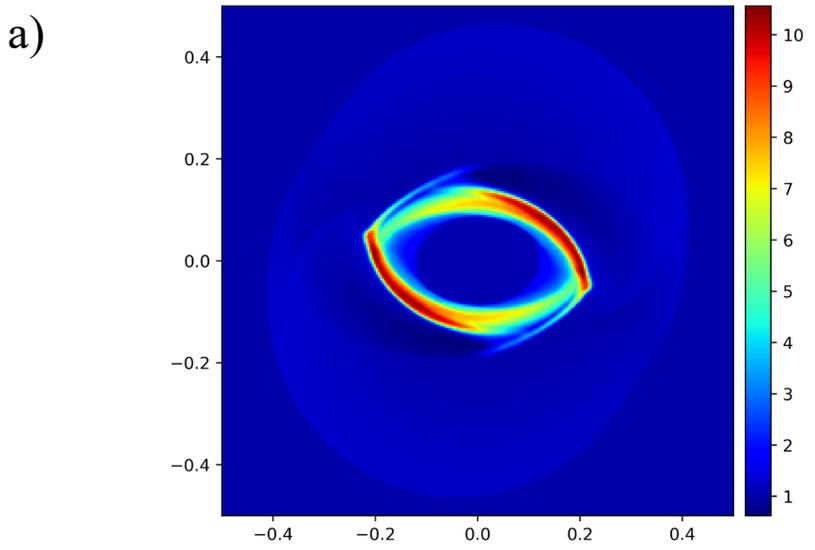 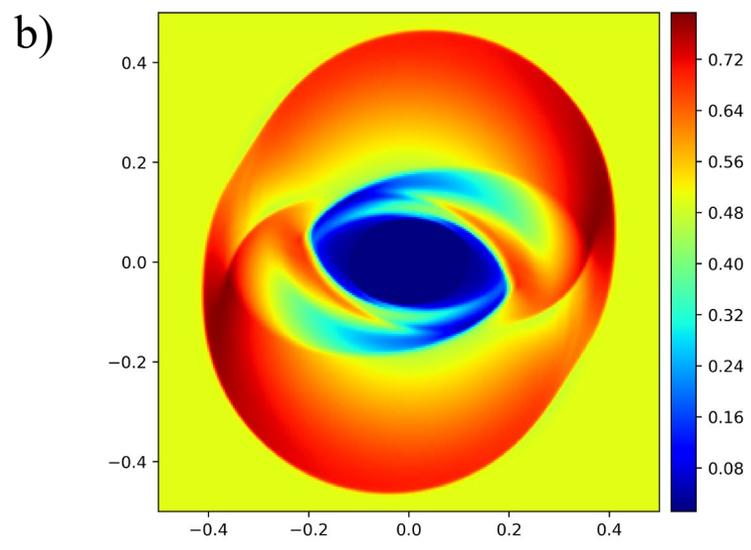
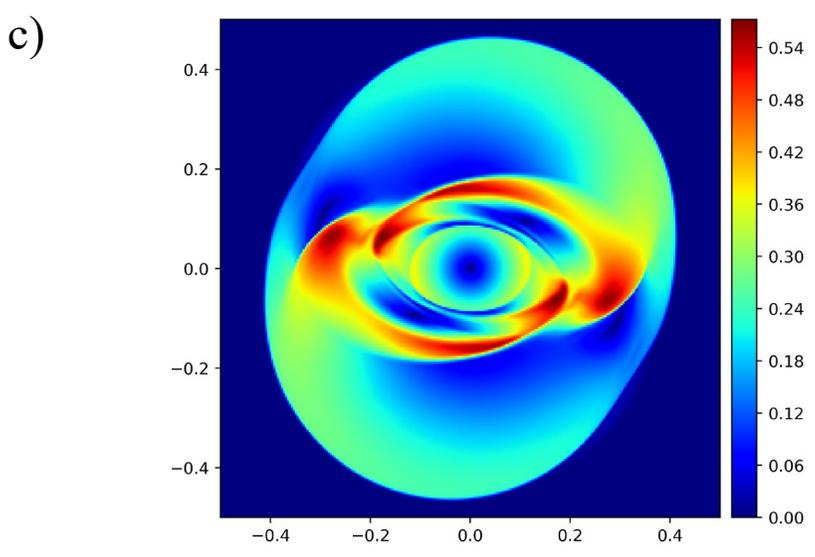 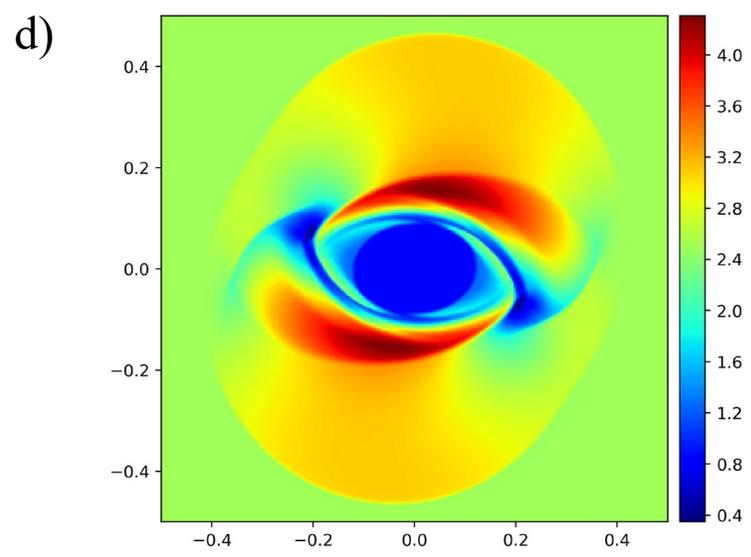

*Fig. 13) MHD: Rotor Problem. Fig. 13a shows the density, Fig. 13b shows the pressure, Fig. 13c shows magnitude of the velocity and Fig. 13d shows magnitude of the magnetic field vector using the seventh-order accurate scheme.*

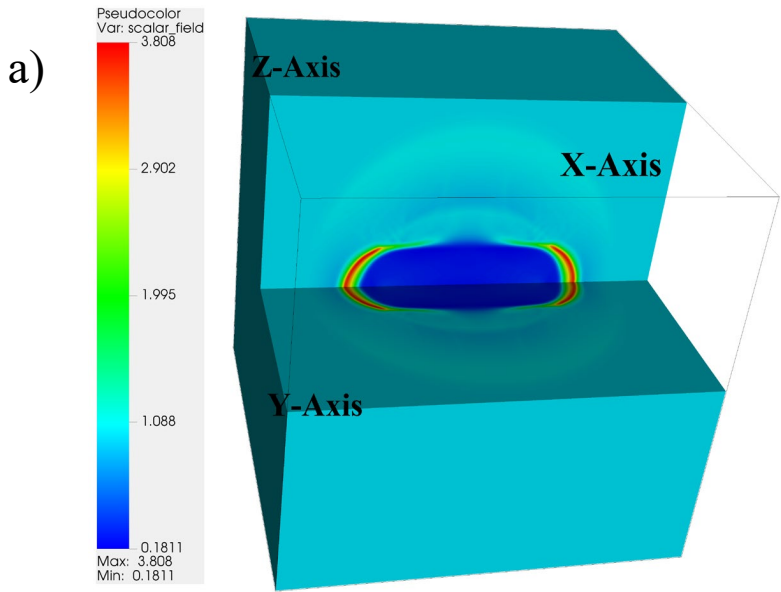 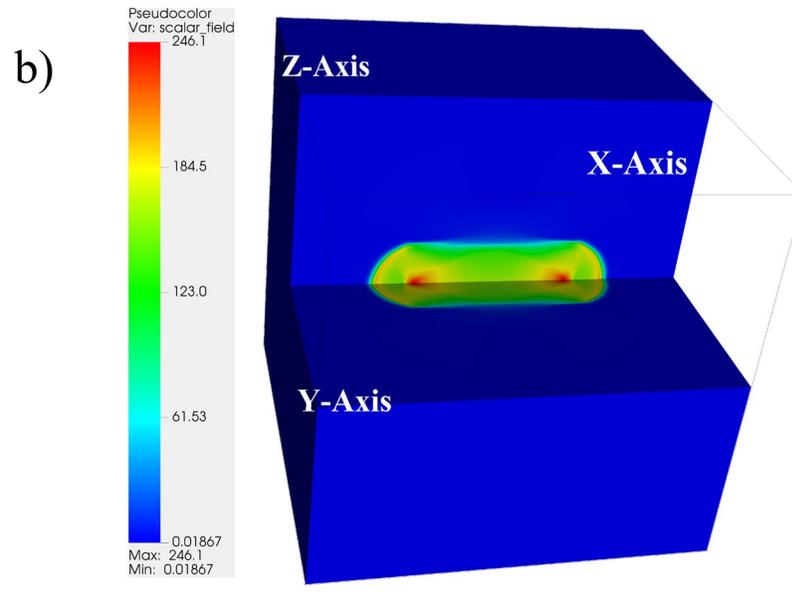
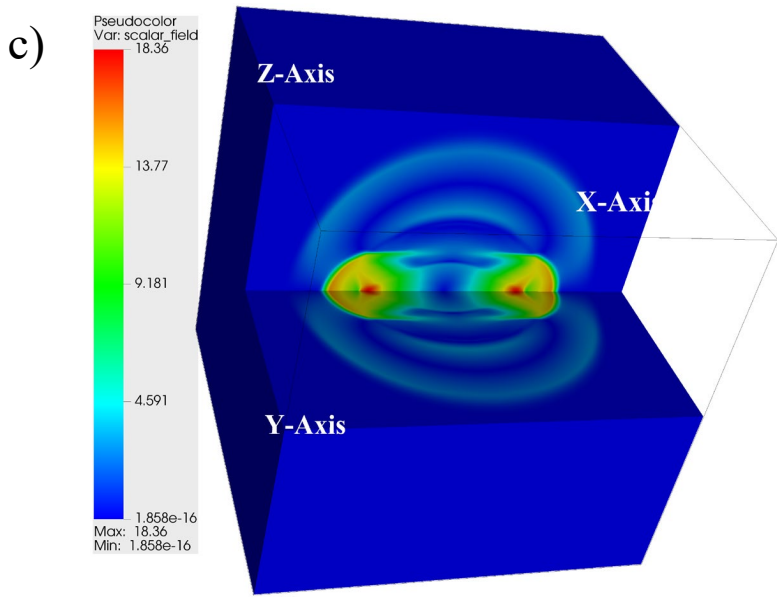 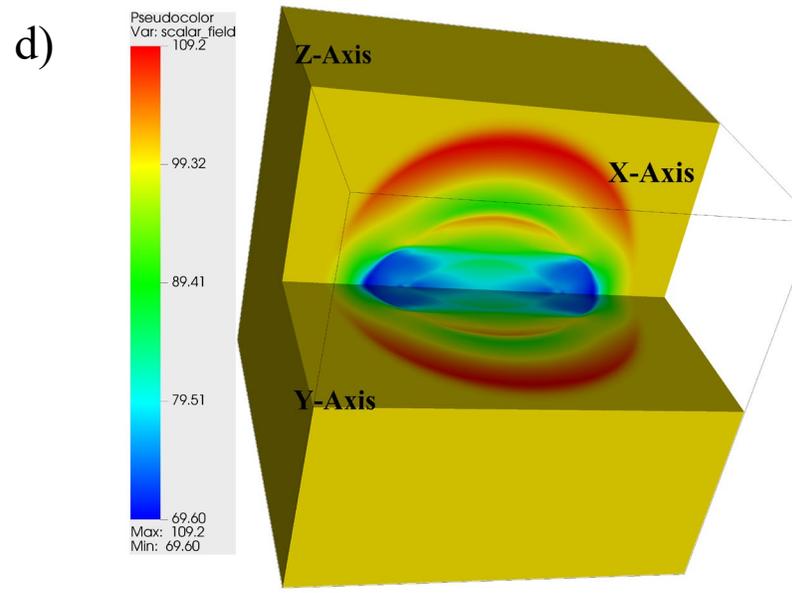

*Fig. 14) MHD: Blast Problem-I (BLAST-I). Fig. 14a shows the density, Fig. 14b shows the pressure, Fig. 14c shows magnitude of the velocity and Fig. 14d shows magnitude of the magnetic field vector using the fifth-order accurate scheme. Here we show the three-dimensional variant of this problem.*

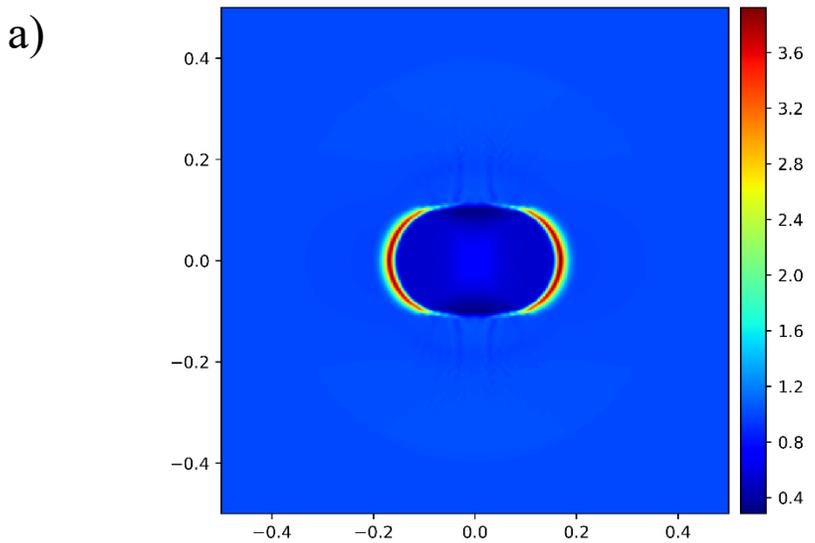 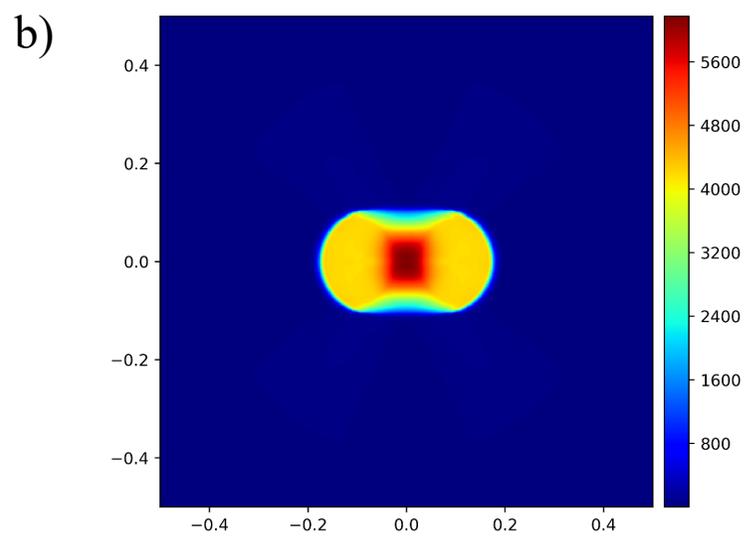
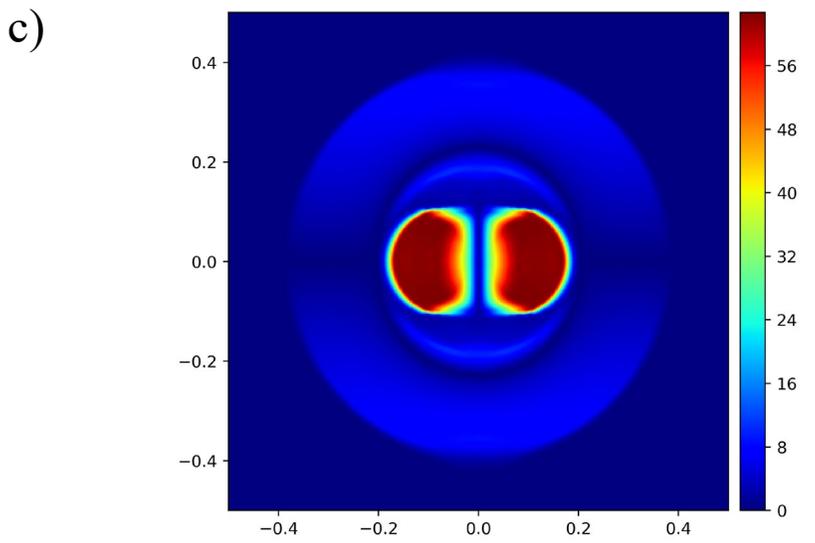 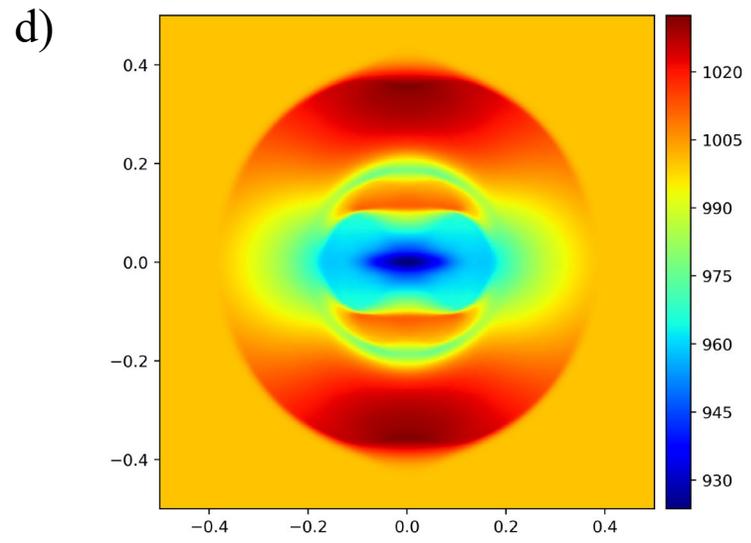

*Fig. 15) MHD: Blast Problem-II (BLAST-II). Fig. 15a shows the density, Fig. 15b shows the pressure, Fig. 15c shows magnitude of the velocity and Fig. 15d shows magnitude of the magnetic field vector using the fifth-order accurate scheme.*

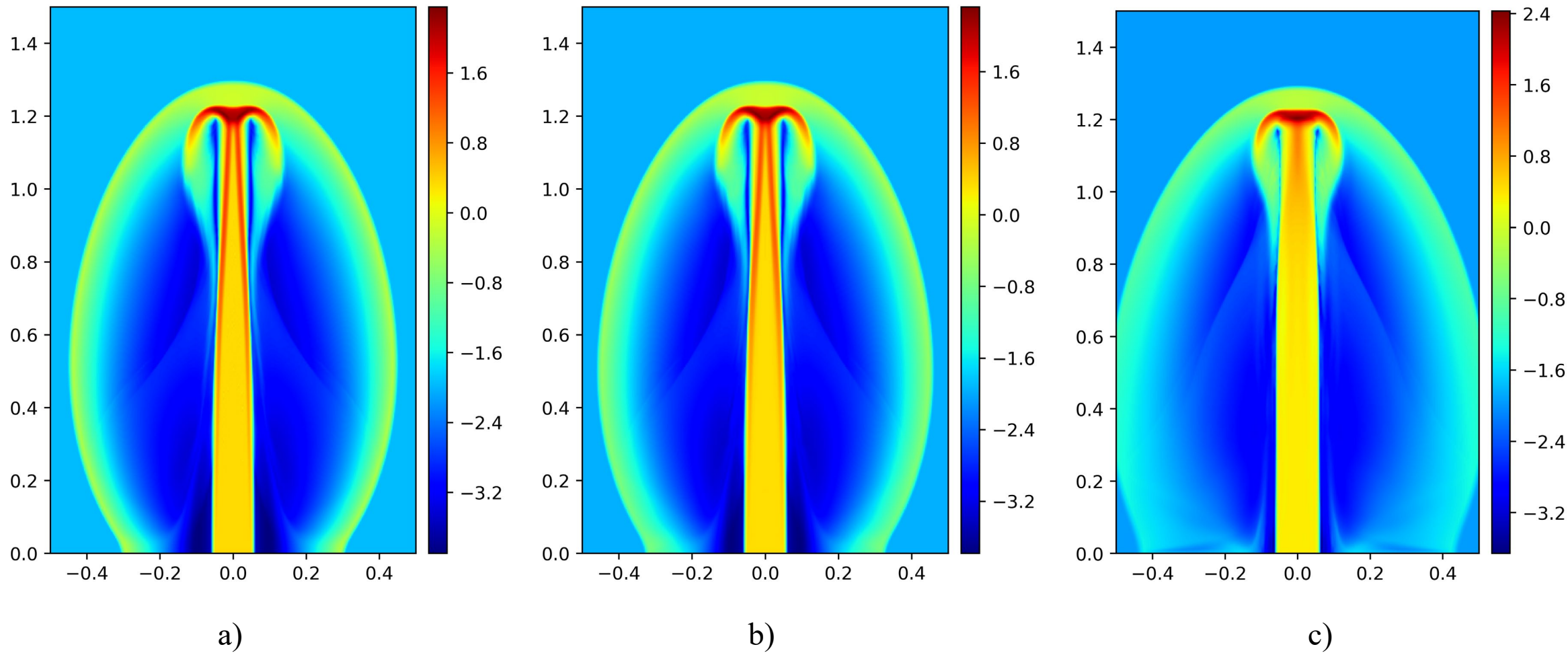

*Fig. 16) MHD: Astrophysical Jet Problem (Jet-I,II and III). Figs. 16a,b,c show the resulting density profiles on logarithmic scales for the Jet-I, Jet-II and Jet-III problems, respectively. Seventh-order accurate schemes is used.*

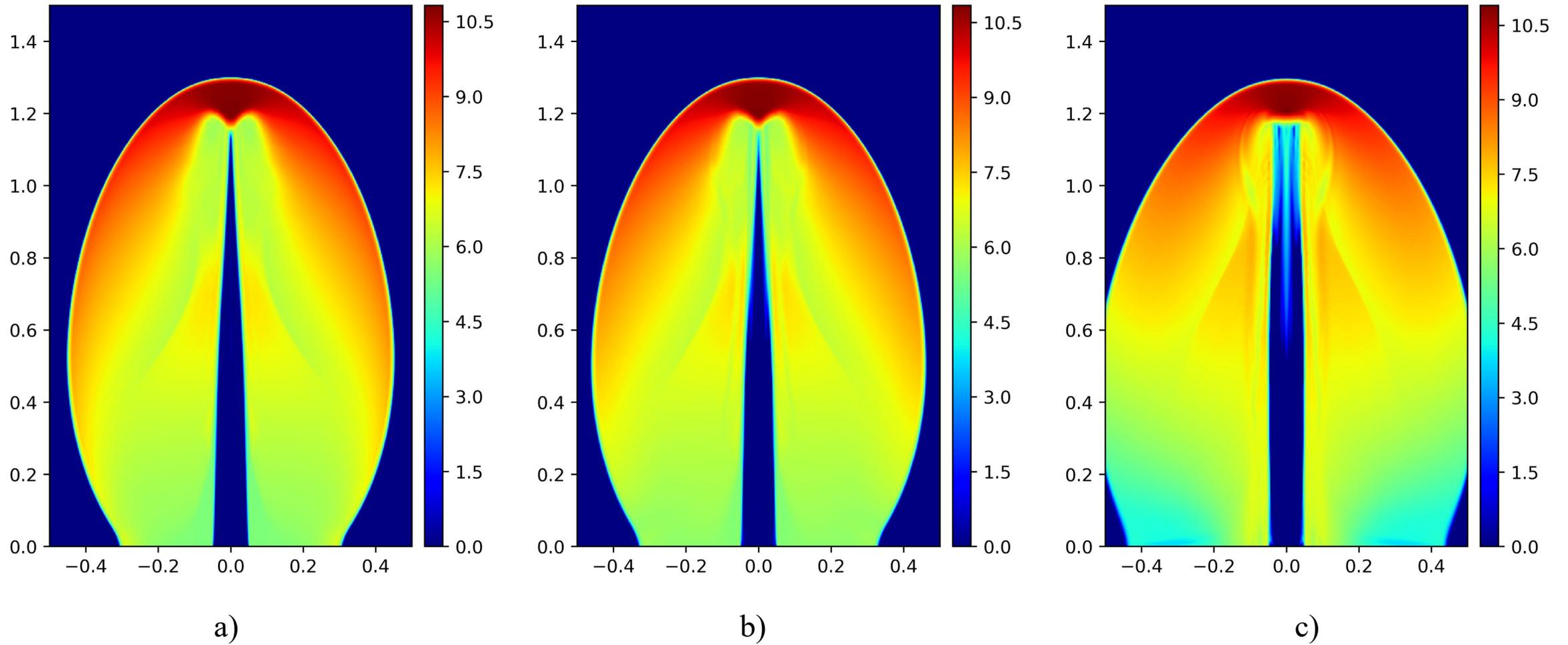

*Fig. 17) MHD: Astrophysical Jet Problem (Jet-I,II and III). Figs. 17a,b,c show the resulting pressure profiles on logarithmic scales for the Jet-I, Jet-II and Jet-III problems, respectively. Seventh-order accurate schemes is used.*

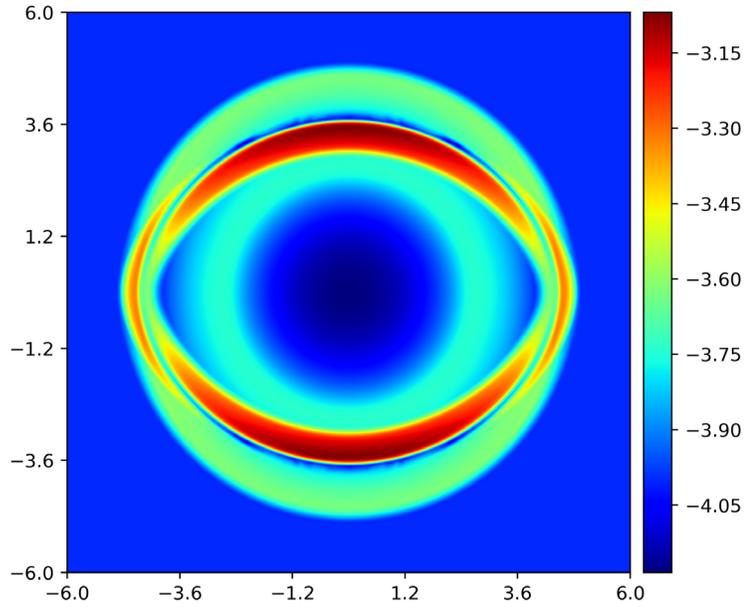 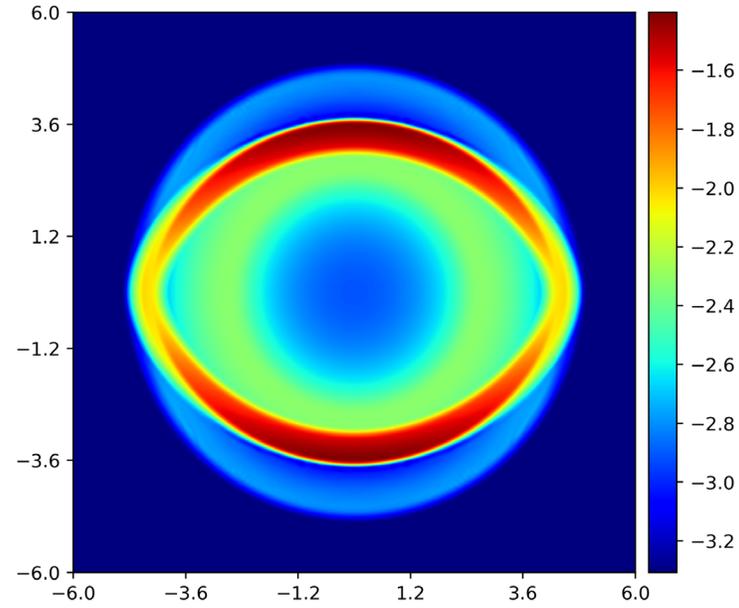 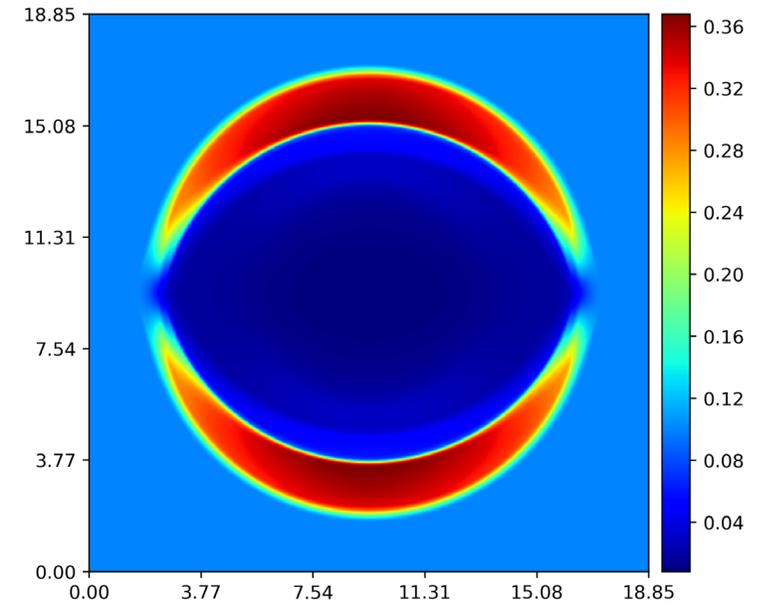

a)  b)  c)

*Fig. 18) RMHD: Blast Problem. Fig. 18a shows the logarithm of density, Fig. 18b shows the logarithm of pressure, and Fig. 18c shows the magnetic pressure using the fifth-order accurate scheme.*

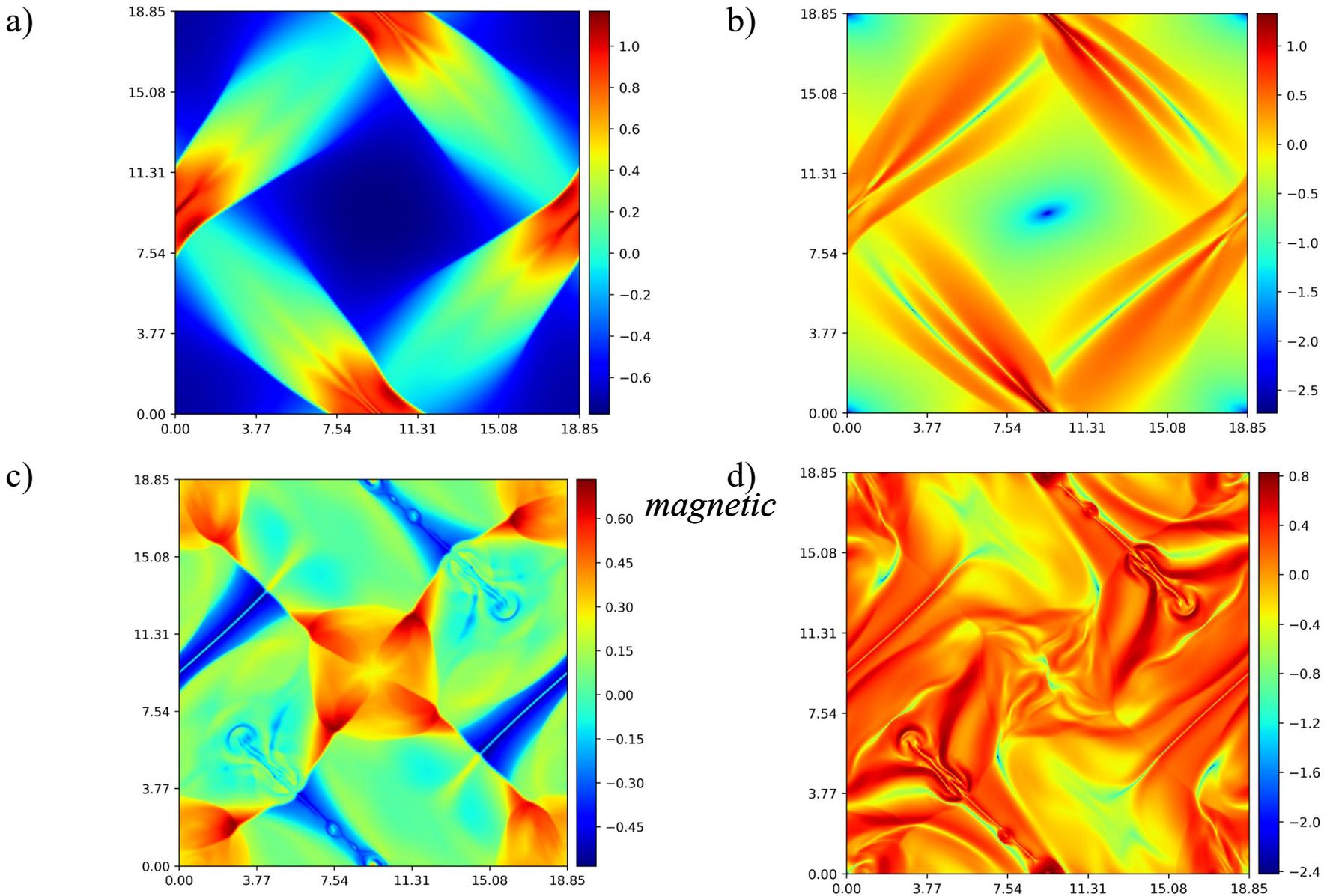

Fig. 19) RMHD: Orszag-Tang Problem. Fig. 19a, b show the logarithm of density and the logarithm of magnetic pressure at time t=2.818127; and Fig. 17c, d show the logarithm of density and the logarithm of magnetic pressure at time t=6.8558. The seventh-order accurate scheme is used.

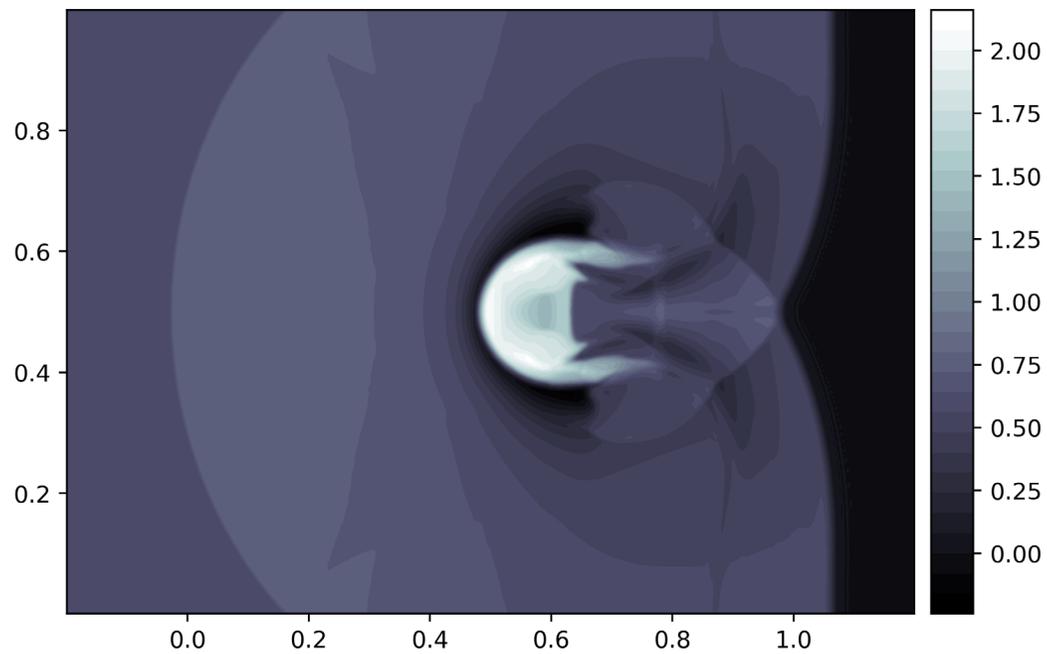 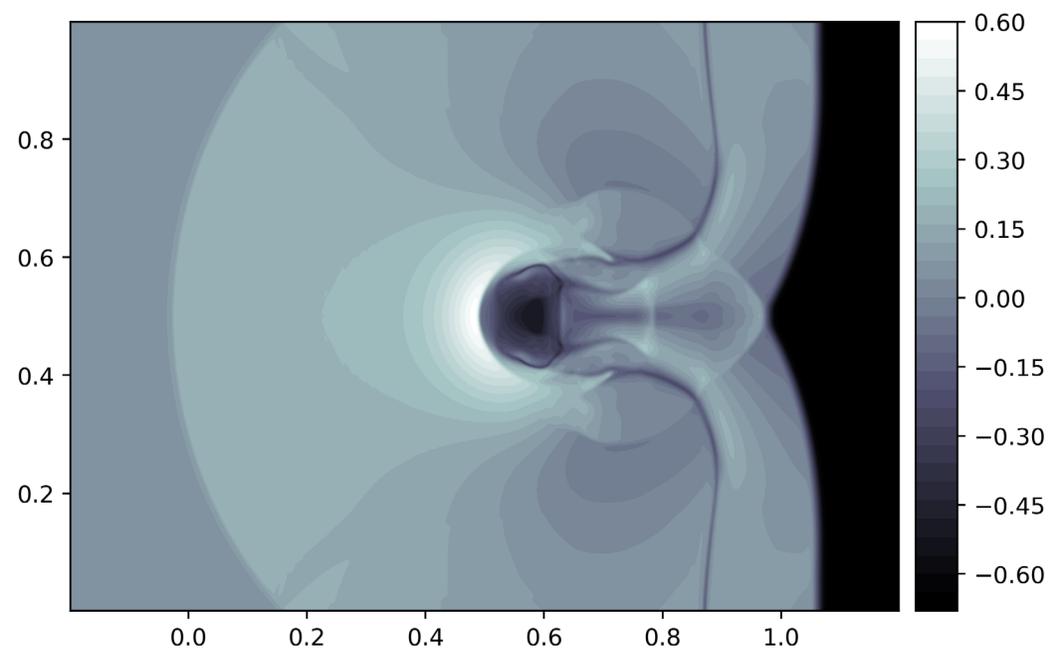

*Fig. 20) RMHD: Shock-Cloud interaction problem. Fig. 20a shows the logarithm of density and Fig. 20b shows the logarithm of magnetic pressure using the ninth-order accurate scheme.*